\documentclass[11pt,a4paper]{amsart}
\usepackage{upref}
\usepackage{enumitem,dsfont,amssymb,bm}
\usepackage{graphicx,color}
\usepackage{mathrsfs}
\usepackage{float}
\usepackage{marginnote}
\usepackage[colorlinks=true, pdftitle={A first hitting time approach}, 
pdfauthor={Daniel Hauer, David Lee}]{hyperref}

\title[A first hitting time approach]{Functional Calculus via the extension
  technique: a first hitting time approach}

\author{Daniel Hauer}
\address[Daniel Hauer]{The University of
  Sydney, School of Mathematics and Statistics, NSW 2006, Australia}
\email{\href{mailto:daniel.hauer@sydney.edu.au}{\nolinkurl{daniel.hauer@sydney.edu.au}}}
\author{David Lee}
 \address[David Lee]{Sorbonne Universit\'e, Laboratoire de
   Probabilit\'es Statistique et Mod\'elisation, 4 Pl. Jussieu, 75005 Paris,France}
\email{\href{mailto:david.lee@upmc.fr}{\nolinkurl{david.lee@upmc.fr}}}

\makeatletter

\newcommand{\Rmnum}[1]{\expandafter\@slowromancap\romannumeral #1@}
\makeatother

\makeatletter
\@namedef{subjclassname@2010}{$2020$ Mathematics Subject Classification}
\makeatother

\thanks{This project has received funding from the
European Union's Horizon 2020 research and innovation program under the
Marie Sk\l{}odowska-Curie grant agreement No 754362 and by the Australian
Research Council via the Discovery Project No DP200101065. In particular, the second
author warmly thanks his PhD supervisor Lorenzo Zambotti for the fruitful
discussions during the development of this paper.
}

\subjclass[2010]{60H30, 47A60, 60B15, 47D07}

\keywords{Complete Bernstein functions, Dirichlet-to-Neumann,
  Dirichlet-to-Robin, Dirichlet-to-Wentzell,
fractional powers, stochastic process, semigroups, Phillips subordination}

\numberwithin{equation}{section}

\theoremstyle{theorem}
\newtheorem{theorem}{Theorem}[section]
\newtheorem{proposition}[theorem]{Proposition}
\newtheorem{lemma}[theorem]{Lemma}
\newtheorem{corollary}[theorem]{Corollary}

\theoremstyle{definition}
\newtheorem{definition}[theorem]{Definition}

\newtheorem{example}[theorem]{Example}

\theoremstyle{remark}
\newtheorem{remark}[theorem]{Remark}
\newtheorem{notation}[theorem]{Notation}

\theoremstyle{plain}

\newcommand\R{{\mathbb{R}}}
\newcommand\N{\mathbb{N}}
\newcommand\C{{\mathbb{C}}}
\newcommand\E{\mathcal{E}}

\newcommand\X{{\mathds{X}}}

\newcommand\td{\mathrm{d}}
\newcommand\dx{\mathrm{d}x }
\newcommand\du{\mathrm{d}u }
\newcommand\dm{\mathrm{d}m }
\newcommand\dy{\mathrm{d}y }
\newcommand\dr{\mathrm{d}r }

\newcommand\dz{\mathrm{d}z }
\newcommand\ds{\mathrm{d}s }
\newcommand\dt{\mathrm{d}t }

\newcommand\CBF{\mathcal{C}\mathcal{B}\mathcal{F}}

\DeclareMathOperator{\Rg}{Rg}

\DeclareMathOperator{\Capacity}{Cap}

\DeclareMathOperator{\supp}{supp}

\def\1{\raisebox{2pt}{\rm{$\chi$}}}

\newcommand\abs[1]{\lvert#1\rvert}

\newcommand\norm[1]{\lVert#1\rVert}
\newcommand\lnorm[1]{\left\lVert#1\right\rVert}

\definecolor{darkred}{rgb}{0.7,0.1,0.1}

\newcounter{fpcounter}
\renewcommand{\thefpcounter}{\arabic{fpcounter}}

\newenvironment{fp}[2]{%
\refstepcounter{fpcounter}%
\label{#1}%
\noindent\textbf{Problem~\thefpcounter.}%
}%
{}%

\newtheorem*{op}{Open Problem}

\begin{document}
\begin{abstract}
 In this article, we present a solution to the problem:
  \begin{center}
    \emph{Which type of linear operators can be realized by the
      Dirichlet-to-Neumann operator associated with the operator
      $-\Delta-a(z)\frac{\partial^{2}}{\partial z^2}$ on an extension
      problem?}
  \end{center}
  which was raised in the pioneering work~[Comm. Par.Diff. Equ. 32
  (2007)] by Caffarelli and Silvestre. In fact, we even go a step
  further by replacing the negative Laplace operator $-\Delta$ on
  $\R^{d}$ by an $m$-accretive operator $A$ on a general Banach space
  $X$ and the Dirichlet-to-Neumann operator by the
  Dirichlet-to-Wentzell operator. We establish uniqueness of solutions
  to the extension problem in this general framework,
  which seems to be new in the literature and of independent
  interest. Our aim of this paper is to provide a new Phillips-Bochner
  type functional calculus which uses probabilistic tools from
  excursion theory. With our method, we are able to characterize all
  linear operators $\psi(A)$ among the class $\CBF$ of complete Bernstein
  functions $\psi$, resulting in a new characterization of the famous
  \emph{Phillips' subordination theorem} within this class $\CBF$.
\end{abstract}

\maketitle

\tableofcontents

\section{Introduction and main results}

\subsection{Introduction}
In~\cite{MR2354493}, Caffarelli and Silvestre
showed that the fractional Laplace operator
$\psi(-\Delta)=(-\Delta)^{\sigma}$ on $\R^{d}$, $d\ge 1$,
$\psi(s)=s^{\sigma}$, $0<\sigma<1$, defined as the singular integral
operator
\begin{displaymath}
  (-\Delta)^{\sigma}f(x)=C_{\sigma,d}\,\textrm{C.V.}
  \int_{\R^{d}}\frac{f(x)-f(\xi)}{\abs{x-\xi}^{d+2\sigma}}\td\xi,
\end{displaymath}
$f\in C^{\infty}_{c}(\R^{d})$, 
 can be characterized by the 
 \emph{Dirichlet-to-Neumann operator}
\begin{equation}
  \label{eq:94}
  f_{\vert \R^{d}}\mapsto \Lambda_{\sigma}f:=\lim_{y\rightarrow
    0^+}-y^{1-2\sigma}u_{y}(.,y)_{\vert \R^{d}}
\end{equation}
corresponding to the \emph{incomplete Dirichlet/extension problem}
 \begin{equation}
   \label{eq:92}
   \begin{cases}
     -\Delta u-\frac{1-2\sigma}{y}u_{y}-u_{yy}=0 & \text{on
       $\R^{d+1}_{+}$,}\\
     \hspace{3.3cm}u=f & \text{on $\partial\R^{d+1}_{+}=\R^{d}$,}
   \end{cases}
 \end{equation}
 for the negative \emph{Bessel operator}
 $\mathcal{A}:=-\left(\Delta+\frac{1-2\sigma}{y}\frac{\partial}{\partial
   y}+\frac{\partial^2}{\partial y^2}\right)$ on the half space
 $\R^{d+1}_{+}:=\R^{d}\times (0,\infty)$. 
 In particular, it was observed that by applying
 the change of variable $z=(y/2\sigma)^{2\sigma}$ to the solution
 $u(y):=u(\cdot,y)$ of the \emph{Bessel equation}
 \begin{equation}
   \label{eq:138}
  -\Delta u-\frac{1-2\sigma}{y}u_{y}-u_{yy}=0\qquad\text{in $\R^{d+1}_{+}$,}
\end{equation}
then $u(z)$ becomes a solution of the elliptic \emph{extension equation}
\begin{equation}
  \label{eq:93}
  -\Delta u-a(z)u_{zz}=0\qquad\text{in $\R^{d+1}_{+}$,}
\end{equation}
for the coefficient $a(z)=z^{-\frac{1-2\sigma}{\sigma}}$. Moreover, one has that
\begin{displaymath}
  y^{1-2\sigma}u_{y}(y)=(2\sigma)^{1-2\sigma}u_{z}(z)
\end{displaymath}
and hence, the
Dirichlet-to-Neumann operator $\Lambda_{\sigma}$ given by
\eqref{eq:94} reduces to
\begin{displaymath}
  \Lambda_{\sigma}f=-(2\sigma)^{1-2\sigma}u_{z}(0).
\end{displaymath}
Subsequently, the following problem was posed
in~\cite[Section~7.]{MR2354493}:\medskip

\begin{fp}{fp:A}{}
\emph{What other linear operators $\psi(-\Delta)$ can be obtained from the
    Dirich\-let-to-Neumann operator associated with the operator $-\Delta
  -a(z)\frac{\partial^{2}}{\partial z^{2}}$?}
\end{fp}
\mbox{}\bigskip 

The aim of this paper, is to provide an answer to this problem, and
even to go beyond. We replace the negative Laplace operator $-\Delta$
(with vanishing conditions at infinity) on $\R^{d}$ by a general
closed linear $m$-accretive operator $A$ defined on a Banach space
$X$. Then $A$ admits the property that
$-A$ generates a $C_{0}$-semigroup $\{e^{-t A}\}_{t\ge 0}$ of
contractions (see Section~\ref{subsec:semigroups} for further
details). Further, if for a given \emph{string $m$ on $\R$ of
infinite length} (see Definition~\ref{def:string}), $\psi_{m}$ is the
associated \emph{complete Bernstein function} with L\'evy triple
$(0,m(0+),\nu_{m})$ (see Definition~\ref{def:1} and~\eqref{eq:48ast}), and
$\psi_{m}(A)$ is the operator defined by~\eqref{eq:109} via the L\'evy
triple $(0,m(0+),\nu_{m})$, then the main result of this article (Theorem~\ref{thm:1}) states
that the characterization
 \begin{equation}
   \label{eq:110}
   \psi_{m}(A)=\Lambda_m
 \end{equation}
 holds, where $\Lambda_{m}$ is the \emph{Dirichlet-to-Wentzell
   operator} given by
 \begin{equation}
   \label{eq:91}
  \Lambda_{m}f:=m(0+)Au(0)-\frac{1}{2}\frac{\td u}{\dz}_{\! +}\!(0)
 \end{equation}
 for every $f\in D(A)$ with corresponding unique bounded \emph{weak solution}
 $u$ of the \emph{(incomplete) Dirichlet problem} (see Definition~\ref{def:weak-solution Dirichlet-problem} below)
 \begin{equation}\label{eq:16}
    \begin{cases}
      \hspace{0,25cm}\mathcal{A}_{m}u(z)=0 &\qquad \text{for $z\in
        (0,\infty)$,}\\
      \hspace{0,25cm}\phantom{\mathcal{A}_{m}}u(0)=f,& 
    \end{cases}
 \end{equation}
 for the \emph{extension operator} 
 \begin{equation}
   \label{eq:137}
 \mathcal{A}_{m}:=A+B_{m}
\end{equation}
on the \emph{extended space} $\mathcal{X}_{+}=X\times (0,\infty)$. We
also refer to Dirichlet problem~\eqref{eq:16} as the \emph{extension
  problem} since the operator $\mathcal{A}_{m}$ extends the operator
$A$ acting on $X$ by the $2$nd-order differential operator
 \begin{displaymath}
   B_{m}:=-\frac{1}{2}\frac{\td}{\dm}\frac{\td
   }{\dz}\qquad\text{acting on $(0,\infty)$,}
 \end{displaymath}
 where we denote by
\begin{equation}
  \label{eq:107}
  \frac{\td u}{\td m}(z):=\lim_{h\to 0}\frac{u(z+h)-u(z)}{m(z+h)-m(z)}
\end{equation}
(if the limit exits in $X$) the \emph{$m$-derivative} at
$z\in (0,\infty)$ of a function $u : [0,\infty)\to X$. We refer to
Definition~\ref{def:weak-dfdm-intro} below for a formulation of the
this notion in the \emph{sense of distributions} and also to
Appendix~\ref{sec:def-weak} for further discussions.

For the proof of the characterization~\eqref{eq:110}, we employ an
intermezzo of probabilistic tools and functional analysis. We begin by
proving existence and uniqueness of bounded weak solutions of the
Dirichlet problem~\eqref{eq:16}. Motivated from the interpretation
that problem~\eqref{eq:16} can be considered as a classic elliptic
Dirichlet problem on the ``extended region'' $\mathcal{X}_{+}$, our
existence proof relies on the classical approach of stopping a related
stochastic process $\{(X_{t},Z_{t}\}_{t\ge 0}$ in $\mathcal{X}_{+}$ at $z=0$. To be more precise,
for a given string $m$ on $\R$, the operator $-B_{m}$ generates a
stochastic process $\{Z_{t}\}_{t\ge 0}$ on $[0,\infty)$, which we call
\emph{generalized diffusion associated with $m$} (see
Definition~\ref{def:gen-diffusion}). Thus, one can define the \emph{first
  hitting time} $\tau$ of zero by $\{Z_{t}\}_{t\ge 0}$ starting from
$z\in (0,\infty)$. If one assumes that the Banach space $X=X(\Sigma)$ is a function
space with domain $\Sigma$ and $-A$ generates a Markov process
$\{X_{t}\}_{t\ge 0}$ in the state space $\Sigma$, then it is clear
that for every $f\in X(\Sigma)$, the unique solution $u$ of Dirichlet
problem~\eqref{eq:16} is given by
\begin{equation}
  \label{eq:111}
  u(x,z)=\mathds{E}_{(x,z)}\left(f\big(X_{\tau}\big)\right)
\end{equation}
for every $(x,z)\in \Sigma\times [0,\infty)$. Since the first hitting time
$\tau$ admits a \emph{density} (see
Section~\ref{SEC:PRIMER_STRINGS} for details) 
\begin{displaymath}
  \omega_{\tau}(t,z):=\frac{\mathbb{P}_{z}(\tau \in \dt)}{\dt}
\end{displaymath}
and since the process $\{X_{t}\}_{t\ge 0}$ in $\Sigma$ is one-to-one related to
the transition semigroup $\{e^{-tA}\}_{t\ge 0}$ on $X(\Sigma)$. Thus the
representation~\eqref{eq:111} of a weak solution $u$ of~\eqref{eq:16}
is equivalent to the \emph{Poisson formula}
 \begin{equation}
   \label{eq:17}
   u(z)=\int_{0}^{\infty}(e^{-tA}f)\,\omega_{\tau}(t,z)\,\dt 
 \end{equation}
 for every $z\in [0,\infty)$. One crucial advantage of formula~\eqref{eq:17} is that the process
 $\{X_{t}\}_{t\ge 0}$ is not directly involved anymore to describe a
 weak solution $u$ of~\eqref{eq:16}. Hence, formula~\eqref{eq:17} provides a
 strong candidate for being a solution of Dirichlet problem~\eqref{eq:16}
 even though $-A$ may not necessarily generate a stochastic process, but
 still generate a $C_{0}$-semigroup of contractions on $X$. By
 using Kent's theorem~\cite{MR637387} on the spectral decomposition of
 the first hitting time density, we are able to derive a link between
 the first hitting time density $\omega_{\tau}$ in~\eqref{eq:17} and
 the L\'evy measure density of the inverse local time
 $\{\tilde{L}_{t}^{-1}\}_{t\ge 0}$ at zero of the generalized
 diffusion $\{Z_{t}\}_{t\ge 0}$. This connections allow us
to easily calculate the Dirichlet-to-Wentzell map $\Lambda_{m}$ given
by~\eqref{eq:91}, to show that the characterization~\eqref{eq:110}
holds, and to establish a Phillips-Bochner type functional calculus.

\subsection{Main results}\label{subsec:main-results}






In order to state the main theorems of this paper, we need first to
introduce some classic definitions and helpful notations. Throughout
this section, $X$ denotes a Banach space and $A$ a closed, linear
operator on $X$ with dense domain $D(A)$, and $\{e^{-t A}\}_{t\ge 0}$
the $C_{0}$-semigroup of contractions on $X$ generated by $-A$.

\begin{definition}\label{def:string}
 A non-decreasing, right-continuous function
 $m: \R\rightarrow [0,\infty)$ is called a \emph{string on
 $\R$ of infinite length} provided $m$ has the following properties
 \begin{enumerate}[label=(\roman*)]
  \item $m(x)=0$ for all $x<0$ and $m(-0)=0$, 
  \item $m(x_0)<\infty$ for some $x_0\ge 0$, and
  \item $m(x)>0$ for all $x>0$.
  \end{enumerate}
  We denote the family of all strings of infinite length by $\mathfrak{m}_{\infty}$.
\end{definition}

In the following, if nothing else is said, we always refer to $m$ as
a string on $\R$ of infinite length. The case of strings of
\emph{finite length} shall be studied in a forthcoming work.  Then, there is a unique Radon measure
$\mu_{m} : \mathcal{B}(\R)\to [0,+\infty]$ on the Borel
$\sigma$-algebra $\mathcal{B}(\R)$ of $\R$ such that
\begin{equation}
  \label{eq:56}
 \mu_{m}((a,b])=m(b+)-m(a+)\qquad\text{for every $a$, $b\in \R$ with $a<b$}
\end{equation}
(cf,~\cite[Theorem~6.7]{MR3726909} and see also~\cite[Remark~6.11]{MR3726909}).
Since $m(x)=0$ for all $x<0$, \eqref{eq:56} yields that the
$\mu_{m}((a,b])=0$ for all $a<b<0$. Further, we denote by
$E_{m}=\supp(\mu_{m})$ the \emph{support} of $\mu_{m}$.


 For a better understanding of the functional calculus, we develop
 here, but also, in order to illustrate that our results generalize
 previous ones obtained for the fractional power case (cf , for example, \cite{MR2354493,
   MR2754080,MR3056307,MR3772192}), we provide the following example.

 \begin{example}\label{ex:1}
   Let $A=-\Delta$ be the negative Laplace operator on 
   $X=L^{2}(\R^{d})$ equipped with vanishing conditions at
   infinity. Then, in the case of the
   fractional power $\psi_{m}(-\Delta)=(-\Delta)^{\sigma}$,
   $\sigma\in (0,1)$ in $X=L^{2}(\R^{d})$, 
   in the extension equation~\eqref{eq:93} the coefficient $a(z)$ is
   given by
   \begin{equation}
     \label{eq:8}
     a_{\sigma}(z)=\tfrac{1}{2m_{\sigma}'(z)}\qquad\text{ and }\qquad
   a_{\sigma}(z)u_{zz}=\frac{1}{2}\frac{\td}{\dm_{\sigma}}\frac{\td u}{\dz}(z),\;
   z\in (0,\infty),
  \end{equation}
  for the string $m_{\sigma}\in \mathfrak{m}_{\infty}$ on $\R$
   given by
   \begin{equation}
     \label{eq:113}
      m_{\sigma}(z):=
      \begin{cases}
        \frac{1}{2}\frac{\sigma}{1-\sigma}
        z^{\frac{1-\sigma}{\sigma}} &\quad\text{if $z>0$,}\\
        0&\quad\text{if $z\le 0$,}
      \end{cases}
 \end{equation}
 for every $z\in \R$. For convenience, we write in~\eqref{eq:8} $m'_{\sigma}(z)$
 to denote $\frac{\td m_{\sigma}}{\dz}$.Thus, equation~\eqref{eq:93} can be rewritten in
 $X=L^{2}(\R^{d})$ by
 \begin{equation}
  \label{eq:108}
  \mathcal{A}_{m_{\sigma}}u(z)=0\qquad\text{ 
    for $z>0$.}
 \end{equation}
 It is worth mentioning that the generalized diffusion
 $\{Z_{t}\}_{t\ge 0}$ generated by $-B_{m_{\sigma}}$, in this case, coincides
 up to a multiple constant with the scaled Bessel process
 $\{Y_{t}^{2\sigma}\}_{t\ge 0}$. We continue discussing the fractional
 power case $A^\sigma$ in the Sections~\ref{subsec:represenation-formula} \&
 \ref{subsec:DtN-map}, and in Example~\ref{ex-scale-of-Bessel-process}
 of Section~\ref{subsec:speed-measures}. For the convenience of the
 reader, we provide in Appendix~\ref{subsec:bessel-processes} of this paper a brief
review of the Bessel process and related properties, which are relevant here.
 \end{example}

 For $1\le q< \infty$, let $L^{q}_{\mu_{m}}(0,\infty;X)$ (respectively,
 $L^{q}_{loc,\mu_{m}}((0,\infty);X)$) denote the weighted Lebesgue
 spaces of all $\mu_{m}$-a.e. rest-classes of measurable functions
 $u : [0,\infty)\to X$ with finite integral
 $\int_{[0,\infty)}\norm{u}^q_{X}\td\mu_{m}$ (respectively, finite
 $\int_{K}\norm{u}^q_{X}\td\mu_{m}$ for every compact subset
 $K\subseteq (0,\infty)$).

\begin{definition}
  \label{def:weak-dfdm-intro}
  For a given $u\in L^{1}_{loc}((0,\infty);X)$, one calls a function
  $g\in L^{1}_{loc, \mu_{m}}((0,\infty);X)$ a \emph{weak $m$-derivative
    of $f$} provided $u$ and $g$ satisfy
\begin{equation}
\label{eq:85}
  \int_{0}^{\infty}g(z)\,\xi(z)\,\td\mu_{m}(z)=-\int_{0}^{\infty}
  u(z)\,\frac{\td \xi}{\dz}(z)\,\dz
\end{equation}
for every $\xi\in C^{\infty}_{c}((0,\infty))$. Due to
Lemma~\ref{lem:variational-lemma} (in the appendix of this paper), a
function $g$ is uniquely defined through~\eqref{eq:85}. Thus, we can
call $g$ \emph{the weak $m$-derivative of $u$} and set
$\frac{\td u}{\td m}=g$.
\end{definition}

For $1\le q< \infty$, we write
$W^{1,q}_{loc,\mu_{m}}((0,\infty);X)$ to denote the \emph{mixed
  $1^{\textrm{st}}$-Sobolev space} of all functions
$u\in L^{q}_{loc}((0,\infty);X)$ with weak $m$-derivative
$\frac{\td u}{\td m}\in L^{q}_{loc,\mu_{m}}((0,\infty);X)$ and by
$W^{1,q}_{\mu_{m}}((0,\infty);X)$ the Sobolev space of all functions
$u\in L^{q}(0,\infty;X)$ with
$\frac{\td u}{\td m}\in L^{q}_{\mu_{m}}(0,\infty;X)$. If
$\mu_{m}=\mu_{Leb}$ is the Lebesgue measure on $(0,\infty)$, then we use the standard notation $W^{1,q}_{loc}((0,\infty);X)$
instead of $W^{1,q}_{loc,\mu_{Leb}}((0,\infty);X)$ and $W^{1,q}(0,\infty;X)$ instead
of $W^{1,q}_{\mu_{Leb}}((0,\infty);X)$.

\begin{remark}
  \label{rem:m-derivative}
  On the other hand, for given $u\in L^{1}_{loc}((0,\infty);X)$ and $g\in
L^{1}_{loc, \mu_{m}}((0,\infty);X)$, \eqref{eq:85} means that
the \emph{regular} (vector-valued) distribution $[u]$ given by
\begin{displaymath}
  \langle [u],\xi\rangle:=\int_{0}^{\infty}u(z)\,\xi(z)\,\dz, \quad \xi \in C^{\infty}_{c}((0,\infty)),
\end{displaymath}
has the 
\emph{(vector-valued) measure} $g\,\mu_{m}$ as its \emph{distributional derivative}
\begin{displaymath}
  \langle [u]',\xi\rangle:=-\int_{0}^{\infty}u(z)\,\frac{\td \xi}{\dz}(z)\,\dz=\int_{0}^{\infty} \xi(z)\,
  g(z)\,\td\mu_{m}(z)
\end{displaymath}
for every $\xi \in C^{\infty}_{c}((0,\infty))$. Thus, the weak
$m$-derivative $\frac{\td u}{\td m}$ of $u$ can be characterized by
\begin{displaymath}
  [u]'=\frac{\td u}{\td m}\,\mu_{m}\qquad\text{in $\mathcal{D}'((0,\infty);X)$.}
\end{displaymath}
It is worth noting that Revuz and Yor~\cite{MR1725357} employed the
notation $g\,\mu_{m}$ instead of $\frac{\td u}{\td m}$ to study the
\emph{generalized Sturm-Liouville problem}~\eqref{eq:54Intro} below
for given Radon measure $\mu_{m}$. 
\end{remark}

Now, we are ready to introduce the notion of a \emph{weak solution} of the
incomplete Dirichlet problem~\eqref{eq:16} associated with the extension operator
$\mathcal{A}_{m}$ on $\mathcal{X}_{+}$. 

\begin{definition}[{\em Weak solution of the incomplete Dirichlet problem}]
\label{def:weak-solution Dirichlet-problem}
 A function $u : (0,\infty)\to X$ is called a \emph{weak solution} of
  the extension equation 
  \begin{equation}
    \label{eq:71}
        \mathcal{A}_{m}u=0\qquad\text{in $(0,\infty)$,}
  \end{equation}
  if $u\in W^{1,1}_{loc}((0,\infty);X)$ with $\tfrac{\td u}{\dz}\in
  W^{1,1}_{loc,\mu_{m}}((0,\infty);X)$ satisfying
  \begin{displaymath}
    u(z)\in D(A)\quad\text{ and }\quad B_{m}u(z)=A(u(z))\quad
    \text{ for $\mu_{m}$-a.e. $z\in (0,\infty)$.}
\end{displaymath}
Further, for given $f\in X$, we define a function
$u\in C([0,\infty);X)$ to be a \emph{weak solution} of Dirichlet
problem~\eqref{eq:16} for the extension operator $\mathcal{A}_{m}$
defined by~\eqref{eq:137} if $u(0)=f$ in $X$,
$u\in W^{1,1}_{loc}([0,\infty);X)$ with
$\tfrac{\td u}{\dz}\in W^{1,1}_{loc,\mu_{m}}([0,\infty);X)$, and $u$ is
a weak solution of~\eqref{eq:71}.
\end{definition}

For our next definition, it is worth noting that for a function
\begin{displaymath}
u\in W^{1,1}_{loc}([0,\infty);X)\quad\text{ with }\quad
\tfrac{\td u}{\dz}\in W^{1,1}_{loc,\mu_{m}}([0,\infty);X),
\end{displaymath}
the right hand-side derivative $\frac{\td u}{\td z}_{\! +}(0)$
exists. We refer to Remark~\ref{REMARK:SOLN} for more details on this.

\begin{definition}[{The Dirichlet-to-Wentzell operator associated with $\mathcal{A}_{m}$}]
  \label{def:DtW-operator}
  Let $D(\Lambda_{m})$ be the set of all $f\in D(A)$ such that there
  exists a unique weak solution $u\in C([0,\infty);X)$ of Dirichlet
  problem~\eqref{eq:16} with boundary value $u(0)=f$. Then, we call
  the linear operator $\Lambda_{m} : D(\Lambda_{m})\to X$ defined
  by~\eqref{eq:91}, where $u$ is the weak solution of~\eqref{eq:49}
  with boundary value $u(0)=f$, the \emph{Dirichlet-to-Wentzell
    operator associated with $\mathcal{A}_{m}$}.
\end{definition}

\begin{remark}
  For strings $m\in \mathfrak{m}_{\infty}$ with right hand-side limit
  $m(0+)=0$, the Dirichlet-to-Wentzell
operator $\Lambda_{m}$ reduces to the classical
\emph{Dirichlet-to-Neumann operator} associated with $\mathcal{A}_{m}$.
\end{remark}

It is convenient to apply the preceding two definitions at a fundamental example.

\begin{example}\label{EX:dirac}
  Let the string $m\in \mathfrak{m}_{\infty}$ on $\R$ be the \emph{Heaviside step
    function}
  \begin{equation}
    \label{eq:150}
    m(x):=
    \begin{cases}
      1 & \text{if $x\ge 0$,}\\
      0 & \text{if $x<0$.}
    \end{cases}
  \end{equation}
  The associated measure
  $\mu_{m}=\delta_{0}$ is the \emph{Dirac measure} at $x=0$. Hence,
  according to the two Definitions~\ref{def:weak-dfdm-intro}
  and~\ref{def:weak-solution Dirichlet-problem}, a function
  $u\in C([0,\infty);X)\cap W^{1,1}_{loc}((0,\infty);X)$ is a weak
  solution of extension equation~\eqref{eq:71} provided
  $\tfrac{\td u}{\dz}\in W^{1,1}_{loc,\mu_{m}}((0,\infty);X)$,
  $u(z)\in D(A)$ for every $z>0$, and
  \begin{equation}
    \label{eq:60}
  -\frac{1}{2}\int_{0}^{\infty}
  \frac{\td u}{\dz}\,\frac{\td \xi}{\dz}\,\dz=0 \qquad
  \text{for every $\xi\in C^{\infty}_{c}((0,\infty))$.}
  \end{equation}
  By Lemma~\ref{lem:variational-lemma-xi-dash} and \eqref{eq:60} 
  implies that $u$ is constant on $[0,\infty)$. Thus, for given $f\in X$,
  the unique bounded weak solution $u$ of Dirichlet problem~\eqref{eq:16} is given
  by $u(z)\equiv f$. This implies that the associated Dirichlet-to-Wentzell
  operator $\Lambda_{m}$ associated with $\mathcal{A}_{m}$ reduces to
  \begin{displaymath}
    \Lambda_{m}f=m(0+)Au(0)
    =Af \quad \text{for all $f\in D(A)$.}
  \end{displaymath}
  In other words, the Dirichlet-to-Wentzell operator $\Lambda_{m}$
  associated with $\mathcal{A}_{m}$ coincides with the operator $A$ if
  $m$ is the Heaviside step function. This can be understood as a time
  change with the trivial subordinator $\textrm{id}_{\R_{+}}(s):=s$,
  $s\in \R_{+}:=[0,\infty)$. We refer the interested reader to
  Subsection \ref{subsec:inverse-local-times} for further details
  regarding the notion of \emph{inverse local times}. In
  Example~\ref{ex:dirac-revisited} we provide an alternative proof of
  this case by using one of our main results in this paper.
\end{example}

Our first theorem provides sufficient conditions to ensure
the existence, uniqueness, and a Poisson formula of bounded weak solutions of
the Dirichlet problem~\eqref{eq:16} for the extension operator
$\mathcal{A}_{m}$ defined by~\eqref{eq:137} on $\mathcal{X}_{+}$.

\begin{theorem}
  \label{thm:1DP}
  Let $\{e^{-tA}\}_{t\ge 0}$ be a $C_{0}$-semigroup of contractions on
  $X$ and $-A$ its infinitesimal generator on $X$. Further, let
  $m\in \mathfrak{m}_{\infty}$ be a string on $\R$, $\omega_{\tau}$ and
  $\psi_{m}$ be the first hitting time density and the complete
  Bernstein function associated with the string $m$. Then, for every
  $f\in D(A)$, the Dirichlet problem~\eqref{eq:16} admits a unique
  bounded weak solution $u$, and this solution $u$ is given by
  the Poisson formula~\eqref{eq:17}.
\end{theorem}

We outline the proof of Theorem~\ref{thm:1DP} in
Section~\ref{sec:proofs}, in which we combine arguments from
stochastic analysis with tools from nonlinear functional analysis in a
refined way. We begin in Section~\ref{subsec:uniqueness} by
establishing uniqueness of bounded weak solutions of Dirichlet
problem~\eqref{eq:16}; see Theorem~\ref{thm:uniqueness}. 

To the best
of our knowledge, uniqueness results for Dirichlet
problem~\eqref{eq:16} are only known in the case when the string $m$
is given by~\eqref{eq:113}, which corresponds to the fractional power
case $A^{\sigma}$; either $A$ being a sectorial operator
(see~\cite{MR3772192} and~\cite{MR4151098}), or $A$ being an
$m$-accretive (possibly, nonlinear and multi-valued) operator on a
Hilbert space (see~\cite{MR4026441}). Thus
Theorem~\ref{thm:uniqueness} essentially improves this result within
the class of accretive operators $A$ on a Banach space $X$.

Our proof
of this result uses arguments from nonlinear functional analysis and
exploits the fact that $A$ is accretive. In
Section~\ref{General_CASE}, we establish existence of bounded weak
solutions of Dirichlet problem~\eqref{eq:16} by simply verifying that
if $u$ is given by the integral~\eqref{eq:17} then $u$ satisfies
another integral representation (see Theorem~\ref{prop:existence}),
which implies that $u$ satisfies all conditions of being a weak
solutions of Dirichlet problem~\eqref{eq:16}. Then, by our uniqueness
result (Theorem~\ref{thm:uniqueness}), every weak solution
of~\eqref{eq:16} can be represented by the Poisson
formula~\eqref{eq:17}. Since the weight $\omega_{\tau}(t,z)$ in the
integral~\eqref{eq:17} is the density function of the first hitting
time $\tau$ by a given generalized diffusion, this method
provides a \emph{first hitting time approach}.

\begin{remark}
  In order to keep this paper well organized, we separated our
  research results obtained in the interesting Hilbert case $X=H$ (see
  \cite{hauerLee2021}) from
  the one presented here. By focusing on sectorial operators $A$ on a
  Hilbert space $H$ which are defined by a continuous, coercive form
  $\mathcal{E} : V\times V\to \C$, where $V$ is another Hilbert space
  continuously and densely embedded into $H$, we obtain existence of a
  weak (variational) solution $u$ of Dirichlet problem~\eqref{eq:16},
  admitting stronger regularity properties. In addition, in the case $H$ is
  separable and $A$ admits a compact resolvent, applying the
  Fourier series yields that finding a weak solution $u$ of
  Dirichlet problem~\eqref{eq:16} becomes equivalent to determining for
  every eigenvalue $\lambda\ge 0$ of $A$ the unique bounded weak
  solution $\phi \in C([0,\infty);\R)$ of the \emph{generalized
    Sturm-Liouville problem}
\begin{equation}\label{eq:54Intro}
\begin{cases}
   \;-\tfrac{1}{2}\phi'' +\lambda\,\phi\,\mu_{m}=0 &
    \quad \text{in $(0,\infty)$,}\\
    \;\hspace{1.8cm}\phi(0)=1.
  \end{cases}
\end{equation}
We note that similar ideas as given in \cite{hauerLee2021}
can be applied to pseudo-differential operators $A$ on $L^2(\R^d)$
with a strictly positive symbol. By using the Fourier transform, one
reduces the Dirichlet problem~\eqref{eq:16} to a Sturm-Liouville
problem. In fact, this was essentially done in \cite[Section
7]{MR2354493}.
\end{remark}

Next, we intend characterizing the operator $\psi_{m}(A)$ among the
class of complete Bernstein functions $\psi_{m}$ in terms of the
Dirichlet-to-Wentzell operator $\Lambda_{m}$. Before doing this, we
briefly recall from~\cite{MR2978140} the following definition.

\begin{definition}\label{def:1}
 A function $\psi : (0,\infty)\to \R$ is called a \emph{Bernstein
   function with L\'evy triple $(a,b,\nu)$}
 provided $\psi$ is given by
 \begin{displaymath}
   \psi(\lambda)=a+ b\lambda+\int_{0}^{\infty}\Big(1-e^{-\lambda
     r}\Big)\,\td\nu(r)\qquad\text{for all $\lambda> 0$,} 
 \end{displaymath}
 for some $a$, $b\ge 0$ and a L\'evy measure $\nu$ on $(0,\infty)$
 with finite integral
 $\int_{0}^{\infty}(r\wedge
 \mathds{1}_{(0,\infty)}(r))\,\td\nu(r)$. The triple $(a,b,\nu)$ is
 referred to as a L\'evy triple. Now, a Bernstein function $\psi$ is
 called \emph{complete} if it has the L\'evy triple
 $(0,b,\nu)$ and admits the following two properties
 \begin{enumerate}[label={(\roman*)}]
 \item the L\'evy measure $\nu$ is absolutely continuous with respect to the
   Lebesgue measure;
   \item the density $h:=\frac{\td\nu}{\dr}$ of the L\'evy measure
     $\nu$ is a \emph{completely monotone} function; that is,
   $h$ is smooth and its derivatives satisfy
   \begin{displaymath}
     (-1)^n\, h^{(n)}(r)\ge 0 \qquad \text{for all}\,\,n\in \N\cup\{0\},\text{
       and }r>0.
   \end{displaymath}
 \end{enumerate}
  We denote the set of complete Bernstein functions 
  by $\mathcal{C}\mathcal{B}\mathcal{F}$.
\end{definition}

Due to Kre\u{\i}n's theorem~\cite{MR0054078} (see
also~\cite[Theorem~1.1]{MR661628}), for every string
$m \in \mathfrak{m}_{\infty}$ on $\R$, there is a unique complete Bernstein
function
\begin{equation}
  \label{eq:48ast}
  \psi_{m}(\lambda):=m(0+)\lambda + \int_{0}^{\infty}\Big(1-e^{-\lambda
   r}\Big)\,\td\nu_{m}(r),\qquad\text{$\lambda> 0$.}  
\end{equation}
Moreover, the mapping
$\Psi : \mathfrak{m}_{\infty}\to \mathcal{C}\mathcal{B}\mathcal{F}$
$m\mapsto \Psi(m)=\psi_{m}$ with $\psi_{m}$ defined
by~\eqref{eq:48ast} is bijective, and called the \emph{Kre\u{\i}n's correspondence}.

\begin{notation}\label{def:associate-CBF}
  For a given string $m \in \mathfrak{m}_{\infty}$ on $\R$, we refer to the
  function $\psi_{m}: [0,\infty)\to\R$ given in~\eqref{eq:48ast} as
  the \emph{complete Bernstein function associated  with $m$}, $\nu_{m}$
  the \emph{L\'evy measure associated with $m$},
  $h_{m}=\frac{\td\nu_{m}}{\dr}$ the density of the L\'evy measure
  $\nu_{m}$, and by $(0,m(0+),\nu_{m})$ the \emph{L\'evy triple} associated with
  $\psi_{m}$.
\end{notation}

According to~\cite[Theorem~5.2]{MR2978140}, $\psi$ is a Bernstein
function if and only if there exists a unique \emph{vaguely continuous
convolution semigroup $\{\gamma_{t}\}_{t\ge 0}$ of sub-probability
measures} $\gamma_{t}$ on $[0,\infty)$ (see
Definition~\ref{def:family-of-sub-probabilities}) such that the
Laplace transform of $\gamma_{t}$
\begin{equation}
    \label{eq:50}
    	\int_{0}^{\infty}e^{-\lambda
          s}\,\td\gamma_{t}(s)=e^{-t\psi(\lambda)}\qquad
        \text{for all $\lambda> 0$, $t\ge 0$.}
\end{equation}
In addition, by Knight's theorem (see Theorem~\ref{thm:2} in
Section~\ref{subsec:inverse-local-times}) for a given string
$m\in \mathfrak{m}_{\infty}$ on $\R$, if $\{Z_{t}\}_{t\ge 0}$ denotes the
generalized diffusion associated with $m$,
$\{\tilde{L}^{-1}_{t}\}_{t\ge 0}$ the \emph{local inverse time at
  zero} of $\{Z_{t}\}_{t\ge 0}$ (see
Section~\ref{subsec:inverse-local-times} for the construction of this
notion), and if $\psi_m$ is the complete Bernstein function associated
with $m$, then the Laplace transform determines uniquely
(see~\eqref{eq:134}) that the convolution semigroup $\{\gamma_{t}\}_{t\ge 0}$ of
sub-probability measures $\gamma_{t}$ on $[0,\infty)$ associated with
$\psi_m$ has to be given by the push-forward measure
  \begin{equation}
    \label{eq:135}
    \gamma_{t}((a,b])
     =\mathbb{P}(\tilde{L}^{-1}_{t}\in (a,b])
    \qquad\text{for all $a,b \in [0,\infty)$, and $t\ge 0$.}
\end{equation}

Now, for a given $C_{0}$-semigroup $\{e^{-tA}\}_{t\ge 0}$ of
contractions $e^{-tA}\in \mathcal{L}(X)$ with infinitesimal generator
$-A$, and a vaguely continuous convolution
semigroup $\{\gamma_{t}\}_{t\ge 0}$ of sub-probability measures on
$[0,\infty)$, the family $\{e^{-t\psi(A)}\}_{t\ge 0}$ of operators
$e^{-t\psi(A)}$ on $X$ defined by the Bochner integral
\begin{equation}
  \label{eq:112}
  e^{-t\psi(A)}f:=\int_{[0,\infty)}e^{-sA}f\,\td\gamma_{t}(s)\qquad\text{for
    every $t\ge 0$ and $f\in X$}
\end{equation}
defines a $C_{0}$-semigroup of contractions on $X$
(see~\cite[Theorem~13.1]{MR2978140}).  Thanks to Phillips'
subordination theorem~\cite{MR2512800} (cf, Theorem~\ref{thm:Phillips}
in Section~\ref{subsec:subordination}), the abstractly defined
infinitesimal generator $-\psi(A)$ of the semigroup
$\{e^{-t\psi(A)}\}_{t\ge 0}$ can be expressed by
\begin{equation}
  \label{eq:109}
    \psi(A)f=af+bAf+\int_{0}^{\infty}\big(f-e^{-tA}f\big)\,\td\nu(t)
\end{equation}
for all $f\in D(A)$, where the integral~\eqref{eq:109} is to be understood in the
Bochner sense, $\psi$ is the unique Bernstein function
associated with the given vaguely continuous convolution semigroup
$\{\gamma_{t}\}_{t\ge 0}$ through~\eqref{eq:50}, and $(a,b,\nu)$ the
corresponding L\'evy triple.\medskip 

Now, for a given string $m\in \mathfrak{m}_{\infty}$ on $\R$ with associated
Bernstein function $\psi_{m}$, our second
main result provides an alternative characterization of the operator
$\psi_{m}(A)$ in terms of the Dirichlet-to-Wentzell operator
$\Lambda_{m}$.

\begin{theorem}\label{thm:1}
  Let $A$ be an $m$-accretive operator on a Banach space $X$.  Given
  a string $m\in \mathfrak{m}_{\infty}$ on $\R$, let $\psi_{m}$
  be the corresponding complete Bernstein function with L\'evy triple
  $(0,m(0+),\nu_{m})$, and $\psi_{m}(A)$ the operator given
  by~\eqref{eq:109} for this triple. Then, the
  operator $\psi_{m}(A)$ coincides with the Dirichlet-to-Wentzell
  operator $\Lambda_m$ given by~\eqref{eq:91}. 
  Moreover, the semigroup $\{e^{-t\psi_{m}(A)}\}_{t\ge 0}$ generated
  by $-\psi_{m}(A)$ is given by
  \begin{equation}
    \label{eq:105}
    e^{-t\psi_{m}(A)}f
    =\mathds{E}\left(e^{-\tilde{L}^{-1}_{t}\! A}f\right)=\int_{[0,\infty)}e^{-sA}f\,\td
    \gamma_{t}(s)
  \end{equation}
  for every $t\ge 0$, $f\in X$, where the convolution semigroup
  $\{\gamma_{t}\}_{t\ge 0}$ is given by~\eqref{eq:135} involving
  the local inverse time $\{\tilde{L}^{-1}_{t}\}_{t\ge 0}$ at zero of
  the generalized diffusion process $\{Z_{t}\}_{t\ge 0}$ associated
  with $m$.
\end{theorem}


%
%
%
%

\begin{example}[{Example~\ref{EX:dirac} revisited}]\label{ex:dirac-revisited}
  Thanks to \eqref{eq:105} in Theorem \ref{thm:1}, we can now give an alternative proof
  of the fact that the the Dirichlet-to-Wentzell operator $\Lambda_{m}$
  for the Heaviside step function $m$ given by~\eqref{eq:150}
  coincides with $A$. Namely, for the semigroup
$\{e^{-t\Lambda_{m}}\}_{t\ge 0}$ generated by $-\Lambda_{m}$ on $X$,
one has then that
\begin{align*}
	e^{-t\Lambda_m}f
	&=\mathbb{E}\Big (e^{-\frac{1}{2}\int_{[0,\infty)}L_{L_{t}^{-1}}(z)\,\td \delta_{0}(z)A}f\Big )\\
	&=\mathbb{E}\Big (e^{-\tfrac{t}{2}A}f\Big )\\
	&=e^{-\tfrac{t}{2}A}f
\end{align*}
for every $f\in X$ and $t\ge 0$, yielding that
$\Lambda_m=\frac{1}{2}A$.
\end{example}

Due to Kre\u{\i}n's correspondence $\Psi$ and by
Theorem~\ref{thm:1}, we obtain the following new characterization of
Phillips's subordination theorem characterizing $\psi(A)$ for any
$\psi$ of the class $\CBF$.

\begin{corollary}
  \label{cor:Phillips-for-CBF} Let $A$ be an $m$-accretive operator on
  a Banach space $X$.  If $\psi$ is a complete Bernstein function with
  L\'evy triple $(0,b,\nu)$ and $m\in \mathfrak{m}_{\infty}$ the
  unique string on $\R$ given by the Kre\u{\i}n's correspondence
  $\Psi(m)=\psi$, then the two operators $\psi(A)$ given
  by~\eqref{eq:109} and the Dirichlet-to-Wentzell operator $\Lambda_m$
  given by~\eqref{eq:91}
  coincide. 
\end{corollary}

\begin{remark}[{\em Problem~\ref{fp:A} \& the Dirichlet-to-`Wentzell-Robin' operator.}]\label{rem:solution1}
  Since for the string $m$ introduced by~\eqref{eq:113}, the
  Dirichlet-to-Wentzell operator~\eqref{eq:91} reduces to the
  Dirichlet-to-Neumann operator~\eqref{eq:94}, Theorem~\ref{thm:1}
  provides an answer to {\bfseries Problem~\ref{fp:A}}. Moreover,
  given a string $m\in \mathfrak{m}_{\infty}$ on $\R$ and $m$-accretive
  operator $A$ on a Banach space $X$, translating the
  associated Dirichlet-to-Wentzell operator $\Lambda_{m}$ by $\alpha\in \R$
  leads to the \emph{Dirichlet-to-`Wentzell-Robin' operator}
  $W_{\alpha,m}$ (cf, \cite{MR1987500}) given by
  \begin{equation}\label{eq:127}
    W_{\alpha,m}f= (\alpha\,\textrm{id}_{X}+\Lambda_{m})f\qquad
    \text{for every $f\in D(\Lambda_{m})$.}
  \end{equation}
  If $\psi_{m}$ is the complete Bernstein
  function associated with $m$, having L\'evy triple
  $(0,m(0+),\nu_{m})$, then
  $\psi_{\alpha,m}(\lambda):=\alpha+\psi_{m}(\lambda)$, $\lambda\ge 0$, is a complete Bernstein
  function with L\'evy triple $(\alpha,m(0+),\nu_{m})$. Thus,
  according to Theorem~\ref{thm:1}, the
  operator $\psi_{\alpha,m}(A)$ defined by~\eqref{eq:109} can be
  characterized by
  \begin{displaymath}
    \psi_{\alpha,m}(A)f=(\alpha\,\textrm{id}_{X}+\Lambda_{m})f\qquad\text{for every $f\in D(A)$.}
  \end{displaymath}
  \end{remark}

 The main results of this paper and and the preceding Remark~\ref{rem:solution1} lead
  naturally to the following open problem.
  
\begin{op}
      Given a general Bernstein function $\psi$ on $\R$ and an
        accretive operators $A$ on a Banach space $X$, can the
        operator $\psi(A)$ still be characterized as an operator
        similar to the Dirichlet-to-`Wentzell-Robin' operator
        $W_{\alpha,m}$, which is associated with an extension
        problem similar to~\eqref{eq:16}?
\end{op}

We conclude this paper with a stability result of the operator
$\psi_{m}(A)$ by varying the string $m$.




\begin{theorem}
  \label{prop:limits_str}
  For a given sequence
  $\{m_{n}\}_{n\ge 1}\subseteq \mathfrak{m}_{\infty}$ of strings
  $m_{n}$ on $\R$ of infinite length, let $\psi_{m_{n}}$ be the
  corresponding complete Bernstein functions with L\'evy triple
  $(0,m_n(0+),\nu_{m_{n}})$. Further, for a given $m$-accretive
  operator $A$ on $X$, let $\psi_{m_{n}}(A)$ be the operator given
  by~\eqref{eq:109} for these triples, and
  $\{e^{-t\psi_{m_{n}}(A)}\}_{t\ge 0}$ the semigroups generated by
  $-\psi_{m_{n}}(A)$. 
    If there is a string
    $m \in \mathfrak{m}_{\infty}$ on $\R$ of infinite length such that
    \begin{equation}
      \label{eq:116}
    \lim_{n\to\infty}m_{n}(z)= m(z)\qquad \text{pointwise for every
    continuity point $z$ of $m$,}
   \end{equation}
   then there is a complete Bernstein function $\psi_{m}$ with L\'evy triple
  $(0,m(0+),\nu_{m})$ such that for every $f\in D(A)$,
    there exists an $f_{n}\in D(\psi_{m_{n}}(A))$ such that
    \begin{equation}
      \label{eq:143}
      \lim_{n\to \infty}x_{n}=x\quad\text{ in $X$}\quad\text{ and }\quad 
      \lim_{n\to \infty}\psi_{m_{n}}(A)f_{n} =\psi_{m}(A)f\quad\text{ in $X$.}
    \end{equation}
\end{theorem}

The next remark is a reminder for later reference.

\begin{remark}\label{rem:Trotter-Kato}
  For a better understanding of the type of convergence obtained in~\ref{eq:143} of
  Theorem~\ref{prop:limits_str}, it worth recalling Trotter-Kato's
  first approximation theorem (cf, \cite[Theorem~1.8 in Chapter
  IV.]{MR2229872}). Accordingly to this theorem, the following
  statements are equivalent:
  \begin{enumerate}
  \item $\psi_{m_{n}}(A)$ converges to $\psi_{m}(A)$ in the
    \emph{graph sense}, that is, for every $f\in D(A)$, there exists
    an $f_{n}\in D(\psi_{m_{n}}(A))$ such that \eqref{eq:143} holds;
  \item $\psi_{m_{n}}(A)$ converges to $\psi_{m}(A)$ \emph{strongly in
      the resolvent sense}, that is, for every $f\in X$, and some or all
    $\lambda>0$,
    \begin{displaymath}
    \lim_{n\to
      \infty}R(\lambda,\psi_{m_{n}}(A))f=R(\lambda,\psi_{m}(A))f\qquad\text{in $X$,}
  \end{displaymath}
  where for every $\lambda>0$ and $n\ge 1$,
  $R(\lambda,\psi_{m_{n}}(A)):=(\lambda+\psi_{m_{n}}(A))^{-1}$
  denotes the \emph{resolvent} operator of $\psi_{m_{n}}(A)$ and
  $R(\lambda,\psi_{m}(A))$ the resolvent operator of $\psi_{m}(A)$;
  \item For every $f\in X$, $e^{-t\psi_{m_{n}}(A)}f\to e^{-t \psi_{m}(A)}f$ in $X$
    uniformly for $t$ in compact subintervals of $[0,\infty)$.
  \end{enumerate}
\end{remark}

We outline the proof of Theorem~\ref{prop:limits_str} in Section~\ref{subsec:stability}. 



%

\subsection{Organization of this paper}
The structure of this paper is as follows. In the subsequent
Section~\ref{sec:motiviation}, we provide a historical development of
the so-called \emph{extension technique}, mention important related
contributions, and discuss various approaches to establish this
technique. Section~\ref{sec:preliminaries} is dedicated to collect
various notions and intermediate results, which are necessary to prove
our main results in Section~\ref{sec:proofs}. In particular, we recall
the construction of \emph{generalized diffusion processes}
(Section~\ref{SEC:PRIMER_STRINGS}), and we briefly review the notion
of an \emph{inverse local time} and it relation to complete Bernstein
functions (Section~\ref{subsec:inverse-local-times}). Since
in~\cite{MR2354493}, the idea to derive the extension
equation~\eqref{eq:93} relies on applying the \emph{change of
  variable} $z=(y/2\sigma)^{2\sigma}$ to the solution $u$ of the
Bessel equation~\eqref{eq:138}, we provide in the two
Sections~\ref{subsec:scale-functions} \&~\ref{subsec:speed-measures} a
stochastic point of view of the impact of such a change of variable on
the extension equation~\eqref{eq:137} and illustrate its usefulness
later in Section~\ref{section:application} on two
applications. Section~\ref{subsec:hitting-times} is dedicated to the
hitting time $\tau$ of generalized diffusion $\{Z_{t}\}_{t\ge 0}$
associated with a string $m$ on $\R$, the probability density
$\omega_{\tau}$ of $\tau$, the transition density $\hat{p}$ of the
killed process $\{\hat{Z}_{t}\}_{t\ge 0}$ of $\{Z_{t}\}_{t\ge 0}$ and
derive the spectral representations of the hitting time density
$\omega_{\tau}$. The intermediate results gathered in this subsection
are crucial to prove our main theorems in this paper. In the
sections~\ref{subsec:semigroups} \&~\ref{subsec:subordination}, we
briefly review the notions of $C_{0}$-semigroups, $m$-accretive
operators on Banach spaces, and subordination of semigroups. As
mentioned above, Section~\ref{sec:proofs} is dedicated to the proof of
our two main results Theorem~\ref{thm:1DP} and Theorem~\ref{thm:1}.  In
Section~\ref{section:application}, we apply our main theorems for
providing a short proof of the classic limit $A^{\sigma}\to A$ in the
\emph{graph sense} as $\sigma\to 1-$. We also discuss briefly the case
$A^{\sigma}\to \textrm{id}_{X}$ in the \emph{graph sense} as
$\sigma\to 0+$.

For the reader, who is not familiar with the mixed framework of
stochastic analysis and PDEs, we provide in the
appendix of this paper a short primer on Bessel processes
(Appendix~\ref{subsec:bessel-processes}). In addition, in
Appendix~\ref{sec:def-weak}, we provide important
properties of the $m$-derivative~\eqref{eq:107}, which are not
available in the literature, but are necessary for the proofs in
this paper.

%
%

\section{Historical development of the extension technique}
\label{sec:motiviation}

Throughout this section, let $\Sigma$ be an open subset of $\R^{d}$, $d\ge 1$,
$\td\eta$ a positive measure defined on $\Sigma$, and $A$ a second
order partial differential operator such
 that $A$ is positive, densely defined, and self-adjoint in
 $L^{2}(\Sigma,\eta)$. Let $\Sigma_{+}$ denote the half-open cylinder
 $\Sigma\times (0,\infty)$ with pairs $(x,y)\in \Sigma_{+}$. Then, the
 main object of this section are the Dirichlet problem
  \begin{equation}\label{eq:1}
    \begin{cases}
      Au-y^{2\sigma-1}\{y^{1-2\sigma}u_{y}\}_{y}&=0 \qquad \text{ in
        $\Sigma_{+}$,}\\
    \hspace{2.5cm} u(x,0)&=f\qquad \text{on $\Sigma$,}
    \end{cases}
 \end{equation}
and the associated Dirichlet-to-Neumann operator
 \begin{equation}
   \label{eq:27}
   f\mapsto \Lambda_{\sigma}f:=\frac{(-1)}{2\sigma}\lim_{y\rightarrow
     0^{+}}y^{1-2\sigma}u_{y}(\cdot,y).
 \end{equation}
 We note that the above framework is the same as in Stinga and
 Torrea~\cite{MR2754080}, who extended some of the results
 in~\cite{MR2354493} on fractional powers $A^{\sigma}$ to this more
 general framework.

%
%
\subsection{A trace process approach} 
\label{subsec:trace-process}

To the best of our knowledge, Sato and Ueno~\cite[Theorem
9.1]{MR198547} were the first, who studied the \emph{trace process}
$\{X^{\ast}_{t}\}_{t\ge 0}$ on a smooth boundary $\partial\Sigma$ of a
bounded domain $\Sigma$ in a sufficiently smooth $N$-dimensional
manifold, whose generator is the negative Dirichlet-to-Neumann
operator $-\Lambda$ given by 
\begin{equation}
  \label{eq:7}
  f_{\vert \partial\Sigma}\mapsto \Lambda f:=\tfrac{\partial u_{f}}{\partial \nu}_{\vert \partial\Sigma}
\end{equation}
associated with a second-order uniformly elliptic differential
operator $A$ with symmetric coefficients on $\Sigma$. Their original aim was to
analyze second-order uniformly elliptic differential
operator operators equipped with Wentzell boundary conditions and their
corresponding stochastic processes. They showed that $\{X^{\ast}_{t}\}_{t\ge 0}$
is the Markov process obtained from a reflecting diffusion $\{X_{t}\}_{t\ge 0}$
through the \emph{time change} 
\begin{equation}
  \label{eq:9}
  X^{\ast}_{t}:=X_{L^{-1}(t)},\qquad\text{($t\ge0$),}
\end{equation}
by the local time $\{L_{t}\}_{t\ge 0}$ of $\{X_{t}\}_{t\ge 0}$ on the
 boundary $\partial\Sigma$ with
\begin{equation}
  \label{eq:10}
  L^{-1}_{t}:=\sup\{s\ge 0\,\vert\,L_{s}\le t\}\qquad\text{for all
    $t\ge 0$.}
\end{equation}

In the mid 80s, Hsu~\cite{MR837810} constructed a Skorokhod-type lemma
(see \cite[Lemma 2.1]{MR3616274} for the original result by Skorokhod
which applies to the half-space) for bounded domains $\Sigma$ in
$\R^{N}$ with a $C^{2}$-boundary $\partial\Sigma$. By applying this
lemma to a Brownian motion $\{B_{t}\}_{t\ge 0}$ inside $\Sigma$, Hsu
obtained a reflecting Brownian motion $\{X_{t}\}_{t\ge 0}$ inside of
$\Sigma$. Then by using It\^{o}'s lemma, Hsu showed that the
Dirichlet-to-Neumann operator~\eqref{eq:7} associated with the
(scaled) Laplacian $\frac{1}{2}\Delta$ is the infinitesimal generator
of the trace process $\{X^{\ast}_{t}\}_{t\ge 0}$ again obtained from
$\{X_{t}\}_{t\ge 0}$ via~\eqref{eq:9} through the time
change~\eqref{eq:10} (see \cite[Proposition 4.1]{MR837810}).

To see the two-step construction of the Dirichlet-to-Neumann operator
via stochastic analysis, it is worth recalling that a reflecting
Brownian motion $\{X_{t}\}_{t\ge 0}$ inside a domain $\Sigma$ is a
classical Brownian motion before it hits the first time the boundary
$\partial\Sigma$. But Brownian motion killed at the boundary
$\partial\Sigma$ corresponds to the classical Dirichlet problem for the
(scaled) Laplacian $\frac{1}{2}\Delta$ on $\Sigma$, and reflecting
Brownian motion corresponds to the Neumann problem involving the
(scaled) Laplacian $\frac{1}{2}\Delta$.\medskip

The trace process $\{X^{\ast}_{t}\}_{t\ge 0}$ on the bottom $\Sigma$
of the cylinder $\Sigma_{+}$ generated by the Dirichlet-to-Neumann
operator $\Lambda_{\sigma}$ given in~\eqref{eq:27} is obtained by starting
with the \emph{joint process}
\begin{displaymath}
\{\X_{t}\}_{t\ge
  0}:=\{(X_{t},Y_{t})\}_{t\ge0}\qquad\text{ in $\Sigma_{+}$}
\end{displaymath}
for given processes $\{X_{t}\}_{t\ge 0}$ generated by the operator $A$
in $\Sigma$ and $\{Y_{t}\}_{t\ge 0}$ generated by the Bessel-operator
$\mathcal{B}_{1-2\sigma}=\tfrac{1}{2}\left (\tfrac{\td^2}{\td
    y^2}+\tfrac{1-2\sigma}{y} \right )$ on the half-line
$(0,+\infty)$.

To the best of our knowledge, the joint process $\{\X_{t}\}_{t\ge 0}$
occurred the first time in the short paper~\cite{MR0247668} by Mol\v{c}anov and
Ostrovski\u{\i}. In the case $\Sigma$ is the Euclidean space
$\R^{N}$, they proved that the trace process $\{X^{\ast}_{t}\}_{t\ge 0}$ on
$\Sigma=\R^{N}$ obtained from~\eqref{eq:9} through the time
change~\eqref{eq:10} of the local time $\{L_{t}\}_{t\ge 0}$ of
$\{Y_{t}\}_{t\ge 0}$ at $y=0$ is generated by the fractional Laplace
operator $(-\Delta)^{\sigma}$ on $\R^{N}$.

In the literature, one finds claims that instead of the highly-cited
paper~\cite{MR2354493} by Caffarelli and Silvestre, the approach by
Mol\v{c}anov and Ostrovski\u{\i} was the first (stochastic) proof to
the extension property of $(-\Delta)^{\sigma}$ (i.e.,
Theorem~\ref{thm:1} applied to $A=-\Delta$ on $\Sigma$).  But, the
final step, the identification of the infinitesimal generator of the
trace process with the Dirichlet-to-Neumann map $\Lambda_{s}$, was not
a considered aim in~\cite{MR0247668}. We agree with the sentiment in
\cite{assing2019extension} that the proof of this identification is
far from obvious. In fact, if one intends to prove the result in
\cite{MR2354493} by using the trace process of the joint process given
by $\{(B_{t},Y_{t})\}_{t\ge 0}$ in $\R^{d+1}_{+}$, where
$\{B_{t}\}_{t\ge 0}$ is a Brownian motion, then It\^o's lemma cannot
be applied immediately since the Bessel process $\{Y_{t}\}_{t\ge 0}$
is no longer a semi-martingale when $\delta \in (0,1)$ (see
Appendix~\ref{subsec:bessel-processes} for more details). But a
potentially more fruitful approach would instead be to consider the
trace process of the joint process
$\{(B_{t},Y^{2\sigma}_{t})\}_{t\ge 0}$ in $\R^{d+1}_{+}$ since
according to Theorem \ref{thm:Donati-yor}, the process
$\{Y^{2\sigma}_{t}\}_{t\ge 0}$ is a submartingale. Nevertheless, in
this setting one still needs to overcome certain technical details, in
order to apply the ideas in \cite{MR837810}. We direct the reader to
\cite{assing2019extension} for details. The scaling of the process
$\{Y_{t}\}_{t\ge 0}$ to $\{Y^{2\sigma}_{t}\}_{t\ge 0}$ is a common
  technique in stochastic analysis outlined in
  Section~\ref{subsec:scale-functions} and
  Section~\ref{subsec:speed-measures} and was also used
  in~\cite{MR2354493}.

%
%

\subsection{A Fourier approach}
\label{subsec:fourier}

A first approach to use the Fourier transform to prove
Theorem~\ref{thm:1} in the case of the negative Laplace operator
$A=-\Delta$ on $\Sigma=\R^{N}$ was provided by Caffarelli and
Silvestre~\cite{MR2354493}. Under the additional assumption that the
operator $A$ has a discrete spectrum
$\sigma(A)=\{\lambda_{k}\}_{k\ge 0}$, Stinga and
Torrera~\cite{MR3709888} outlined the following Fourier-series
approach to the existence and uniqueness of solutions $u$ to the
Dirichlet problem~\eqref{eq:1}. Under this hypothesis, there is an
orthonormal basis $\{\phi_{k}\}_{k\ge 0}$ in $L^2(\Sigma,\mu)$ such
that $A\phi_{k}=\lambda_{k}\phi_{k}$ for all $k\ge 0$. Then, for given
$f\in L^2(\Sigma,\mu)$, the Fourier-series
\begin{equation}
  \label{eq:12}
  u(\cdot,y)=\sum_{k\ge 0}c_{k}(y)\phi_{k}\qquad\text{converging in
  $L^{2}(\Sigma,\eta)$ for every $y>0$,}
\end{equation}
is the unique solution of~Dirichlet
problem~\eqref{eq:1} if and only if for every $k\ge 0$, $c_{k}$ is the
unique solution of the Dirichlet problem
\begin{equation}
  \label{eq:11}
  \begin{cases}
    \hspace{1cm}\mathcal{L}_{\sigma,k}c_{k}&=0, \qquad  \qquad  \qquad\text{on $(0,\infty)$,}\\
    \hspace{1.2cm} c_k(0)&=\langle
      f,\phi_{k}\rangle_{L^{2}(\Sigma,\eta)},\\[2pt]
    \hspace{0.2cm}\displaystyle\lim_{y\to+\infty}c_{k}(y)& = 0.
  \end{cases}
\end{equation}
for the \emph{Sturm-Liouville operator}
\begin{displaymath}
  \mathcal{L}_{\sigma,k}c(y):=-\frac{1}{y^{1-2\sigma}}\{y^{1-2\sigma}c'\}'+\lambda_{k}\, c.
\end{displaymath}
If $K_{\sigma}$ denotes the \emph{modified Bessel function of the third kind}, then the
unique solution $c_{k}$ to Dirichlet problem~\eqref{eq:11} is given by
\begin{displaymath}
  c_{k}(y)=y^\sigma\frac{2^{1-\sigma}}{\Gamma(\sigma)}\lambda_{k}^{\sigma/2}
  \langle f,\phi_{k}\rangle_{L^{2}(\Sigma,\eta)} K_{\sigma}(\lambda_{k}^{1/2}y).
\end{displaymath}
By using the asymptotic of $K_{\sigma}$ as $y\to0+$ and $y\to+\infty$, one
sees that the series~\eqref{eq:12} is a solution of Dirichlet
problem~\eqref{eq:1}.\medskip

Motivated by these ideas, we provide in this paper a stochastic proof
to our main result (Theorem~\ref{thm:1}). This approach simplifies
essentially the recently appeared proofs related to
Theorem~\ref{thm:1} by Assing and Herman~\cite{assing2019extension},
and Kwa\'{s}nicki and Mucha~\cite{MR3859452}.

%
%

\subsection{A stochastic representation formula}
\label{subsec:represenation-formula}

From a stochastic analytical point of view, it is worth mentioning
that the \emph{Sturm-Liouville equation} 
\begin{equation}
\label{eq:15}
   -y^{2\sigma-1}\{y^{1-2\sigma}c'\}'+\lambda_{k}\, c=0\qquad\text{on $(0,\infty)$}
\end{equation}
related to Dirichlet problem~\eqref{eq:11} is studied heavily in
relation with the \emph{first hitting time} $\tau$ at zero of the
$2(1-\sigma)$-Bessel Process $\{Y_{t}\}$ starting at $y>0$. To
simplify notation, we suppress the $y$-dependence of $\tau$ and we
refer to Appendix~\ref{subsec:bessel-processes}. In fact,
this is quite natural, since in the \emph{Poisson formula} 
\begin{equation}
    \label{eq:13}
    u(x,y)=\int_{0}^{\infty}\left(e^{-tA}f\right)\!\!(x)\,
    \omega_{\sigma}(t,y)\,dt\qquad\text{on
    $\Sigma_{+}$,}
  \end{equation}
 of the weak solution $u$ to the Dirichlet problem~\eqref{eq:1}, the density
  \begin{equation}
    \label{eq:14}
    \omega_{\sigma}(t,y)=\frac{y^{2\sigma}}{2^{2\sigma}\Gamma(\sigma)}
    e^{-\frac{y^2}{4t}}\frac{1}{t^{1+\sigma}}\qquad\text{for
    every $t$, $y>0$,} 
  \end{equation}
coincides exactly with the
\emph{first hitting time density} $\tau$ (cf.,~\cite{MR1997032}, and
\cite[No. 43 \& 44, p.75]{MR1912205}). We note that the Poisson
formula~\eqref{eq:13} was obtained in~\cite{MR2754080} by Stinga and Torrea.
In other words, for~\eqref{eq:13}, we have \allowdisplaybreaks
\begin{align*}
  u(x,y)
  &=\frac{y^{2\sigma}}{2^{2\sigma}\Gamma(\sigma)}\int_{0}^{\infty}\,(e^{-tA}f)(x)\,
    e^{-\frac{y^2}{4t}}\,\frac{\dt}{t^{\sigma+1}}\\
  &= \frac{y^{2\sigma}}{2^{2\sigma}\Gamma(\sigma)}\int_{0}^{\infty}\mathds{E}^{x}\left(f\big(X_{t}\big)\right)\,
   e^{-\frac{y^2}{4t}}\,\frac{\dt}{t^{\sigma+1}}\\
 &=\int_{0}^{\infty}\mathds{E}_{x}\left(f\big(X_{t}\big)\right)\;\td\mathds{P}_{\tau_{y/\sqrt{2}}}(t)\\
&=\int_{0}^{\infty}\mathds{E}_{x}(f(X_{t})\,|\tau_{y}\in
  dt) \,\mathds{P}(\tau_{y}\in dt)\\
&=\int_{0}^{\infty}\mathds{E}_{x}(f(X_{\tau_{y}})\,|\tau\in
  dt)\,\mathds{P}(\tau\in dt)\\
&=\mathds{E}_{(x,y)}\left(\mathds{E}_{(x,y)}\left(f\big(X_{\tau})\,
  \Big\vert\,\tau\right)\right)
=\mathds{E}_{(x,y)}\left(f\big(X_{\tau}\big)\right),
\end{align*}
where we used that $\{e^{-tA}\}_{t\ge 0}$ with
  $e^{-tA}f(x):=\mathds{E}_{x}\left(f\big(X_{t}\big)\right)$ is the
  transition semigroup of the process $\{X_{t}\}_{t\ge 0}$, and the
  independence of $\{X_{t}\}_{t\ge 0}$ and $\tau$. Thus, the
  Poisson formula~\eqref{eq:13} is nothing less than a stochastic
  representation formula to Dirichlet problem~\eqref{eq:1}. In this
  paper, we outline how to derive from Sturm-Liouville
  equation~\eqref{eq:15} the stochastic representation
  formula~\eqref{eq:17} to the more general Dirichlet
  problem~\eqref{eq:16}.

%
%

\subsection{The L\'evy measure of the inverse local time}
\label{subsec:DtN-map}

The following observation is critical for the complete Bernstein
case. Let $s(y)=y^{2\sigma}$, which is also the \emph{scale function}
of the Bessel process $\{Y_{t}\}_{t\ge 0}$ (see
Section~\ref{subsec:scale-functions} for more details). Then, for
solutions $u$ to the extension problem~\eqref{eq:1}, Stinga and
Torrera~\cite{MR3709888} showed that
\begin{displaymath}
  \lim_{y\rightarrow 0^{+}}\frac{1}{s'(y)}\frac{\partial u}{\partial
    y}(x,y)=\lim_{y\rightarrow 0^{+}}\frac{1}{s(y)}(u(x,y)-u(x,0)).
\end{displaymath}
Inserting~\eqref{eq:13} into $u(x,y)$ and $u(x,0)$, and by using that 
\begin{equation}
  \label{eq:18}
    \td \nu(t)=\frac{\mathds{1}_{\{t>0\}}}{2^{2\sigma}\Gamma (\sigma)}\frac{\dt}{t^{\sigma+1}}
\end{equation}
is a L\'evy measure, dominated convergence and Phillips' subordination theorem
yield that
\begin{align*}
	\lim_{y\rightarrow 0^{+}}\frac{1}{s'(y)}\frac{\partial
  u}{\partial y}(x,y)&=\lim_{y\rightarrow
                       0^{+}}\int_{0}^{\infty}(e^{-tA}f(x)-f(x))\,\Big(\frac{1}{s(y)}\omega_{\sigma}(t,y)\Big
                       )\,dt\\
	&=\int_{0}^{\infty}(e^{-tA}f(x)-f(x))\,\frac{1}{2^{2\sigma}\Gamma(\sigma)}\frac{1}{t^{1+\sigma}}\,dt\\
	&=-\frac{\Gamma(1-\sigma)}{2^{2\sigma}\Gamma(1+\sigma)}A^{\sigma}f(x).
\end{align*}
Here, the important observation is that for every $t>0$,
\begin{displaymath}
  \lim_{y\rightarrow
    0^{+}}\frac{1}{s(y)}\omega_{\sigma}(t,y)=\frac{1}{2^{2\sigma}\Gamma(\sigma)}\frac{1}{t^{1+\sigma}}
  =\lim_{y\rightarrow 0^{+}}\frac{1}{s'(y)}\frac{\partial\omega_{\sigma}}{\partial y}(t,y),
\end{displaymath}
where $\frac{1}{2^{2\sigma}\Gamma(\sigma)}\frac{1}{t^{1+\sigma}}$ is the density of
the measure $\nu$ given in~\eqref{eq:18}. However, it is well-known for
quite a while (cf~\cite{MR637387} or \cite{MR2417969}) that the measure $\nu$ is
the L\'{e}vy measure of the inverse local time process of the $2\sigma^{\textrm{th}}$-powered
process $\{Y_{t}^{2\sigma}\}_{\ge 0}$ of the $2(1-\sigma)$-Bessel process
$\{Y_{t}\}_{\ge 0}$ (cf, Theorem~\ref{thm:Donati-yor}).\medskip 

The goal of this paper is to exploit
the idea of the first hitting time density $\tau$ for proving the
extension technique for operators $\psi(A)$ provide $\psi$ is a complete
Bernstein function.

\section{Preliminaries}
\label{sec:preliminaries}

In this section, we introduce several important notions and collect
related results, which are required to proof our main theorems stated
in Section~\ref{subsec:main-results}. 
%
%

\subsection{A primer on generalized diffusion processes}
\label{SEC:PRIMER_STRINGS}
In this section, we briefly collect some basic notions from the theory
of generalized diffusion processes and strings needed throughout this
paper. For an excellent review of the subject, we refer the interested reader
to~\cite[Chapter 14]{MR2978140} and for a more in depth resource
see~\cite[Chapter 5 \& 6]{MR0448523}.\medskip

After providing the definition of a string $m\in \mathfrak{m}_{\infty}$ on
$\R$ (Definition~\ref{def:string}), and the associated measure
$\mu_{m}$, we still need to introduce to the following notion.

\begin{definition}\label{def:support-of-mu}
 For every string $m\in \mathfrak{m}_{\infty}$ we denote by
 \begin{displaymath}
   E_{m}:=\supp(\mu_{m})
 \end{displaymath}
 the \emph{support of $\mu_{m}$}. 
\end{definition}

\begin{remark}
We note that in the language of Feller's theory of diffusion processes
(cf., \cite[Chapter III and VII]{MR1725357}),
 the measure $\mu_{m}$ plays the role of a \emph{speed
   measure}. 
\end{remark}

With this preliminary in mind, we can start introducing generalized
diffusion processes. Fix a string $m\in \mathfrak{m}_{\infty}$, and for
given $z\in [0,\infty)$, let $\{B^{+}_{t}\}_{t\ge 0}$ be a reflecting
Brownian motion on $[0,\infty)$ starting at $z$. We refer to $\{B^{+}_{t}\}_{t\ge 0}$
as the \emph{reflecting Brownian motion}. 

It is well known that the reflecting Brownian motion
$\{B^{+}_{t}\}_{t\ge 0}$ is a Hunt process, which is symmetric with
respect to the Lebesgue measure and hence, the process
$\{B^{+}_{t}\}_{t\ge 0}$ has an associated Dirichlet form $\E$ on
$L^2((0,\infty),\dx)$. We refer to ~\cite[Appendix A.2]{MR2978140} and
the classical textbook~\cite{MR2778606} concerning Dirichlet forms and
symmetric Hunt processes. To the Dirichlet form $\E$, one can define
the \emph{capacity} $\Capacity$ by
\begin{displaymath}
\Capacity(O):=\inf\Big\{
\E(\phi,\phi)+\norm{\phi}^{2}_{L^2((0,\infty),\dx)}\,\vert\, \phi\ge
1\text{ a.e.}\big\}
\end{displaymath}
for every open subset $O$ of $[0,\infty)$. Now, it is worth mentioning
that for the string $m$, the associated measure $m$ is \emph{smooth}
in the sense that it does not charge any set of zero capacity and if
there is an increasing sequence $(A_{n})_{n\ge 1}$ of closed subsets
$A_{n}\subseteq [0,\infty)$ satisfying $m(A_{n})<\infty$ for all
$n\ge 1$ and $\lim_{n\to\infty}\Capacity(K\setminus A_{n})=0$ for all
compact subsets $K$ of $[0,\infty)$.

For every $z\in (0,\infty)$, let $\{L_{t}(z)\}_{t\ge 0}$ be the local
time process at level $z$ of the reflecting
Brownian motion $\{B^{+}_{t}\}_{t\ge 0}$. Then, the family $\{L_{t}(z)\}_{t\ge 0, 0\le
 z<\infty}$ is jointly continuous and (after normalization) the following
\emph{occupation time formula}
\begin{equation}\label{EQN:BROWNIAN_OCC_TIMES}
 \int_{0}^{t}g(B_{r}^{+})\,\dr=\int_{0}^{\infty}g(z)\, L_{t}(z)\, \dz
 \qquad \text{holds for all $t\ge0$,}
\end{equation}
and all Borel functions $g : [0,\infty)\to [0,\infty)$. Since for the
string $m$, the associated measure $m$ is \emph{smooth} and for the
reflecting Browning motion $\{B^{+}_{t}\}_{t\ge 0}$ only the empty set has capacity zero,
we can define a positive continuous additive process
$\{A_{t}\}_{t\ge 0}$ by setting
\begin{equation}
 \label{eq:21}
 A_{t}=\int_{[0,\infty)}L_{t}(z)\,\td\mu_{m}(z)\qquad\text{for all $t\ge 0$.} 
\end{equation}

It follows from our construction (cf., \cite[Chapter 5]{MR2778606})
that the \emph{Revuz measure} of the process $\{A_{t}\}_{t\ge 0}$
is the measure $\mu_{m}$.  Now, let $\{Z_{t}\}_{t\ge 0}$ be the process
defined by the \emph{time change} $\{A^{-1}_{t}\}_{t\ge 0}$, that is,
\begin{equation}
 \label{eq:19}
 Z_{t}:=B^{+}_{A^{-1}_{t}} \qquad \text{for every $t \ge 0$,}
\end{equation} 
where $\{A^{-1}_{t}\}_{t\ge 0}$ is the right-continuous inverse of
$\{A_{t}\}_{t\ge 0}$ given by
\begin{displaymath}
 A^{-1}_{t}=\inf\{s>0\,|\,A_{s}>t\},\qquad (t\ge0).
\end{displaymath}
Then, $\{Z_{t}\}_{t\ge 0}$ is an $m$-symmetric, continuous Hunt
process in $[0,\infty)$ with infinitesimal generator (see, for example, \cite[Theorem 6.2.1]{MR2778606})
\begin{equation}
 \label{eq:32}
   -B_{m}=\frac{1}{2}\frac{\td}{\dm}\frac{\td}{\dz}\qquad\text{ on
     $L^{2}(E_{m},\mu_{m})$.}
\end{equation}

\begin{definition}\label{def:gen-diffusion}
 For given $m \in \mathfrak{m}_{\infty}$, the process $\{Z_{t}\}_{t\ge 0}$
 defined by~\eqref{eq:19} is called a \emph{generalized
   diffusion process associated with $m$}.
\end{definition}

Here, the operator $B_m$ given in~\eqref{eq:32} also occurs in the next definition (cf, \cite[p. 271]{MR2978140}).

\begin{definition}\label{def:operator-B}
 Let $E$ be a non-empty, connected subset of $\R$ with $-\infty\le l:=\inf E<r:=\sup
 E\le\infty$, $m$ be a Radon measure on
 $E$, and $s: E\to \R$ a continuous, strictly increasing
 function. Then, let $W^{2,1}_{loc,ds,\mu}(E)$ be the set of functions
 $u : \R\to \C$ of the form
 \begin{displaymath}
   u(z)=\alpha+\beta z+\int_{l}^{z}\int_{(l,y]}g(r)\,\td\mu_{m}(r)\,\ds(y),\qquad(z\in \R),
 \end{displaymath}
 for some $\alpha$, $\beta\in \C$ and locally $m$-integrable
 $g : \R\to \C$. Note, the second integral w.r.t. $\ds(y)$ is a
 Lebesgue-Stieltjes integral. If $-\infty<l$ then for every $z\le l$, one interprets the interval
 $[l,y]=\emptyset$ for $y<z$ and hence, $u(z)=\alpha+\beta z$ for
 $z<l$. Then, the 2nd-order differential operator
 \begin{displaymath}
   -B_{m,s}=\frac{\td}{\dm}\frac{\td}{\ds}
 \end{displaymath}
 is defined by $-B_{m,s}u=g$ for every $u\in W^{2,1}_{loc,ds,\mu_{m}}(E)$.
\end{definition}

Moreover, the following properties hold.

\begin{proposition}\label{prop:domain-of-B}
 Let $m\in \mathfrak{m}_{\infty}$ be a string, $s=\textrm{id}_{\R}$ the
 identity on $\R$, and $B_{m}$ be the operator given by~\eqref{eq:32}
 with domain $D(B_{m})$. Then, for every $u\in D(B_{m})$, one has that
 \begin{enumerate}
   \item $u$ is locally absolutely continuous on $\R$,
   \item $u$ is linear on $\R\setminus E_{m}$,
   \item the right derivative $\frac{\td u}{\dz}_{\! +}\!$ and left derivative $\frac{\td u}{\dz}_{\! -}\!$
     exist on $\R$, and
     \begin{displaymath}
       \frac{\td u}{\dz}_{\! +}\!(x)-\frac{\td u}{\dz}_{\! -}\!(x)=\mu_m\{x\}B_{m}u(x)\qquad\text{for every $0\le
         x\le \sup E_{m}$.}
     \end{displaymath}
 \end{enumerate}
\end{proposition}

We omit the elementary proofs to the above proposition.\medskip

For a generalized diffusion process $\{Z_{t}\}_{t\ge 0}$, we can
associate for every $z\in E_{m}$, a local time process
$\{\tilde{L}_{t}(z)\}_{t\ge0}$ at level $z$, which can be realized as
a time-change of the local time process $\{L_{t}(z)\}_{t\ge0}$ of the
reflecting Brownian motion $\{B^{+}_{t}\}_{t\ge 0}$. Namely, we have that
\begin{equation}
 \label{eq:24}
 \tilde{L}_{t}(z)=L_{A^{-1}_{t}}(z)\qquad\text{for every $t\ge 0$,
   $z\in [0,\infty)$,}
\end{equation}
and that the following \emph{occupation times formula} 
\begin{equation}
 \label{eq:20}
 \int_{0}^{t}g(Z_{r})\,\dr=\int_{[0,\infty)}g(z)\tilde{L}_{t}(z)\,\td \mu_{m}(z)\qquad
\text{for all $t\ge 0$}
\end{equation}
holds for all essentially bounded Borel functions
$g:[0,\infty)\rightarrow [0,\infty)$.


\subsection{Inverse local times and complete Bernstein functions}
\label{subsec:inverse-local-times}

For a generalized diffusion process $\{Z_{t}\}_{t\ge 0}$, let
$\{\tilde{L}_{t}\}_{t\ge0}$  be the
local time $\{\tilde{L}_{t}(0)\}_{t\ge0}$ of $\{Z_{t}\}_{t\ge 0}$ at
the level $z=0$, and 
$\{\tilde{L}_{t}^{-1}\}_{t\ge0}$ be the inverse local time give by
\begin{displaymath}
 \tilde{L}^{-1}_{t}=\inf\Big\{r>0\,\Big\vert\,
 \tilde{L}_{r}>t\}\qquad\text{for every $t\ge 0$.}
\end{displaymath}
Then, by~\eqref{eq:21} and~\eqref{eq:24}, we have that
\begin{equation}\label{EQN:INV_LOC_TIME}
\tilde{L}^{-1}_{t}=A_{{L}^{-1}_{t}}=\int_{[0,\infty)}L_{L^{-1}_{t}}(z)\,\td\mu_{m}(z)\qquad\text{for
  every $t\ge 0$.}
\end{equation}
%
%
It is worth noting that the inverse local time $\{\tilde{L}_{t}^{-1}\}_{t\ge0}$
is a \emph{subordinator} (cf.,~\cite[p. 114]{MR1406564}), that is, a
one-dimensional L\'evy process that is non-decreasing (a.s.). Hence
(see~\cite[Theorem~1.3.15]{MR2512800}), there is a L\'evy measure
$\nu_m$ on $\mathcal{B}(\R\setminus\{0\})$ satisfying $\nu_{m}(-\infty,0)=0$,
\begin{equation}
 \label{eq:22}
 \int_{0}^{\infty}\Big(r\wedge \mathds{1}_{(0,\infty)}(r)\Big)\,\td\nu_m(r)<\infty,
\end{equation}
and there is some $b\ge 0$ such that the \emph{L\'evy symbol} $\eta$
of $\{\tilde{L}_{t}^{-1}\}_{t\ge0}$ is given by
\begin{displaymath}
 \eta_m(s)=i bs+\int_{0}^{\infty}\Big(e^{i sr}-1\Big)\,\td\nu_m(r)
\qquad\text{for all $s\in \R$.}
\end{displaymath}
From this, one finds that the \emph{Laplace exponent $\psi_{m}$}
of $\{\tilde{L}_{t}^{-1}\}_{t\ge0}$ is given by
\begin{equation}
 \label{eq:23}
 \psi_m(\lambda)=m(0+)\lambda+\int_{0}^{\infty}\Big(1-e^{-\lambda
   r}\Big)\,\td\nu_m(r)\qquad\text{for all $\lambda\geq 0$.} 
\end{equation}
The right hand-side in~\eqref{eq:23} of the Laplace
exponent $\psi_{m}$ provides an important example
of a \emph{Bernstein function} (cf, Definition~\ref{def:1} and, cf.,~\cite[Theorem~3.2]{MR2978140}).

\begin{remark}
  \label{rem:5}
  Due to the Bernstein's theorem of monotone functions \cite[Theorem 4.8]{MR2978140}, one has
  that for every completely monotone function $h$, there is an
  associated measure $\hat{\Delta}$ satisfying
  \begin{equation}
    \label{eq:25}
    h(r)=\int_{0}^{\infty}e^{-r\,s}\,\td \hat{\Delta}(s) \qquad\text{for all $r>0$.} 
  \end{equation}
\end{remark}


\begin{definition}\label{def:Levy-density-principal-measure}
  For given $m \in \mathfrak{m}_{\infty}$, let $\{\hat{Z}_{t}\}_{t\ge 0}$
  be the generalized process $\{Z_{t}\}_{t\ge 0}$ associated with $m$
  killed at level $0$, and $h_{m}=\frac{\td\nu_{m}}{\dr}$ be the
  density of the L\'evy measure $\nu_{m}$ associated with $m$. Then,
  we call the measure $\hat{\Delta}_{m}$ satisfying~\eqref{eq:25} for
  $h_{m}$ the \emph{principal measure} of the process
  $\{\hat{Z}_{t}\}_{t\ge 0}$.
\end{definition}



The Laplace exponent of the \emph{inverse local time processes} of a
given generalized diffusion process can be described in terms of
Bernstein functions thanks to the famous result~\cite{MR647781} by
Knight. 

\begin{theorem}[{Knight's theorem} {\cite[Theorem 3.1 and Theorem 1.2]{MR647781}}]
  \label{thm:2}
 For given string $m \in \mathfrak{m}_{\infty}$ on $\R$, let $\{Z_{t}\}_{t\ge 0}$ be a
 generalized diffusion associated with
 $m$ and $\{\tilde{L}^{-1}_{t}\}_{t\ge 0}$ be the corresponding
 inverse local time at $0$. Then the Laplace transform
 of $\tilde{L}^{-1}_{t}$ is given by
 \begin{equation}
   \label{eq:134}
	\mathbb{E}\Big(e^{-\lambda
         \tilde{L}^{-1}_{t}}\Big)=e^{-t\psi_{m}(\lambda)}\qquad\text{for
       all $\lambda> 0$,}
  \end{equation}
  where $\psi_{m}: (0,\infty)\to\R$ is a complete Bernstein function
  given by
  \begin{equation}
    \label{eq:48}
   \psi_{m}(\lambda)=m(0+)\lambda + \int_{0}^{\infty}\Big(1-e^{-\lambda
   r}\Big)\,\td\nu_{m}(r)\qquad\text{for all $\lambda> 0$.}  
  \end{equation}
  Moreover, the mapping $m\mapsto \psi_{m}$ with $\psi_{m}$ defined
  by~\eqref{eq:48} is a bijection
  $\Psi : \mathfrak{m}_{\infty}\to \mathcal{C}\mathcal{B}\mathcal{F}$,
  $\Psi(m):=\psi_{m}$ for every $m\in \mathfrak{m}_{\infty}$ known as the
  Kre\u{\i}n's correspondence.
\end{theorem}

We shall discuss properties of the principal measure
$\hat{\Delta}_{m}$ in Section \ref{subsec:hitting-times}.

%
%
%
%


\subsection{Substitution and Scale functions}
\label{subsec:scale-functions}
The aim of this section is to outline the connection between the
substitution used in~\cite{MR2354493} to simplify the extension
problem~\eqref{eq:92} and the scaling of associated stochastic
processes and to give an overview of the impact on the corresponding
Dirichlet-to-Wentzell operator (see
Section~\ref{subsec:speed-measures}).\medskip

It is clear that the substitution
\begin{displaymath}
 z=\left(\frac{y}{2\sigma}\right)^{2\sigma}\qquad
 \text{for $y>0$ and fixed $0<s<1$,}
\end{displaymath}
transforms the differential equation in the Dirichlet
problem~\eqref{eq:1} to the non-divergence form equation
\begin{displaymath}
 Au-\frac{1}{z^{-\frac{2\sigma-1}{\sigma}}}u_{zz}=0 \qquad
 \text{ in $\Sigma \times (0,\infty)$;}
\end{displaymath}
see~\cite{MR2354493}, or alternatively for more details,
see~\cite[Lemma~3.4]{MR4026441}. Moreover, Stinga and Torrea~\cite{MR2754080}
proved that for the scale function
\begin{equation}
  \label{eq:28}
  s(y)=y^{2\sigma}\qquad\text{for $y>0$ and fixed $0<\sigma<1$,}
\end{equation}
the co-normal derivative~\eqref{eq:27} can be rewritten
as 
\begin{displaymath}
 \frac{\td u}{\ds}(\cdot,0):=
  \lim_{y\rightarrow 0^{+}}\frac{u(\cdot,y)-u(\cdot,0)}{s(y)}
  =\frac{(-1)}{2\sigma}\lim_{y\rightarrow
    0^{+}}y^{1-2\sigma}u_{y}(\cdot,y).
\end{displaymath}
On the other hand, in stochastic analysis the function $s$ given
by~\eqref{eq:28} is used to \emph{scale} the $2(1-\sigma)$-Bessel
process $\{Y_{t}\}_{t\ge 0}$ in $[0,\infty)$ starting at $y=0$ in
order to derive stronger properties of the process (see
Theorem~\ref{thm:Donati-yor} in the appendix). In order to give more details t
this, we briefly need to review the necessary definition and results
on \emph{scale functions} from~\cite[Chapter 7]{MR1725357}.\medskip

For given $-\infty\le l<r\le \infty$, let $E_{l,r}$ be either a
closed, open or semi-closed interval in $\R$, and let
$\{Y_{t}\}_{t\ge 0}$ be a continuous, regular, strong Markov process
on the state space $E_{l,r}$ with a killing time $\zeta$. Further,
suppose, $\{Y_{t}\}_{t\ge 0}$ can only be killed at the end-points $l$
and $r$ of $E_{l,r}$, provided they do not belong to $E_{l,r}$. For
given $c\in E_{l,r}$, we define the \emph{first time}
$\{Y_{t}\}_{t\ge 0}$ hits $y=c$ by
\begin{equation}
  \label{eq:26}
  \tau^{c}_{Y}:=\inf\Big\{t>0\,\Big\vert\;Y_{t}=c\Big\}.
\end{equation}

In this framework the following important existence result of a
\emph{scale function} holds.

\begin{proposition}[{\cite[Chapter VII, Proposition~3.2]{MR1725357}}]
 \label{prop:existence-scale-functions}
 Under the hypothesis of this subsection, there exists a continuous,
 strictly increasing function $s: E_{l,r}\to \R$ with the property
 that for every $a$, $b$, $y\in E_{l,r}$ satisfying $l\le a<y<b\le r$,
 one has that
 \begin{equation}
   \label{eq:29}
   \mathbb{P}_{y}(\tau^b_{Y}<\tau^{a}_{Y})=\frac{s(y)-s(a)}{s(b)-s(a)}, 
\end{equation}
where $\mathbb{P}_{y}(A)$ denotes the conditional probability of the
event $A$ under the condition that the process $\{Y_{t}\}_{t\ge
 0}$ starts at $y$.
\end{proposition}

With this in mind, we can now give the definition of a scale function.

\begin{definition}\label{def:scale-function}
 A continuous,
 strictly increasing function $s: E_{l,r}\to \R$ is called a
 \emph{scale function} of a continuous, regular, strong Markov process
 $\{Y_{t}\}_{t\ge 0}$ on the state space $E_{l,r}$ if for every $a$, $b$,
 $y\in E_{l,r}$ satisfying $l\le a<y<b\le r$, one has
 that~\eqref{eq:29} holds.
\end{definition}

We note that if $s: E_{l,r}\to \R$ is a scale function of the process
$\{Y_{t}\}_{t\ge 0}$, then it is not hard to see that for every
$\alpha$, $\beta \in \R$, the function $\tilde{s}: E_{l,r}\to \R$
given by the \emph{affine transformation}
 \begin{displaymath}
   \tilde{s}(y) =\alpha s(y)+\beta\qquad\text{for every $y\in E_{l,r}$}
 \end{displaymath}
 is another scale function of the process $\{Y_{t}\}_{t\ge 0}$. Thus,
 if there is a scale function $s$ of the process
 $\{Y_{t}\}_{t\ge 0}$, then there are infinitely many scale functions
 of $\{Y_{t}\}_{t\ge 0}$, and moreover, there is no loss of generality to
 assume that $s(0)=0$.\medskip

Now, if there is a scale function $s$ of the process
$\{Y_{t}\}_{t\ge 0}$ on $E_{l,r}$, then we can define the \emph{scaled process}
$\{Z_{t}\}_{t\ge0}$ on the state space $E_{s(l),s(r)}$ by setting
\begin{displaymath}
 Z_{t}:=s(Y_{t})\qquad\text{for every $t\ge 0$,}
\end{displaymath}
where $s(l)$ and $s(r)$ are the possibly improper limits
\begin{displaymath}
s(l):=\lim_{y\to l+}s(y)\quad\text{ and }\quad
s(r):=\lim_{y\to r-}s(y).
\end{displaymath}
Now, for every $a$, $b$, $z\in E_{s(l),s(r)}$ satisfying
$l\le a<y<b\le r$,  one has that 
\begin{displaymath}
  \mathbb{P}_{z}(\tau^b_{Z}<\tau^{a}_{Z})=\frac{z-a}{b-a}, 
\end{displaymath}
showing that the identity $s(z)=z$ is a scale function
of the scaled process $\{Z_{t}\}_{t\ge 0}$.

\begin{definition}
 One says that a process $\{X_{t}\}_{t\ge 0}$ on the state space $E_{l,r}$ is of
 \emph{natural scale} if the identity $s(x)=x$, ($x\in E_{l,r}$), is a scale function
 of the the process $\{X_{t}\}_{t\ge 0}$. 
\end{definition}

Processes of natural scale have the following characterization.

\begin{proposition}[{\cite[Chapter VII,  Proposition~3.5]{MR1725357}}]
  \label{prop:natural-scale}
 A continuous, regular, strong Markov process process
 $\{X_{t}\}_{t\ge 0}$ on state space $E_{l,r}$ is of natural scale if and only if 
 $\{X_{t}\}_{0\le t\le \tau^{l}_{X}\wedge\tau_{X}^{r}}$ is a local martingale.
\end{proposition}

From this proposition, we can conclude the following characterization
of a scale function $s$.

\begin{corollary}\label{cor:charact-scale}
 A continuous, strictly increasing function $s: E_{l,r}\to \R$ is called a
 \emph{scale function} of a continuous, regular, strong Markov process
 $\{Y_{t}\}_{t\ge 0}$ on the state space $E_{l,r}$ if and only if for
 the scaled process $\{Z_{t}\}_{t\ge 0}$, 
 \begin{equation}
   \label{eq:30}
 \{Z_{t}\}_{0\le t\le \tau^{s(l)}_{Z}\wedge\tau_{Z}^{s(r)}}\qquad\text{ is a
 local martingale.}
\end{equation}
\end{corollary}

With this preliminary, we can now come back to the example of the
Bessel process.

\begin{example}\label{ex:bessel-process-part1}
 For $\sigma\in (0,1)$, the $2(1-\sigma)$-Bessel
 process $\{Y_{t}\}_{t\ge 0}$ starting at $y=0$ for 
 is a continuous, regular, strong Markov process has state space
 $E_{0,\infty}=[0,\infty)$. By Theorem~\ref{thm:Donati-yor}, the
 $2\sigma^{\textrm{th}}$-powered process
 $\{Y^{2\sigma}_{t}\}_{t\ge 0}$ is sub-martingale with the Doob-Meyer
 decomposition $Y^{2\sigma}_{t}=M_{t}+L_{t}$, consisting of a
 continuous martingale $\{M_{t}\}_{t\ge0}$ and a continuous,
 non-decreasing process $\{L_{t}\}_{t\ge 0}$ carried by the zeros of
 $\{Y_{t}\}_{t\ge 0}$. Thus, for the function $s(y)=y^{2\sigma}$, the
 scaled process $\{Z_{t}\}_{t\ge 0}$ give by $Z_{t}=Y^{2\sigma}_{t}$,
 ($t\ge 0$), satisfies~\eqref{eq:30} and hence by
 Corollary~\ref{cor:charact-scale}, $\{Z_{t}\}_{t\ge 0}$ is of
 natural scale; or in other words, $s$ is a scale function of the
 $2(1-\sigma)$-Bessel process $\{Y_{t}\}_{t\ge 0}$.
\end{example}


\subsection{Speed measures}
 \label{subsec:speed-measures}
In this section, we intend to outline the impact of scale function
$s$ on speed measures $\hat{m}$ and the associated generalized
diffusion $\{Z_{t}\}_{t\ge 0}$ associated with
$m\in \mathfrak{m}_{\infty}$.\medskip

We built upon the scenario from
Section~\ref{subsec:scale-functions}. For given
$-\infty\le l<r\le \infty$, let $E_{l,r}$ be either a closed, open or
semi-closed interval in $\R$, and $\{Y_{t}\}_{t\ge 0}$ be a
continuous, regular, strong Markov process on the state space
$E_{l,r}$ with a killing time $\zeta$. Further, suppose,
$\{Y_{t}\}_{t\ge 0}$ can only be killed at the end-points $l$ and $r$
of $E_{l,r}$, provided they do not belong to $E_{l,r}$.\medskip

Let $I=(a,b)$ be an open interval such that the closure $[a,b]\subset
E_{l,r}$, and let $\sigma_I:=\inf\{t>0\,\vert\,Y_{t}\not\in I\}$ be
the \emph{first time $\{Y_{t}\}_{t\ge 0}$ exiting} the interval
$I$. Then, for the starting point $y$ of $\{Y_{t}\}_{t\ge 0}$, one has that
$y\in I$ yields that the first exit time $\sigma_{I}=\tau^{a}_Y\wedge
\tau_{X}^b$ almost surely, and $y\notin I$ implies that the first
exit time $\sigma_I=0$ almost surely. Further, let 
\begin{displaymath}
  m_I(y):=\mathbb{E}_{y}(\sigma_I)\qquad\text{for every $y\in \R$.}
\end{displaymath}

Now, let us take $J=(c,d)$ to be an open subinterval of $I$. Then by
the strong Markov property of $\{Y_{t}\}_{t\ge 0}$ and due to
property~\eqref{eq:29}, one obtains that for every scale function $s :
E_{l,r} \to \R$ of $\{Y_{t}\}_{t\ge 0}$, one has that
\begin{align*}
	m_{I}(y)
        &=m_J(y)+\mathbb{E}_{y}(\mathbb{E}_{Y_{\sigma_J}}(\sigma_I))\\
	&=m_J(y)+\frac{s(d)-s(y)}{s(d)-s(c)}m_I(c)+\frac{s(y)-s(c)}{s(d)-s(c)}m_I(d)
\end{align*}
for all $a<c<x<d<b$. 

Given a scale function $s :E_{l,r} \to \R$ of $\{Y_{t}\}_{t\ge 0}$ and an open interval $I=(a,b)$
with closure $[a,b]$ in $E_{l,r}$, let
\begin{displaymath}
	G_{s,I}(x,y):=
	\begin{cases}
	\displaystyle\frac{(s(x)-s(a))(s(b)-s(y))}{s(b)-s(a)} &\quad \text{if}\,\,a\le x\le y\le b,\\[9pt]
	\displaystyle\frac{(s(y)-s(a))(s(b)-s(x))}{s(b)-s(a)} &\quad \text{if}\,\,a\le x\le y\le b,\\
		0	&\quad \text{otherwise}.
	\end{cases}
\end{displaymath}

With this preliminary, we can now recall the following existence
theorem of the \emph{speed measure} $\hat{m}$. 

\begin{theorem}[{\cite[Chapter VII, Proposition~3.5]{MR1725357}}]
 Let $\{Y_{t}\}_{t\ge 0}$ be a
continuous, regular, strong Markov process on the state space
$E_{l,r}$. Then, there exists a unique Radon
 measure $\mu_{\hat{m}}$ on the interior $\mathring{E}_{l,r}$ of $E_{l,r}$
 and a scale function $s :E_{l,r} \to \R$ of $\{Y_{t}\}_{t\ge 0}$
 such that for every open sub-interval $I=(a,b)$ of $E$, one has
 \begin{equation}
   \label{eq:31}
   m_I(y)=\int_{I}G_{s,I}(y,r)\,\td\mu_{\hat{m}}(r)\qquad\text{for every $y\in I$.}
 \end{equation}
\end{theorem}

The above existence theorem leads to the following definition.

\begin{definition}
 For a continuous regular, strong Markov process $\{Y_{t}\}_{t\ge 0}$
 on the state space $E_{l,r}$, one calls the Radon measure $\hat{m}$
 satisfying~\eqref{eq:31} the \emph{speed measure} of the process
 $\{Y_{t}\}_{t\ge 0}$.
\end{definition}

For our next theorem, we recall the following elementary definition.

\begin{definition}
 For a strictly increasing continuous function $s : (a,b)\to \R$
 defined on an open interval $(a,b)$ with range $\Rg(s)\subseteq \R$,
 a function $f : \Rg(s)\to \R$ is called \emph{$s$-differentiable at
   $y\in \Rg(s)$} if the limit
 \begin{displaymath}
   \lim_{\hat{y}\to y}\frac{f(\hat{y})-f(y)}{s(\hat{y})-s(y)}
 \end{displaymath}
 exists and then, we denote by
 \begin{displaymath}
   \frac{\td f}{\ds}(y):=\lim_{\hat{y}\to y}\frac{f(\hat{y})-f(y)}{s(\hat{y})-s(y)}
 \end{displaymath}
 the \emph{$s$-derivative of $f$ at $y$.}
\end{definition}

\begin{remark}
 Recall, every continuous, strictly increasing function
 $s : (a,b)\to \R$ is a.e. differentiable in $(a,b)$. Thus, if $f$ is
 differentiable at $y\in \Rg(s)$ and $s'(y)>0$ exists, then 
 the $s$-derivative of $f$
 \begin{equation}
   \label{eq:33}
   \frac{\td f}{\ds}(y)=\frac{\td f}{\dy}(y)\frac{1}{s'(y)}.
 \end{equation}
\end{remark}

Now, our next theorem outlines the relation between the speed measure
$\hat{m}$ and \emph{scale function} $s : E_{l,r} \to \R$ in sense of
infinitesimal generator of the process $\{Y_{t}\}_{t\ge 0}$. Here, we
denote by $A_{Y}$ the infinitesimal generator of $\{Y_{t}\}_{t\ge 0}$ and by $D(A_{Y})$
its domain.

\begin{theorem}[{\cite[Chapter VII, Theorem~3.12 and
   Proposition~3.13]{MR1725357}}]\label{thm:characterization-of-B}
 Let $\{Y_{t}\}_{t\ge 0}$ be a continuous, regular, strong Markov
 process on the state space $E_{0,r}$, ($0<r\le \infty$), with infinitesimal generator
 $B$. Then, there is a scale function $s : E_{0,r}\to \R$ of
 $\{Y_{t}\}_{t\ge 0}$ such that
 \begin{displaymath}
   Bf=\frac{\td}{\td\hat{m}}\frac{\td f}{\ds}\qquad\text{for every
     $f\in D(B)$,}
 \end{displaymath}
where $\hat{m}$ is the speed measure of $\{Y_{t}\}_{t\ge 0}$. Moreover,
 if $E_{0,r}=[0,r)$, then one has that
   \begin{displaymath}
     \hat{m}(0+)Bf(0)-\frac{\td f}{\ds}_{+}(0)=0\qquad\text{for
       every $f\in D(B)$.}
   \end{displaymath}
\end{theorem}

From this theorem, one see that the pair $(s,\hat{m})$ of the scale
function $s$ and the speed measure $\hat{m}$ given by
Theorem~\ref{thm:characterization-of-B} characterizes the infinitesimal
generator $B$ of the continuous, regular, strong Markov
process $\{Y_{t}\}_{t\ge 0}$.\medskip 

Let us consider the $2(1-\sigma)$-Bessel process $\{Y_{t}\}_{t\ge 0}$ as an example.

\begin{example}\label{ex-scale-of-Bessel-process}
 We recall (cf.,~\cite[Sect. 1.2]{MR2417969}) that for
 $\sigma\in (0,1)$, the infinitesimal generator $B_{Y}$ of the
 $2(1-\sigma)$-Bessel process $\{Y_{t}\}_{t\ge 0}$ on
 $[0,\infty)$ is given by
 \begin{align*}
   B_{Y}u&=\frac{1}{2}\frac{\td^2u}{\td
          y^2}+\frac{1-2\sigma}{2y}\frac{\td u}{\dy}\\
   &=\frac{1}{2y^{1-2\sigma}}
     \frac{\td}{\dy}\Big(y^{1-2\sigma}\frac{\td u}{\dy}\Big)\\
   &=\frac{1}{2}\Big
     (\frac{y}{2\sigma}\Big)^{2\sigma-1}
     \frac{\td}{\dy}\Big(\Big(\frac{y}{2\sigma}\Big)^{1-2\sigma}\frac{\td
     u}{\dy}\Big).
 \end{align*}
 From the above computation, it is natural that the scale function
 $s$ and speed measure $\mu_{\hat{m}_{Y}}$ are given by
 \begin{displaymath}
   s(y)=\Big (\frac{y}{2\sigma}\Big )^{2\sigma}\qquad\text{and}\qquad 
   \td\mu_{\hat{m}_{Y}}(y)={2}\Big (\frac{y}{2\sigma}\Big )^{1-2\sigma}\dy, 
 \end{displaymath}
 and so, by~\eqref{eq:33}, one sees that
 \begin{displaymath}
   B_{Y}u=\frac{1}{2}\Big
   (\frac{y}{2\sigma}\Big)^{2\sigma-1}\frac{\td}{\dy}\frac{\td
     u}{\ds}
   =\frac{\td}{\td\hat{m}_{Y}}\frac{\td u}{\ds}.
 \end{displaymath}
 
 Next, if we apply the substitution
 \begin{displaymath}
   z=s(y)=\Big (\frac{y}{2\sigma}\Big )^{2\sigma}
 \end{displaymath}
 to the process $\{Y_{t}\}_{t\ge 0}$; that is,
 one sets $Z_{t}=s(Y_{t})$, then the process $\{Z_{t}\}_{t\ge 0}$ is
 the (up to a multiple constant) $2\sigma^{\textrm{th}}$-power process
 (cf., Theorem~\ref{thm:Donati-yor}), which is of natural scale (cf.,
 Example~\ref{ex:bessel-process-part1}) with speed measure
 $\mu_{\hat{m}_{Z}}$ given by
 \begin{displaymath}
   \td\mu_{\hat{m}_{Z}}(z)=2z^{-\frac{2\sigma-1}{\sigma}}\dz. 
 \end{displaymath}
 Hence, the infinitesimal generator $B_{Z}$ of $\{Z_{t}\}_{t\ge 0}$
 is given by
 \begin{displaymath}
   B_{Z}=\frac{1}{2}z^{\frac{2\sigma-1}{\sigma}}\frac{\td^2u}{\td
     z^2}=\frac{1}{2z^{\frac{1-2\sigma}{\sigma}}}\frac{\td}{\dz}\frac{\td
   u}{\dz}=\frac{\td}{\td \hat{m}_{Z}}\frac{\td u}{\dz}.
\end{displaymath}
On the other hand, the process $\{Z_{t}\}_{t\ge 0}$ is a generalized
diffusion associated with the string $m$ given by
\begin{displaymath}
 m(z)=
 \begin{cases}
   \displaystyle\int_{0}^{z} \td\mu_{\hat{m}_{Z}}(s)=
   \frac{2\sigma}{1-\sigma}z^{\frac{1-\sigma}{\sigma}} &\text{if $z\ge 0$,}\\[7pt]
   0 &\text{otherwise.}
 \end{cases}
\end{displaymath}
\end{example}

In the context of local times and generalized diffusions, it is natural
to ask what the effect is for a given continuous, regular, strong Markov process
$\{Z_{t}\}_{t\ge 0}$ associated with a given string
$m\in \mathfrak{m}_{\infty}$, if one switches from \emph{natural scale} to the
\emph{original scale}, that is, from $\{Z_{t}\}_{t\ge 0}$ to
$\{Y_{t}\}_{t\ge 0}$ by setting
\begin{equation}
 \label{eq:34}
 Y_t=s^{-1}(Z_t)\qquad\text{for all $t\ge 0$,}
\end{equation}
for a scale function $s : E_{l,r}\to \R$ of $\{Y_{t}\}_{t\ge 0}$.
To see this, we employ the occupation times
formula~\eqref{eq:20} involving the local time process
$\{\tilde{L}_t(z)\}_{t\ge 0}$ of $\{Z_{t}\}_{t\ge 0}$ at level $z$. Then,
\begin{displaymath}
 \int_{0}^{t}g(s^{-1}(Z_{r}))\,\dr=\int_{E_{m}}g(s^{-1}(z))\tilde{L}_t(z)\,\td
 \mu_m(z)
\end{displaymath}
for every $g\in L^{\infty}(E_{m})$. Applying the change of variable
$z=s(y)$ to the integral on the right hand-side of the last equation, yields that
\begin{equation}
 \label{eq:35}
\int_{0}^{t}g(Y_{r})\,\dr=\int_{s^{-1}(E_{m})}g(y)\tilde{L}_{t}(s(y))\,\td
s^{-1}_{\#}\mu_m(y),
\end{equation}
where $s^{-1}_{\#}\mu_{m}$ is the push-forward measure of $m$ by $s^{-1}$.
Since $s(0)=0$, the rescaling~\eqref{eq:34} with the scale function
$s$ does not affect the local time $\tilde{L}_t(0)$ of
$\{Z_{t}\}_{t\ge 0}$ at level $0$. Moreover, the process
$\{\tilde{L}_{t}(s(y))\}$ is a local time process of
$\{Y_{t}\}_{t\ge 0}$ at level $s(y)$. Thus, the occupation times
formula~\eqref{eq:35} suggests that the speed measure $\hat{m}_{Y}$
associated with the process $\{Y_{t}\}_{t\ge 0}$ is the push-forward
measure $s^{-1}_{\#}m$, which certainly is the case (cf.,~\cite[see
Chapter VII, Exercise 3.18]{MR1725357}).\medskip

In Table~\ref{table:1}, we summarize the relation between a scale function
$s$, the speed measure $\mu_{\hat{m}_{Y}}$ associated with $\{Y_{t}\}_{t\ge 0}$,
and the \emph{natural scaled} process $\{Z_{t}\}_{t\ge 0}$
via~\eqref{eq:34}, and its impact to the Dirichlet problem and
corresponding Dirichlet-to-Wentzell operator (DtW operator).

\begin{center}
 \begin{table}
   \begin{tabular}{|c|c|c|cc}\cline{1-3}
    \emph{string} $m$ &  \emph{original scale} $s$ &
                                                           \emph{natural scale} $\textrm{id}_{\R}$	 &  &  \\ \cline{1-3}
     {\small\em  speed measure}	& 	$\mu_{\hat{m}_{Y}}=s^{-1}_{\#}\mu_m$	 &	$\mu_{\hat{m}_{Z}}=\mu_m$&  &  \\ \cline{1-3}
     {\small\em  generator}	& $-\frac{1}{2}\frac{\td}{\td\hat{m}_{Y}}\frac{\td}{\ds}$ &
                                                                    $-\frac{1}{2}\frac{\td}{\td m}\frac{\td}{\dz}$ &  &  \\ \cline{1-3}
     {\small \em Dirichlet prb.  }&
                               $Au-\frac{1}{2}\frac{\td}{\td\hat{m}_{Y}}\frac{\td}{\ds}u=0$
                                                   &
                                                     $Au-\frac{1}{2}\frac{\td}{\td m}\frac{\td}{\dz}u=0$
                                                &  &  \\ \cline{1-3}
     {\small\em   DtW operator } &  {\small  $u(0)\mapsto
                           m(0+)Au(0)-\frac{1}{2}\frac{\td u}{\ds}_{+}(0)$}
                      &  {\small $u(0)\mapsto m(0+)Au(0)-\frac{1}{2}\frac{\td
                        u}{\dz}_{+}(0)$}
                                                &  &  \\ \cline{1-3}
   \end{tabular}
   \mbox{}
   \caption{Scale functions, speed measures,
     DtW operator}
   \label{table:1}
 \end{table}
\end{center}

We note that the generator is slightly different than what is proposed
in Theorem \ref{thm:characterization-of-B} in the sense that we add a
factor of $\frac{1}{2}$. The reason for this is that in Section
\ref{SEC:PRIMER_STRINGS} we obtain a \emph{generalized diffusion} by
considering a time change of reflecting Brownian motion. The generator
of a reflecting Brownian motion is $\frac{1}{2}\frac{d^2}{dx^2}$
however the generator of the generalized diffusion is
$\frac{1}{2}\frac{d}{dm}\frac{d}{dx}$ hence for us it is convenient to
keep a factor of $\frac{1}{2}$.

We conclude this section with the following remark.

\begin{remark}\label{rem:3}
 It is worth noting that for every string $m\in \mathfrak{m}_{\infty}$,
 the generalized diffusion process $\{Z_{t}\}_{t\ge 0}$ associated
 with $m$, is already in \emph{natural scale} $s=2\,\textrm{id}_{\R}$ and the speed
 measure $\hat{m}_{Z}$ coincides with the associated measure
 $\mu_{m}$ of $m$. This follows immediately from the fact that the infinitesimal
 generator $B$ of $\{Z_{t}\}_{t\ge 0}$ is given by~\eqref{eq:32}.
\end{remark}

%
%
%
%

\subsection{Hitting times and spectral representations}
\label{subsec:hitting-times}
Hitting time densities admit the useful property of a spectral
decomposition (see, for instance, the papers \cite{MR637387},
\cite{MR1022933} or \cite{MR661628}, the survey paper \cite[Section
3]{MR2540855}, and~\cite{MR1997032} containing interesting
examples). We intend to use the spectral decomposition of hitting time
densities to directly calculate the Dirichlet-to-Wentzell operator
$\Lambda_m$.\medskip

For a given string $m\in \mathfrak{m}_{\infty}$, let $\{Z_{t}\}_{t\ge 0}$ be
a generalized diffusion process associated with $m$, and 
\begin{equation}
  \label{eq:46}
\tau:=\inf\Big\{t>0\,\Big\vert\,Z_{t}=0\Big\} 
\end{equation}
be the \emph{first hitting time} of $z=0$ by
$\{Z_{t}\}_{t\ge0}$. Then through $\tau$, we can define the
\emph{killed process} $\{\hat{Z}_t\}_{t\geq 0}$ by setting
\begin{displaymath}
\hat{Z}_t:=
\begin{cases}
Z_t & t\le \tau, \\
\partial & t> \tau, \\
\end{cases}
\end{displaymath}
where $\partial$ denotes the \emph{coffin state}.

\begin{remark}\label{rem:infinitesimal-generator-of-killed-process}
  We note that the killed process $\{\hat{Z}_{t}\}_{t\ge 0}$ has the
  same speed measure $m$ and the same scale function
  $s=2\,\textrm{id}_{\R}$ as $\{Z_{t}\}_{t\ge 0}$ (cf.,
  Remark~\ref{rem:3}). Moreover, the \emph{infinitesimal generator}
  $-\hat{B}_{m}$ of the killed process $\{\hat{Z}_{t}\}_{t\ge 0}$ is
  the $2$nd-order differential operator $-B_{m}$ defined
  in~\eqref{eq:32} equipped with the boundary conditions $f(0)=0$.
\end{remark}

Further, the killed process $\{\hat{Z}_{t}\}_{t\ge 0}$ has a \emph{transition density} 
\begin{displaymath}
    \hat{p}(t,z,y):=\frac{\mathbb{P}_{z}(\hat{Z}_{t}\in \dy)}{\td
      \mu_m}
    =\frac{\mathbb{P}_{z}(Z_{t}\in\dy\,|\,t<\tau) }{\td \mu_m}
\end{displaymath}
for every $t\ge 0$, $z$, $y\in [0,\infty)$. In the next proposition, we
recall some important properties from \cite[Section~4.11]{MR0345224}
of the transition density $\hat{p}(t,z,y)$.

\begin{proposition}[{\cite[Section~4.11]{MR0345224} \& \cite{MR1022933}}]
  \label{prop:densities-p-and-omega}
  For a given string $m\in \mathfrak{m}_{\infty}$, let $\{Z_{t}\}_{t\ge 0}$ be
the generalized diffusion process with generator
$\tfrac{1}{2}\tfrac{\td}{\td m}\tfrac{\td}{\dz}$. Then the transition density
  $\hat{p}(t,z,y)$ of the killed process $\{\hat{Z}_{t}\}_{t\ge 0}$ is
  a  jointly continuous function $\hat{p} : [0,\infty)\times [0,\infty) \times [0,\infty)\to
  [0,\infty)$ with the following properties.
  \begin{enumerate}
    \item $\hat{p}(t,z,0)=\hat{p}(t,0,y)=0$ for every $t\ge 0$, $z$, $y\in [0,\infty)$.
  \item $\hat{p}$ is symmetric; that is,
  \begin{displaymath}
    \hat{p}(t,z,y)=\hat{p}(t,y,z)\qquad\text{for all
      $t\in [0,\infty)$ and $y$, $z\in (0,\infty)$.}
  \end{displaymath}
\end{enumerate}
\end{proposition}

The next definition contains the core density of our \emph{Poisson
formula}~\eqref{eq:17}.

\begin{definition}\label{def:proba-density-of-tau}
  For a given string $m\in \mathfrak{m}_{\infty}$, let $\tau$ be the first
  hitting time of $z=0$ of a generalized diffusion process
  $\{Z_{t}\}_{t\ge 0}$ associated with $m$. Then $\tau$ has the
  \emph{probability density}
  \begin{displaymath}
    \omega_{\tau}(t,z):=\frac{\mathbb{P}_{z}(\tau \in
      \dt)}{\dt}\qquad\text{for every $z\in (0,\infty)$ and $t>0$.}
  \end{displaymath}
\end{definition}

  To derive the spectral representation of the first hitting time
  density $\omega_{\tau}$, we begin by recalling
  from~\cite[Section~5]{MR637387} the so-called \emph{eigenfunction
    expansion} $\{C(\cdot,\gamma)\}_{\gamma>0}$ of the operator
  $\hat{B}_{m}$, which is directly related to the killed process
  $\{\hat{Z}_{t}\}_{t\ge 0}$
  (cf.~Remark~\ref{rem:infinitesimal-generator-of-killed-process}). To
  do so, for given $\gamma>0$, let $C(\cdot,\gamma)$ be the unique
  solution of
\begin{equation}
  \label{eq:47}
  \begin{cases}
   \displaystyle -\frac{1}{2}\frac{\td}{\dm}\frac{\td}{\dz}C(z,\gamma)=\gamma\,C(z,\gamma)
    &\qquad\text{for $z\in (0,\infty)$,}\\[7pt]
     C(0;\gamma)=0, \quad \lim_{z\rightarrow 0^{+}}\frac{C(z;\gamma)}{2z}=1.&
  \end{cases}
\end{equation}
Alternatively, for every $\gamma>0$, $C(\cdot,\gamma)$ satisfies 
\begin{equation}\label{eq:86}
  C(z;\gamma)=2z-2\gamma \int_{0}^{z}\int_{(0,x]}C(r;\gamma)\,\td \mu_m(r)\,\dx 
\end{equation}
for every $z\in [0,\infty)$. Next, let $\{C_{n}\}_{n\ge 0}$ be a
sequence of functions $C_{n} : [0,\infty)\to \R$ recursively defined by
$C_0(z):=2z$, and for every $n\ge 1$, 
\begin{equation}
  \label{eq:36}
  C_{n}(z):=2\int_{0}^{z}\int_{(0,x]}C_{n-1}(r)\,\td
  \mu_m(r)\,\dx,\qquad(z\in [0,\infty)). 
\end{equation}
For this sequence $\{C_{n}\}_{n\ge 0}$ (cf.,~\cite[Section
5.4]{MR0448523} for a similar construction), the function
$C(\cdot,\gamma)$ has the following \emph{series representation}
\begin{equation}
  \label{eq:38}
  C(z,\gamma)=\sum_{n=0}^{\infty}(-\gamma)^{n}C_{n}(z)
\end{equation}
for every $z\in [0,\infty)$ and $\gamma>0$. Note, the series on the
right hand-side in~\eqref{eq:38} converges locally uniformly on
$[0,\infty)$. This follows from the Weierstrass M-test, and the
estimates~\eqref{eq:39} in the subsequent lemma. The next statements will be
useful later in the proofs of our theorems to verify the hypothesis of
dominated convergence.

\begin{lemma}[{\cite[Lemma 3.1 and 3.2]{MR2540855}}]
  \label{lem:1}
  For every $n\ge 0$, the functions $C_{n}$ defined in~\eqref{eq:36}
  are positive and increasing along $[0,\infty)$, and satisfy
  \begin{equation}
    \label{eq:39}
   C_{n}(z)\le \frac{2z}{n!}\Big ( 2\int_{0}^{z}M(x)\,\dx\Big)^{n}.
 \end{equation}
 where $M(x):=m(x)-m(0)$ for $x>0$. Moreover, for the unique solution
$C(\cdot,\gamma)$ of~\eqref{eq:48},one has that
 \begin{equation}
   \label{eq:40}
	\abs{C(z,\gamma)}\le 2z\,e^{\gamma 2\int_{0}^{z}M(x)\dx}
  \end{equation}
  for every $\gamma>0$ and $z\in [0,\infty)$.    
\end{lemma}

\begin{proof}
  It follows from~\eqref{eq:36} and since $C_{0}(z)=2z$ that each
  $C_{n}$ is positive and increasing along $[0,r_{m})$. Next, we show
  that inequality~\eqref{eq:39} by an induction over $n\ge
  0$. Obviously,~\eqref{eq:39} is satisfies by $C_{0}(z)=2z$ and hence,
  the case $n=0$ holds. Now, suppose that~\eqref{eq:39} holds for an
  integer $n\ge 1$. Then, it remains to show that~\eqref{eq:39} holds
  for $n+1$. But by induction hypothesis, one sees that
  \begin{align*}
    C_{n+1}(z)	&=2\int_{0}^{z}\int_{(0,x]}C_{n}(r)\,\td \mu_m(r)\,\dx \\
	&\le 2\int_{0}^{z}\int_{(0,x]}\frac{r}{n!}\Big(2
          \int_{0}^{r}M(u)\,\du\Big)^{n}\,\td \mu_m(r)\,\dx \\
	&\le \frac{2z}{n!} \int_{0}^{z}\int_{(0,x]}\Big(
   2\int_{0}^{r}M(u)\,\du\Big)^{n}\,\td \mu_m(r)\,\dx \\
                &\le \frac{2z}{n!} \int_{0}^{z}\int_{(0,x]}\Big(
   2\int_{0}^{x}M(u)\,\du\Big)^{n}\,\td \mu_m(r)\,dx\\
	&\le \frac{2z}{n!}\int_{0}^{z}\Big(2\int_{0}^{x}M(u)\,\du\Big)^{n}M(x)\,\dx\\
	&=\frac{2z}{(n+1)!}\Big(2\int_{0}^{z}M(x)\dx\Big)^{n+1}.
  \end{align*}
  This proves~\eqref{eq:39} for $n+1$ and hence, this inequality holds
  for all integers $n\ge 0$. One sees that estimate~\eqref{eq:40} holds
  after applying~\eqref{eq:39} to the series representation~\eqref{eq:38}. 
  This completes the proof of this lemma.
\end{proof}

By~\cite[Theorem 2.3]{MR1022933}, for the killed process
$\{\hat{Z}_{t}\}_{t\ge 0}$, there exists a $\sigma$-finite measure
$\hat{\Delta}_{m}$, called the \emph{principal measure} of
$\{\hat{Z}_{t}\}_{t\ge 0}$, for which the following holds
\begin{equation}
  \label{eq:42}
  \int_{0}^{\infty}\frac{1}{\gamma(\gamma+1)}\,\td\hat{\Delta}_{m}(\gamma)<\infty
\end{equation}
and
\begin{equation}
  \label{eq:43}	
  \int_{0}^{\infty}\frac{1}{\gamma}\,\td\hat{\Delta}_{m}(\gamma)=\infty.
\end{equation}
From this, one can deduce the following \emph{integral representation}
\begin{equation}
  \label{eq:44}
  \hat{p}(t,z,y)=\int_{0}^{\infty} e^{-\gamma t}C(z,\gamma)
  C(y,\gamma)\,\td\hat{\Delta}_{m}(\gamma),\quad
  (t\ge 0,\, z,\,y\in [0,\infty)),
\end{equation}
of the transition density $\hat{p}$ of $\{\hat{Z}_{t}\}_{t\ge
  0}$. Further, the
following \emph{spectral representation} of the hitting time density
$\omega_{\tau}$ and \emph{representation of the density} $h_{m}$ of
the L\'evy measure $\nu_{m}$ associated with
$m$ (cf, Definition~\ref{def:1}) is available.

  \begin{theorem}[{\cite[Theorem 3.1 and 3.2]{MR1022933}}]
    \label{DTN:Levy Measure}
    For a given string $m\in \mathfrak{m}_{\infty}$, let $\tau$ be the first
  hitting time of $z=0$ by a generalized diffusion process
  $\{Z_{t}\}_{t\ge 0}$  associated with $m$, and $\omega_{\tau}$ the probability
    density of $\tau$. Further, for given $\gamma>0$, let $C(\cdot,\gamma)$ be the
    unique solution of~\eqref{eq:48}. Then, the following statements
    hold.
    \begin{enumerate}
    \item The probability density $\omega_{\tau}$ has the spectral representation
      \begin{equation}
        \label{eq:37}
        \omega_{\tau}(t,z)=\int_{0}^{\infty} e^{-\gamma
          t}C(z,\gamma)\,\td\hat{\Delta}_{m}(\gamma)
      \end{equation}
      for every $t>0$, $z\in (0,\infty)$.
    \item One has that
      \begin{equation}
        \label{eq:45}
        \omega_{\tau}(t,z)=
        \lim_{y\rightarrow 0+}\frac{\hat{p}(t,z,y)}{2y}\qquad\text{for
          every $z\in (0,\infty)$.}
      \end{equation}
      \item The density
      $h_{m}=\frac{d\nu_{m}}{dt}$ of the L\'evy measure $\nu_{m}$
      associated with $m$, one has that
      \begin{equation}
        \label{eq:125}
        h_{m}(t)
        = \lim_{z\rightarrow 0^{+}}\frac{\omega_{\tau}(t,z)}{2z}
	=\int_{0}^{\infty} e^{-\gamma t}\,\td\hat{\Delta}_{m}(\gamma)\qquad\text{for every $t>0$.}
      \end{equation}
    \end{enumerate}
\end{theorem}


   
%

Integrating formula~\eqref{eq:37} of the first hitting density
$\omega_{\tau}$ of $\tau$ over $(t,\infty)$ for $t>0$, then
by~\eqref{eq:40} and Fubini's theorem yields the following.

\begin{corollary}[{\cite[Proposition~3.7]{MR2540855}}]\label{prop:1}
  For a given string $m\in \mathfrak{m}_{\infty}$ on $\R$, let $\tau$ be the first
  hitting time~\eqref{eq:46} of $z=0$ by a generalized
  diffusion $\{Z_{t}\}_{t\ge 0}$ associated with $m$, and for $\gamma>0$, let
  $C(\cdot,\gamma)$ be the unique solution of~\eqref{eq:48}. Then for
  every $z\in (0,\infty)$, one has that
  \begin{displaymath}
    \mathbb{P}_{z}(\tau>t)=\int_{0}^{\infty}\frac{e^{-\gamma
      t}}{\gamma}C(z,\gamma)\,\td\hat{\Delta}_{m}(\gamma)\qquad\text{for all $t\in [0,\infty)$.} 
\end{displaymath}
\end{corollary}

In our next proposition, we collect some important properties of an
anti-derivative $\beta_{m}$ of the density $h_{m}$ of $\nu_{m}$, which
we require later in Section~\ref{subsec:existence} to derive an
\emph{integration by parts} argument (see
Lemma~\ref{lem:ibp}). 

\begin{proposition}
  \label{prop:well-defined-beta-m}
  For a given string $m\in \mathfrak{m}_{\infty}$ on $\R$, let
  $\hat{\Delta}_{m}$ be the principal measure associated with $m$, and
  set
  \begin{equation}\label{eq:118}
    \beta_{m}(t):=\int_0^\infty \frac{1}{\gamma}e^{-\gamma t}\,
    \td\hat{\Delta}_{m}(\gamma) \qquad \text{for every $t>0$.}
  \end{equation}
  Then, the following properties hold.
  \begin{enumerate}
  \item The function $\beta_m: (0,\infty)\to [0,\infty)$ is differentiable
    with derivative
    \begin{equation}
      \label{Cor:int_by_prts}
	\frac{\td\beta_m}{\td t}(t)=h_m(t)\qquad \text{for every $t>0$}.
      \end{equation}
    \item One has that
      \begin{equation}
        \label{COR:LIMITS_INT}
        \lim_{t\rightarrow 0^+}t\beta_m(t)=\lim_{t\rightarrow \infty}\beta_m(t)=0. 
      \end{equation}
  \end{enumerate}
\end{proposition}

\begin{proof}
  By applying the elementary inequality
  \begin{equation}
    \label{eq:126}
    (\gamma+1)e^{-\gamma t}\le
    \frac{1}{t}e^{t}\qquad\text{holding for all $t>0$,} 
  \end{equation}
  to the integrand in~\eqref{eq:118} shows that
  \begin{displaymath}
    \beta_m(t)\le\frac{1}{t}e^{t}
    \int_{0}^{\infty}\frac{1}{\gamma(1+\gamma)}\,
    \td\hat{\Delta}_{m}(\gamma).
  \end{displaymath}
  Since the integral on right hand-side is finite (see \eqref{eq:42}),
  $\beta_{m}$ is a finite positive function on $(0,\infty)$. To see
  that $\beta_{m}$ is differentiable at $t>0$, note that the integrand
  $f : (0,\infty)\times (0,\infty)\to (0,\infty)$ given by
  \begin{displaymath}
    f(\gamma,t)=\frac{1}{\gamma}e^{-\gamma t}\qquad\text{for every
      $(\gamma,t)\in (0,\infty)\times (0,\infty)$,}
  \end{displaymath}
  is continuously differentiable in $t$ with partial derivative
  $\frac{\partial f}{\partial t}(\gamma,t)=- e^{-\gamma t}$. Moreover,
  since $e^{-\gamma t}\le  e^{-\gamma \varepsilon}$ for any
  $0<\varepsilon<t$ and since
  \begin{displaymath}
    0\le e^{-\gamma \varepsilon}\lesssim
    \frac{1}{\gamma(\gamma+1)}\qquad\text{for every $\gamma\in (0,\infty)$,}
  \end{displaymath}
  it follows from~\eqref{eq:42} and Lebesgue's dominated
  convergence theorem that the parameter integral
  \begin{displaymath}
    \beta_{m}(t)=\int_{0}^{\infty}f(\gamma,t)\,\td\hat{\Delta}_{m}(\gamma)
  \end{displaymath}
  is differentiable and $\frac{\td\beta_m}{\td t}(t)=h_m(t)$ with
  $h_{m}$ given by~\eqref{Cor:int_by_prts}. Further, since for every
  $\gamma>0$, $\lim_{t\to \infty}f(\gamma,t)=0$, and since 
  \begin{displaymath}
    0\le f(\gamma, t)\le f(\gamma, M)\qquad\text{for every $t>M>0$,}
  \end{displaymath}
   we can conclude again from Lebesgue's dominated convergence theorem
   that $\lim_{t\rightarrow \infty}\beta_m(t)=0$. Finally, by using
   again~\eqref{eq:126}, one sees that
   \begin{displaymath}
     0\le t\,f(\gamma,t)\le \frac{e^{t}}{\gamma(\gamma+1)}\qquad\text{for every
      $(\gamma,t)\in (0,\infty)\times (0,\infty)$,}
   \end{displaymath}
   and since $\lim_{t\to 0+}t\, f(\gamma,t)=0$, it follows
   from~\eqref{eq:42} and Lebesgue's dominated convergence theorem
   that $\lim_{t\rightarrow 0+}t\beta_m(t)=0$. This completes the
   proof of~\eqref{COR:LIMITS_INT} and thereby the proof of this proposition.
\end{proof}

The following proposition provides a crucial identity to prove existence of a weak
solution to the extension problem~\eqref{eq:16} (see Section~\ref{subsec:existence}). 

\begin{proposition}\label{LEMMA:dou_int}
  For a given string $m\in \mathfrak{m}_{\infty}$ on $\R$, let
  $\omega_{\tau}$ be the probability density~\eqref{eq:45} of the
  first hitting time $\tau$ of $z=0$ by a generalized diffusion
  process $\{Z_{t}\}_{t\ge 0}$ associated with $m$, and $\beta_{m}$
  given by~\eqref{eq:118}. Then, one has that
  \begin{equation}
    \label{eq:129}
	2\,\int_{0}^{z}\int_{(0,x]}\omega_{\tau}(t,v)\,\td \mu_m(v)\dx
        =2\,z\, \beta_m(t) 
        -\mathbb{P}_{z}(\tau>t)
  \end{equation}
  for every $z\in (0,\infty)$ and $t>0$.
\end{proposition}

\begin{proof}
  Let $t>0$ and $z\in (0,\infty)$. According to the spectral representation~\eqref{eq:37}
  of the first hitting time density $\omega_{\tau}$ (Theorem~\ref{DTN:Levy
    Measure}), one sees that
  \begin{equation}\label{eq:128}
    \begin{split}
      &\int_{0}^{z}\int_{(0,x]}\omega_{\tau}(t,v)\,\td\mu_m(v)\,\dx\\
      &\hspace{2cm}
        =\int_{0}^{z}\int_{(0,x]}\int_{0}^{\infty}e^{-\gamma
          t}C(v,\gamma)\,
        \td\hat{\Delta}_{m}(\gamma)\,\td\mu_m(v)\,\dx\\
      &\hspace{2cm}
        = \int_{\Omega_z}\int_{0}^{\infty}e^{-\gamma
        t}C(v,\gamma)\,\td\hat{\Delta}_{m}(\gamma)\,
        \td(\mu_{\text{Leb}}\otimes \mu_m)(x,v).
    \end{split}
  \end{equation}
  where we set $\Omega_z=\{(x,v)\in [0,\infty)^2:0\leq v\le x\le z\}$ and
  $\mu_{\text{Leb}}\otimes \mu_m$ is the product measure of
  the one-dimensional Lebesgue measure $\mu_{\text{Leb}}$ and $\mu_m$.
  Further, if $\hat{p}$ is the
  transition density of the killed process $\{\hat{Z}_{t}\}_{t\ge 0}$
  of the generalized diffusion $\{Z_{t}\}_{t\ge 0}$ and
  $h_{m}$ the density $h_{m}$ of
  the L\'evy measure $\nu_{m}$ associated with $m$, then by Cauchy-Schwarz's
  inequality, by \eqref{eq:44}, and by \eqref{eq:125}, one sees that 
  \begin{align*}
    &\int_{0}^{\infty}e^{-\gamma t}|C(v,\gamma)|\,\td\hat{\Delta}_{m}(\gamma)\\
    &\qquad \le
      \left( \int_{0}^{\infty}e^{-\gamma
      t}\,\td\hat{\Delta}_{m}(\gamma)
      \right)^{\frac{1}{2}}\,\left(
      \int_{0}^{\infty}e^{-\gamma t}|C(v,\gamma)|^2\,\td\hat{\Delta}_{m}(\gamma)
      \right)^{\frac{1}{2}}\\
    &\qquad=\sqrt{h_m(t)}\,\sqrt{\hat{p}(t,v,v)}
  \end{align*}
  and so,
  \begin{align*}
     &\int_{\Omega_z}\int_{0}^{\infty}e^{-\gamma
       t}|C(v,\gamma)|\,\td\hat{\Delta}_{m}(\gamma)\,
       \td(\mu_{\text{Leb}}\otimes \mu_m)(x,v)\\
     &\qquad \le \int_{\Omega_z}\sqrt{h_m(t)}\,
       \sqrt{\hat{p}(t,v,v)}\, \td(\mu_{\text{Leb}}\otimes \mu_m)(x,v)\\
	&\qquad= \sqrt{h_m(t)}\,\int_{\Omega_z}\sqrt{\hat{p}(t,v,v)}\,
   \td(\mu_{\text{Leb}}\otimes \mu_m)(x,v)\\
    &\qquad=
      \sqrt{h_m(t)}\,\int_{0}^{z}\int_{(0,x]}\sqrt{\hat{p}(t,v,v)}\,\td\mu_{m}(v)
      \dx\\
    &\qquad\le
      \sqrt{h_m(t)}\,z\,\norm{\sqrt{\hat{p}(t,\cdot,\cdot)}}_{L^{\infty}([0,z]^2)}
      \,\mu_{m}((0,z]).
  \end{align*}
  Note, the right hand-side of the last estimate above is finite since
  $\hat{p}(t,\cdot,\cdot)$ is continuous on $[0,\infty)^2$ by
  Proposition~\ref{prop:densities-p-and-omega}.  Therefore, we have
  thereby shown that the function
  \begin{displaymath}
  g(\gamma,v):=e^{-\gamma t}C(v,\gamma)\quad\text{ belongs to
  $L^{1}((0,\infty)\times \Omega_z ;\hat{\Delta}_{m}\otimes
  (\mu_{\text{Leb}}\otimes \mu_{m}))$.}
\end{displaymath}
 Hence, by Fubini's theorem and subsequently by applying the identity~\eqref{eq:86} for
 $C(v,\gamma)$, \eqref{eq:118}, and by Corollary~\ref{prop:1}, one sees that
\begin{align*}
	&2 \int_{\Omega_z}\int_{0}^{\infty}e^{-\gamma
   t}C(v,\gamma)\,\td\hat{\Delta}_{m}(\gamma)\, \td(\mu_{\text{Leb}}\otimes \mu_m)(x,v)\\
	&\qquad=\int_{0}^{\infty} \left (2 \int_{\Omega_z}e^{-\gamma
   t}C(v,\gamma)\, \td(\mu_{\text{Leb}}\otimes \mu_m)(y,v)\right
   )\,\td\hat{\Delta}_{m}(\gamma)\\
  &\qquad=\int_{0}^{\infty} \left (2 \int_{0}^{z}\int_{(0,y]}e^{-\gamma
   t}C(v,\gamma)\, \td\mu_m)(v)\,\dy\right)\,\td\hat{\Delta}_{m}(\gamma)\\
	&\qquad=2z\int_0^\infty\frac{1}{\gamma}e^{-\gamma t}
   \td\hat{\Delta}_{m}(\gamma)-\int_0^\infty\frac{1}{\gamma}e^{-\gamma
   t}C(z,\gamma)\td\hat{\Delta}_{m}(\gamma)\\
	&\qquad=2z\,\beta_m(t)-\mathbb{P}_{z}(\tau>t).
\end{align*}
Combining this with~\eqref{eq:128}, one sees that~\eqref{eq:129} holds.
\end{proof}


\subsection{Semigroups and a bit of convex analysis}
\label{subsec:semigroups}

This section is dedicated to recall the definition of a
\emph{$C_{0}$-semigroup of contractions}, Hille-Yosida-Phillips'
characterization of the infinitesimal generator $-A$ of such a
semigroup, and to recall some important definitions and notions from
convex analysis used throughout this paper. Here, let $X$ denote a
Banach space equipped with norm $\norm{\cdot}_{X}$, $X'$ its dual
space and $\langle\cdot,\cdot\rangle_{X',X}$ the duality pairing.

\begin{definition}[{$C_{0}$-semigroup of contractions}]\label{def:c0-semigroup}
  A family $\{T_{t}\}_{t\ge 0}$ of linear operators $T_{t}\in
  \mathcal{L}(X)$ is called a \emph{$C_{0}$-semigroup of contractions}
  on $X$ provided
  \begin{enumerate}[label=(\roman*)]
    \item $T_{t+s}=T_{t}\circ T_{s}$ for every $t$, $s\ge 0$;
   \item $T_{0}=\textrm{id}_{X}$; 
    \item for every $x\in X$, the function $t\mapsto T_{t}x$ belongs
      to $C([0,\infty);X)$;
     \item $\norm{T_{t}}_{\mathcal{L}(X)}\le 1$ for all $t\ge 0$.
     \end{enumerate}
   \end{definition}

   \begin{definition}[{Infinitesimal generator}]\label{def: infinitesimal_generator}
   For a given $C_{0}$-semigroup $\{T_{t}\}_{t\ge 0}$ on $X$,
   one can define an associated \emph{infinitesimal generator} $-A$ on
   $X$ by setting
   \begin{displaymath}
     D(A):=\Big\{x\in X\,\Big\vert\, \lim_{h\to
       0+}\frac{T_{h}x-x}{h}\text{ exists in $X$}\Big\}
   \end{displaymath}
   and
   \begin{displaymath}
     -Ax:=\lim_{h\to
       0+}\frac{T_{h}x-x}{h}\quad\text{for every $x\in D(A)$.}
   \end{displaymath}
 \end{definition}

 \begin{notation}
  To emphasize that for a given operator $A$ and a
   $C_{0}$-semigroup $\{T_{t}\}_{t\ge 0}$, $-A$ is the infinitesimal
   generator of $\{T_{t}\}_{t\ge 0}$, we write
   $\{e^{-tA}\}_{t\ge 0}$ instead of $\{T_{t}\}_{t\ge 0}$.
 \end{notation}

In order to characterize operators $A$, for which $-A$ is the infinitesimal
generator of a $C_{0}$-semigroup of contractions, we require to
recall some notions from convex analysis.

\begin{definition}\label{def:sub-diff}
  For a given convex, proper functional
  $\varphi : X \to (-\infty,\infty]$ with \emph{effective domain}
  $D(\varphi):=\{u\in X\,\vert\,\phi(u)<\infty\}$, the
  \emph{sub-differential} $\partial\varphi : X\to 2^{X}$ is a possibly
  multi-valued mapping given by
  \begin{displaymath}
    \partial\varphi(u)=\Big\{x'\in X'\,\Big\vert\,\langle
    x',v-u\rangle_{X',X}\le \varphi(v)-\varphi(u)\;\forall\,v\in X\Big\}
  \end{displaymath}
  for every $u\in D(\varphi)$.
\end{definition}

\begin{definition}
  \label{def:normalized-duality-map}
  Let $\xi \in C([0,\infty))$ be a continuous, monotone, and
  surjective function satisfying $\omega(0)=0$. Then the (multi-valued) operator
  $J_{\xi} : X\to 2^{X'}$ given by
  \begin{displaymath}
    J_{\xi}(x)=\Big\{x'\in X'\,\Big\vert\, \langle
    x',x\rangle_{X',X}=\norm{x'}_{X'}\,\norm{x}_{X},\;
    \norm{x'}_{X}=\xi(\norm{x}_{X})\Big\}
  \end{displaymath}
 for every $x\in X$, is called the \emph{duality map with gauge
   function $\xi$.} In the case, the gauge function $\xi(r)=r$, $r\in
 [0,\infty)$, we simply write $J$ instead of $J_{\xi}$ and call it
 the \emph{normalized duality map}.
\end{definition}

\begin{remark}
  Since the duality map $J_{\xi}$ is the sub-differential operator
  $\partial \varphi : X\to 2^{X'}$ of the convex functional
  $\varphi(x)=\zeta(\norm{x}_{X})$ for $\zeta(r):=\int_{0}^{r}\xi(s)\,\ds$,
  $r\ge 0$, it follows that $J_{\xi}$ is \emph{monotone}, that
  is, one has that
  \begin{displaymath}
    \langle  x'_{1}-x'_{2},x_{1}-x_{2}\rangle_{X',X}\ge 0\qquad
    \text{for all pairs $(x_{1},x'_{1})$, $(x_{2},x'_{2})\in J_{\xi}$.}
  \end{displaymath}
  Moreover, for every $x\in X$, $J_{\xi}(x)$ is closed, convex and
  non-empty subset of $X'$.
\end{remark}

\begin{definition}
  \label{def:m-accretive}
  Let $J : X\to 2^{X'}$ be the \emph{normalized duality map}
  given by
  \begin{displaymath}
    J(x)=\Big\{x'\in X'\,\Big\vert\, \langle
    x',x\rangle_{X',X}
    =\norm{x}_{X}^{2}=\norm{x'}_{x'}^{2}\Big\}\quad\text{for every
      $x\in X$.}
  \end{displaymath}
  Then, a linear operator $A : D(A)\to X$ with domain $D(A)\subseteq X$ is
  called \emph{accretive} if for every $u\in D(A)$, there is an $x'\in J(u)$
  \begin{displaymath}
    \textrm{Re}\langle x',Au\rangle_{X',X}\ge 0.
  \end{displaymath}
  Further, $A$ is called \emph{$m$-accretive} if $A$ is accretive and
  satisfies the so-called \emph{range condition}; that is, for every
  $f\in X$ and $\lambda>0$, there exists a unique $u\in D(A)$ such
  that $u+\lambda Au=f$.
\end{definition}

\begin{theorem}[{\cite[Corollary 3.3.5]{MR2798103}}]\label{thm:hille-yosida}
  Let $A$ be a linear operator on $X$. Then, $-A$ is the infinitesimal
   generator of a $C_{0}$-semigroup of contractions $\{e^{-t
     A}\}_{t\ge 0}$ on $X$ if and only if $A$ is $m$-accretive on $X$.
\end{theorem}


\subsection{Subordination of semigroups}
\label{subsec:subordination}

In this section, we briefly review the definition of $\psi(A)$ for a given
Bernstein function $\psi$ and a linear $m$-accretive operator $A$ on a
Banach space $X$. Here, we follow closely~\cite[Chapter~13]{MR2978140}\medskip

We begin by recalling the following definitions.

\begin{definition}
  A finite Borel measure $\nu$ on $[0,\infty)$ is called a
  \emph{sub-probability measure} provided $\nu([0,\infty))\le 1$.
\end{definition}

\begin{definition}\label{def:vague-convergence}
  Let $(X,\tau)$ be a locally compact topological Hausdorff space. Then a sequence
  $\{\nu_{n}\}_{n\ge 1}$ of Radon measures $\nu_{n}$ on $X$
  \emph{converges vaguely} to a measure $\nu$ if
  \begin{displaymath}
    \int_{X}\varphi(x)\,\td\nu_{n}(x)\to
    \int_{X}\varphi(x)\,\td\nu(x)\qquad\text{as $n\to \infty$}
  \end{displaymath}
  for all compactly supported, continuous, real-valued functions
  $\varphi$ on $X$ (we denote this set of functions by
  $C_{c}(X)$.
\end{definition}

\begin{definition}\label{def:convolution}
  The \emph{convolution $\mu\ast\nu$ of two sub-probability measures}
  $\mu$ and $\nu$ on $[0,\infty)$ is defined by
  \begin{displaymath}
    \int_{[0,\infty)}\varphi(t)\,\td(\mu\ast\nu)(t)
    =\int_{[0,\infty)}\int_{[0,\infty)}\varphi(t+s)\,\td\mu(t)\,\td\nu(s) 
  \end{displaymath}
  for every bounded continuous function $\varphi$ on $[0,\infty)$ (which
we summarize in the set $C_{c}([0,\infty))$).
\end{definition}

We note that the convolution $\mu\ast\nu$ of two
sub-probability measures $\mu$ and $\nu$ is again a sub-probability
measure on $[0,\infty)$. With the \emph{convolution}-operation, we can
defined now the following type of semigroups. 
        
\begin{definition}\label{def:family-of-sub-probabilities}
  A family $\{\gamma_{t}\}_{t\ge 0}$ of finite Borel measures $\gamma_{t}$ on
  $[0,\infty)$ is called a \emph{vaguely continuous convolution
    semigroup of sub-probability measures} provided the family
  $\{\gamma_{t}\}_{t\ge 0}$ satisfies
\begin{enumerate}[label={(\roman*)}]
  \item $\gamma_{t}([0,\infty))\le 1$ for all
    $t\ge 0$\qquad\hfill (sub-probability condition);
  \item $\mu_{t+s}=\gamma_{t}*\mu_{s}$ for all $t$, $s\ge
        0$\qquad\hfill (semigroup property);
   \item $\displaystyle\lim_{t\rightarrow
         0}\gamma_{t}=\delta_{0}$ vaguely\qquad\hfill  (vague continuity).
  \end{enumerate}
\end{definition}

The next theorem highlights the bijective relation between Bernstein
functions $\psi$ and vaguely continuous convolution
    semigroup of sub-probability measures.

\begin{theorem}[{\cite[Theorem 5.2]{MR2978140}}]\label{thm:charact-Bernstein-conv-semigroup}
  Let $\{\gamma_{t}\}_{t\ge 0}$ be a vaguely continuous convolution semigroup of
  sub-probability measures on $[0,\infty)$. Then there exists a unique
  Bernstein function $\psi$ such that the Laplace transform of
  $\gamma_{t}$
  \begin{equation}
    \label{eq:50bis}
	\int_{0}^{\infty}e^{-\lambda
          s}\,\td\gamma_{t}(s)=e^{-t\psi(\lambda)}\qquad
        \text{for all $\lambda> 0$, $t\ge 0$.}
  \end{equation}
  Conversely, for a given Bernstein function $\psi$, there exists a
  unique vaguely continuous convolution semigroup $\{\gamma_{t}\}_{t\ge 0}$ of sub-probability
  measures on $[0,\infty)$ satisfying~\eqref{eq:50bis}. 
\end{theorem}

The following proposition provides the existence theorem of the
operator $\psi(A)$ via the infinitesimal generator of a semigroup.

\begin{proposition}[{\cite[Proposition 13.1]{MR2978140}}]
\label{prop:semigroup-property-sub-probability}
Let $\{e^{-tA}\}_{t\ge 0}$ be a $C_{0}$-semigroup of
contractions with infinitesimal generator
$A$ on a Banach space $X$, and $\{\gamma_{t}\}_{t\ge 0}$ be a vaguely continuous convolution
semigroup of sub-probability measures on $[0,\infty)$ with the
corresponding Bernstein function $\psi$. Then the family
$\{e^{-t\psi(A)}\}_{t\ge 0}$ defined by the Bochner integral
  \begin{equation}
    \label{eq:51}
    e^{-t\psi(A)} f:=\int_{[0,\infty)}e^{-sA}f\,\td\gamma_{t}(s)\qquad\text{for every
      $t\ge 0$, $f\in X$,}
  \end{equation}
  defines a $C_{0}$-semigroup $\{e^{-t\psi(A)}\}_{t\ge 0}$ of
  contractions $e^{-t\psi(A)}\in \mathcal{L}(X)$. 
\end{proposition}

\begin{notation}[{The operator $\psi(A)$}]
  For a given $m$-accretive operator $A$ on a Banach space $X$, and
  vaguely continuous convolution semigroup $\{\gamma_{t}\}_{t\ge 0}$ of
  sub-probability measures $\gamma_{t}$ associated the corresponding
  Bernstein function $\psi$, we denote by $\psi(A)$ the
  \emph{infinitesimal generator} (see Definition~\ref{def:
    infinitesimal_generator}) of the semigroup
  $\{e^{-t\psi(A)}\}_{t\ge 0}$ given by~\eqref{eq:51}.
\end{notation}

\begin{theorem}[{Phillips' subordination theorem}, {\cite[Theorem
  13.6]{MR2978140}}]\label{thm:Phillips} Let $\{e^{-tA}\}_{t\ge 0}$ be
a $C_{0}$-semigroup of contractions with
infinitesimal generator $-A$ on Banach space $X$, and $\psi$
be a Bernstein function with the L\'evy triple $(a,b,\nu)$. Consider the semigroup
$\{e^{-t\psi(A)}\}_{t\geq 0}$ with the infinitesimal generator
$-\psi(A)$. Then, the domain $D(A)$ of $A$ is an operator core of the
domain $D(\psi(A))$ of $\psi(A)$ and
  \begin{displaymath}
    \psi(A)f=af+bAf+\int_{0}^{\infty}\big(f-e^{-tA}f\big)\,\td\nu(t)
  \end{displaymath}
  for all $f\in D(A)$, where the integral is to be understood in the
  Bochner sense.
\end{theorem}

%
%
%
%
%

\section{Proofs of the Main Results}
\label{sec:proofs}


Throughout this section, let $m\in \mathfrak{m}_{\infty}$ and associated
Lebesgue-Stieltjes measure $\mu_{m}$. Further,
$\{Z_{t}\}_{t\ge 0}$ be a generalized diffusion associated with
$m$, $\tau$ the first hitting time of
$z=0$ by the process $\{Z_{t}\}_{t\ge 0}$ as defined in~\eqref{eq:46},
and $\psi_{m}$ the associated complete Bernstein
function~\eqref{eq:48} from Theorem~\ref{thm:2}. 

For the moment, we assume
that $A : D(A)\to X$ is merely a closed, linear 
operator defined on a Banach $X$ with norm $\norm{\cdot}_{X}$, and
denote by $X'$ the dual space, and $\langle\cdot,\cdot\rangle_{X',X}$
the corresponding duality brackets.\medskip 

Then, the main object of this section is the \emph{Dirichlet-to-Wentzell operator}
\begin{equation}
  \label{eq:59}
  f\mapsto \Lambda_mf:=m(0+)Au(0)-\frac{1}{2}\frac{\td u}{\dz}_{\! +}\!(0)
\end{equation}
associated with the \emph{Dirichlet problem}
\begin{equation}\label{eq:49}
\begin{cases}
Au(z)-\frac{1}{2}\frac{\td}{\dm}\frac{\td}{\dz}u(z)&=0, \quad \quad \text{for $z\in (0,\infty)$,}\\
\phantom{-Au(y)+A_{Y}}u(0)&=f,
\end{cases}
\end{equation}
for given $f\in D(A)$. In the next subsection, we briefly
discuss the notion of a weak solution $u$ to the Dirichlet
problem~\eqref{eq:49} introduced in Definitions~\ref{def:weak-solution
  Dirichlet-problem}. 


\subsection{Weak solutions}

We begin by discussing the  the notion of weak solutions $u$ to the Dirichlet
problem~\eqref{eq:49} as it was introduced in Definition~\ref{def:weak-solution
	Dirichlet-problem}.

\begin{remark}\label{REMARK:SOLN}
  For a given string $m\in \mathfrak{m}_{\infty}$, $\mu_{m}$ the associated measure to $m$, and $E_{m}$
  the support of $\mu_{m}$. Then the
  following comments are worth noting.
  \begin{enumerate}[label=(\alph*)]
    \item \label{REMARK:SOLN-1} We show in the
      Proposition~\ref{propo:chara-weak-m-derivative} (in the appendix)
      that  a function
  $f\in L^{1}_{loc}((0,\infty);X)$ has a weak $m$-derivative
  $g\in L^{1}_{loc, \mu_{m}}((0,\infty);X)$ if and only if $f$ can be
  represented by
  \begin{displaymath}
    f(z_{2})=f(z_{1})+\int_{z_{2}}^{z_{1}}g(r)\,\td\mu_{m}(r)\qquad\text{for a.e. $z_1,z_2\in (0,\infty)$.}
  \end{displaymath}

    \item \label{REMARK:SOLN-2} Due to Remark~\ref{REMARK:SOLN-1}, for
      given $f\in X$, a \emph{weak solution} $u$
      of Dirichlet problem~\eqref{eq:49} is characterized by
      satisfying $u\in C([0,\infty);X)\cap W^{1,1}_{loc}([0,\infty);X)$, the weak derivative
      $\frac{\td u}{\dz}$ is weakly $m$-differentiable, and there is a
      $g\in L^{1}_{loc,\mu_{m}}([0,\infty);X)$ satisfying
      \begin{equation}
        \label{eq:98}
        u(z)=f+\frac{\td u}{\dz}_{\! -}\!(0)\,z
        +\int_{0}^{z}\int_{[0,y]}g(r)\,\td\mu_{m}(r)\,\dy
      \end{equation}
      for every $z\in [0,\infty)$, where $\frac{\td u}{\dz}_{\! -}\!(0)$
      denotes the left hand-side derivative of $u$ at $z=0$, and
      \begin{displaymath}
        2Au(z)=g(z)\qquad\text{ for $\mu_{m}$-a.e. $z\in (0,\infty)$.}
    \end{displaymath}
    This formulation consistent with the real-valued function case described
      in~\cite{MR2978140} provided by Revuz and Yor.

    \item \label{REMARK:SOLN-3} In comparison to
      Remark~\ref{rem:m-derivative}, an alternative characterization
      of a weak solution $u$ is of the \emph{extension equation}
         \begin{equation}
           \label{eq: REMARK:SOLN-3}
           Au(z)-\frac{1}{2}\frac{\td}{\dm}\frac{\td}{\dz}u(z)=0 \qquad\text{in $z\in (0,\infty)$,}
         \end{equation}
         is given by $u$ belonging to $W^{1,1}_{loc}((0,\infty);X)$,
         for $\mu_{m}$-a.e. $z\in (0,\infty)$, one has that $u(z)\in D(A)$,
         $Au\in L^{1}_{loc, \mu_{m}}((0,\infty);X)$, and
         \begin{equation}
           \label{eq-2: REMARK:SOLN-3}
           \frac{1}{2}\frac{\td^{2} u}{\td z^2}=Au\,\mu_{m}\qquad\text{in $\mathcal{D}'((0,\infty);X)$.}
         \end{equation}

    \item \label{Remark:2} A function $u$ given by~\eqref{eq:98} with weak
      $m$-derivative
      \begin{displaymath}
        \frac{\td }{\td m}\frac{\td u}{\dz}\in
      L^{1}_{loc,\mu_{m}}([0,\infty);X)
    \end{displaymath}
    is linear on $[0,\infty)\setminus E_{m}$, the left and
    right hand-side derivative $\frac{\td u}{\dz}_{\! -}\!(z)$, $\frac{\td u}{\dz}_{\! +}\!(z)$ exist at every
    $z\in [0,\infty)$ with
      \begin{align*}
        \frac{\td u}{\dz}_{\! +}\!(z)&=\frac{\td u}{\dz}_{\! -}\!(0)+\int_{[0,z]}g(r)\,\td\mu_{m}(r),\\
        \frac{\td u}{\dz}_{\! -}\!(z)&=\frac{\td u}{\dz}_{\! -}\!(0)+\int_{[0,z)}g(r)\,\td\mu_{m}(r),
      \end{align*}
      and $u'$ is a.e. continuous on $[0,\infty)$.
  \end{enumerate}
\end{remark}

With these comments in mind, we can now start by establishing
uniqueness of bounded solutions of Dirichlet problem~\eqref{eq:49}.

\subsection{Uniqueness of weak solutions of the Dirichlet problem}
\label{subsec:uniqueness}
In this subsection, we outline our uniqueness result of bounded
\emph{solutions} of Dirichlet problem~\eqref{eq:49}. Our method to
prove uniqueness relies essentially on tools and arguments borrowed
from the theory of nonlinear evolution theory (see, for
instance,~\cite{Benilan72,MR4026441}).\medskip

Throughout this section, we
assume that the operator $A$ in the extension problem~\eqref{eq:49} is
\emph{accretive} on $X$ as defined in
Definition~\ref{def:m-accretive}.\medskip

Our first lemma partially generalizes a result by B\'enilan~\cite{Benilan72}.

\begin{lemma}
  \label{lem:differentiability}
  For a given string $m\in \mathfrak{m}_{\infty}$, let $\mu_{m}$ be the
  associated measure to $m$. Let $\varphi : X\to \R$ be a continuous,
  convex functional on a Banach space $X$ satisfying $\varphi(0)=0$, and
  for some $g\in L^{1}_{loc,\mu_{m}}([0,\infty);X)$, $x$, $\beta\in X$,
  let $u\in C([0,\infty);X)$ be given by
  \begin{equation}
    \label{eq:100}
    u(z)=x+\beta\,z+\int_{0}^{z}\int_{[0,y]}g(r)\,\td\mu_{m}(r)\,\dy
  \end{equation}
 for every $z\in [0,\infty)$. Then the following statements hold.
  \begin{enumerate}
  \item\label{lem:differentiability-claim1} The mapping
    $z\mapsto \varphi(u(z))$ is differentiable from the right at every
    $z\in [0,\infty)$ with
    \begin{equation}
      \label{eq:101}
      \frac{\td}{\dz}_{\!
        +}\varphi(u(z))=\max_{w\in \partial\varphi(u(z))}\langle
      w,\frac{\td u}{\dz}_{\! +}\! (z)\rangle_{X',X},
    \end{equation}
    where $\partial\varphi$ denotes the sub-differential of $\varphi$
    (see Definition~\ref{def:sub-diff}).
  \item\label{lem:differentiability-claim2}  The mapping $z\mapsto
    \frac{\td}{\dz}_{\! +}\varphi(u(z))$ is $m$-differen\-tiable
    $\mu_{m}$-a.e. on $[0,\infty)$, and locally of bounded variation.
  \item\label{lem:differentiability-claim3} For every $z\in [0,\infty)$
    such that $\frac{\td }{\td m}\frac{\td u}{\dz}_{\! +}\! (z)$ and
    $\frac{\td }{\td m}\frac{\td}{\dz}_{\! +}\varphi(u(z))$ exist, one
    has that
    \begin{equation}
      \label{eq:104}
        \frac{\td }{\td m}\frac{\td}{\dz}_{\! +}\varphi(u(z))\ge 
        \langle w, \frac{\td }{\td m}\frac{\td u}{\dz}_{\! +}\!(z)\rangle_{X',X}
    \end{equation}
    for every $w\in \partial\varphi(u(z))$.  
  \end{enumerate}
\end{lemma}

\begin{proof}
  Let $u$ be given by~\eqref{eq:100}. Then $u$ is differentiable from
  the right at every $z\in [0,\infty)$. Thus,
  claim~\eqref{lem:differentiability-claim1} follows by the same
  arguments as in~\cite[Lemme~1]{Benilan72} with $u'(z)$ replaced by
  $\frac{\td u}{\dz}_{\! +}\!(z)$.

  To see that claim~\eqref{lem:differentiability-claim2} holds, let
  $0\le z_{1}\le z_{2}\le T<\infty$ and for $i=1, 2$, let
  $w_{i}\in \partial\varphi(u(z_{i}))$ such that~\eqref{eq:101}
  holds. Then,
  \begin{align*}
    &\frac{\td}{\dz}_{\! +}\varphi(u(z_{2}))-\frac{\td}{\dz}_{\!
      +}\varphi(u(z_{1}))\\
    &\qquad= \langle w_{2},\frac{\td u}{\dz}_{\! +}\!(z_{2})-\frac{\td u}{\dz}_{\! +}\!(z_{1})\rangle_{X',X}
      +\langle w_{2}-w_{1},\frac{\td u}{\dz}_{\! +}\!(z_{1})\rangle_{X',X}.
  \end{align*}
  Since
  \begin{displaymath}
    \frac{\td u}{\dz}_{\! +}\!(z)=\beta+\int_{[0,z]}g(r)\,\td\mu_{m}(r),\qquad z\in [0,\infty),
  \end{displaymath}
  we have that
  \begin{align*}
    \langle w_{2},\frac{\td u}{\dz}_{\! +}\!(z_{2})-\frac{\td u}{\dz}_{\! +}\!(z_{1})\rangle_{X',X}
    &=\langle w_{2},
      \int_{(z_{1},z_{2}]}g(r)\,\td\mu_{m}(r)\rangle_{X',X}\\
    &\ge - \norm{w_{2}}_{X'}\, \int_{(z_{1},z_{2}]}\norm{g(r)}_{X}\,\td\mu_{m}(r).
  \end{align*}
  On the other hand, since
  \begin{align*}
   \frac{u(z_{2})-u(z_{1})}{z_{2}-z_{1}}
    &= \frac{\td u}{\dz}_{\! +}\!(z_{1})-\int_{[0,z_{1}]}g(r)
      \,\td\mu_{m}(r)\\
    &\hspace{3cm}+\tfrac{1}{z_{2}-z_{1}}\int_{z_{1}}^{z_{2}}\int_{[0,z]}g(r)\,\td\mu_{m}(r)\,\dz,
  \end{align*}
 one has that 
 \begin{equation}
   \label{eq:102}
   \frac{u(z_{2})-u(z_{1})}{z_{2}-z_{1}}=\frac{\td u}{\dz}_{\! +}\!(z_{1})
   -\tfrac{1}{z_{2}-z_{1}}\int_{z_{1}}^{z_{2}}\int_{(z_{1},z]}g(r)\,\td\mu_{m}(r)\,\dz.
 \end{equation}
 Therefore and since $\partial\varphi$ is monotone, 
  \begin{align*}
    \langle w_{2}-w_{1},\frac{\td u}{\dz}_{\! +}\!(z_{1})\rangle_{X',X}
    &\ge \langle w_{2}-w_{1},
      \tfrac{1}{z_{2}-z_{1}}\int_{z_{1}}^{z_{2}}\int_{(z_{1},z]}g(r)\,\td\mu_{m}(r)\,\dz\rangle_{X',X}\\
    &\ge-\norm{w_{2}-w_{1}}_{X'}\, \int_{(z_{1},z_{2}]}\norm{g(r)}_{X}\,\td\mu_{m}(r).
  \end{align*}
 By hypothesis, the functional $\varphi$ is continuous on $X$. Hence,
 on every bounded subset $U$ of $X$, there is a constant $M>0$ such
 that
 \begin{displaymath}
   \norm{x'}_{X'}\le M\qquad\text{for every $x'\in \partial\varphi(u)$ and $u\in U$.}
\end{displaymath}
Further, by the continuity of $u : [0,T]\to X$, there is a bounded
open subset $U$ of $X$ containing $u([0,T])$. Thus, there is an $M>0$
such that
 \begin{displaymath}
   \langle w_{2},\frac{\td u}{\dz}_{\! +}\!(z_{2})
   -\frac{\td u}{\dz}_{\! +}\!(z_{1})\rangle_{X',X}\ge
   - M\, \int_{(z_{1},z_{2}]}\norm{g(r)}_{X}\,\td\mu_{m}(r).
 \end{displaymath}
 and
 \begin{displaymath}
   \langle w_{2}-w_{1},\frac{\td u}{\dz}_{\!
     +}\!(z_{1})\rangle_{X',X}\ge
   -2\,M\, \int_{(z_{1},z_{2}]}\norm{g(r)}_{X}\,\td\mu_{m}(r)
 \end{displaymath}
 and hence,
 \begin{displaymath}
   \frac{\td}{\dz}_{\! +}\varphi(u(z_{2}))-\frac{\td}{\dz}_{\!
      +}\varphi(u(z_{1}))\ge -3\,M\, \int_{(z_{1},z_{2}]}\norm{g(r)}_{X}\,\td\mu_{m}(r),
 \end{displaymath}
or, equivalently, for
\begin{displaymath}
  h(z):=\frac{\td}{\dz}_{\! +}\varphi(u(z))\quad\text{and}\quad
  H(z):=\int_{[0,z]}\norm{g(r)}_{X}\,\td\mu_{m}(r),
\end{displaymath}
one has that $h(z_{2})+H(z_{2})\ge h(z_{1})+H(z_{1})$. Since
$0\le z_{1}\le z_{2}\le T< \infty$ were arbitrary, we have thereby
shown that $z\mapsto h(z)+H(z)$ is monotonically increasing along
$[0,T]$, and hence $m$-differentiable $\mu_{m}$-a.e. on
$[0,\infty)$. Since $H$ is also monotone, it is also $m$-differentiable
$\mu_{m}$-a.e. on $[0,\infty)$, which implies that
$h(z)=\frac{\td}{\dz}_{\!  +}\varphi(u(z))$ is $m$-differentiable
$\mu_{m}$-a.e. on $[0,\infty)$. This completes the proof of
statement~\eqref{lem:differentiability-claim2}.

Next, let $u$ be given by~\eqref{eq:100}. Then, one has that
  \begin{equation}
    \label{eq:95}
      \frac{\td }{\td m}\frac{\td u}{\dz}_{\! +}\!(z)=g(z)
  \end{equation}
for $\mu_{m}$-a.e. $z\in (0,\infty)$. From this, we can infer that
\begin{equation}
  \label{eq:103}
  \lim_{h\to 0+}\frac{u(z+h)-2u(z)+u(z-h)}{h\,(m(z)-m(z-h))}=\frac{\td
  }{\td m}\frac{\td u}{\dz}_{\! +}\!(z)
  \qquad\text{in $X$}
\end{equation}
for $\mu_{m}$-a.e. $z\in (0,\infty)$. To see that~\eqref{eq:103}
holds, we first write
\begin{displaymath}
  \frac{u(z+h)-2u(z)+u(z-h)}{h\,(m(z)-m(z-h))}=
        \frac{u(z+h)-u(z)}{h\,(m(z)-m(z-h))}-\frac{u(z)-u(z-h)}{h\,(m(z)-m(z-h))}.                                        
\end{displaymath}
Then for $h>0$, \eqref{eq:102} yields that
\begin{displaymath}
  u(z+h)-u(z)=\frac{\td u}{\dz}_{\! +}\!(z)h
  -\int_{z}^{z+h}\int_{(z,y]}g(r)\,\td\mu_{m}(r)\,\dy
\end{displaymath}
and
\begin{displaymath}
  u(z)-u(z-h)=\frac{\td u}{\dz}_{\! +}\!(z-h)h
  -\int_{z-h}^{z}\int_{(z-h,y]}g(r)\,\td\mu_{m}(r)\,\dy.
\end{displaymath}
Thus,
\begin{align*}
  &\frac{u(z+h)-u(z)-(u(z)-u(z-h))}{h\, (m(z)-m(z-h))}\\
  &\qquad =
    \frac{\frac{\td u}{\dz}_{\! +}\!(z)-\frac{\td u}{\dz}_{\! +}\!(z-h)}{m(z)-m(z-h)}-\tfrac{1}{h\,(m(z)-m(z-h))}
    \int_{z}^{z+h}\int_{(z,y]}g(r)\,\td\mu_{m}(r)\,\dy\\
  &\hspace{4cm}+\tfrac{1}{h\,(m(z)-m(z-h))}\int_{z-h}^{z}\int_{(z-h,y]}g(r)\,\td\mu_{m}(r)\,\dy.
\end{align*}
If \eqref{eq:95} holds at $z\in (0,\infty)$ and if $z$ a Lebesgue point
of $\mu_{m}$, then the last two terms on the right hand-side of the
latter equation tend to zero as $h\to 0+$, showing that~\eqref{eq:103}
holds.

Finally, for $w\in \partial\phi(u(z))$, the monotonicity of
$\partial\phi$ implies that
\begin{equation}
  \label{eq:96}
\begin{split}
  &\frac{\phi(u(z+h))-2\phi(u(z))+\phi(u(z-h))}{h\,(m(z)-m(z-h))}\\
  &\qquad \qquad \ge \frac{\langle w, u(z+h)-2u(z)+u(z-h)\rangle_{X',X}}{h\,(m(z)-m(z-h))}\\
  &\qquad \qquad =\langle w, \frac{u(z+h)-2u(z)+u(z-h)}{h\,(m(z)-m(z-h))}\rangle_{X',X} .
\end{split}
\end{equation}
Thus, for every $z\in [0,\infty)$ such that
$\frac{\td }{\td m}\frac{\td u}{\dz}_{\! +}\! (z)$ and
$\frac{\td }{\td m}\frac{\td}{\dz}_{\! +}\varphi(u(z))$ exist, taking
the limit as $h\to 0+$ in~\eqref{eq:96} yields that~\eqref{eq:104}
holds.
\end{proof}

With the help of Lemma~\ref{lem:differentiability}, we obtain
uniqueness of bounded solutions of the Dirichlet
problem~\eqref{eq:49} in a Banach spaces $X$.

\begin{theorem}[{Uniqueness of bounded solutions}]\label{thm:uniqueness}
  Let $m\in \mathfrak{m}_{\infty}$
  with associated measure $\mu_{m}$, and $A$ an accretive
  operator on the Banach space $X$. Then, for every given $f\in X$,
  there is at most one solution $u\in L^{\infty}([0,\infty);X)$ of Dirichlet
  problem~\eqref{eq:49}.
\end{theorem}

\begin{proof}
  Let $u\in L^{\infty}([0,\infty);X)$ be a solution of Dirichlet
  problem~\eqref{eq:49} with initial data $f=0$ and
  $\varphi(x):=\frac{1}{2}\norm{x}_{X}^{2}$ for every $x\in X$. Since the operator
  $A$ is accretive on $X$, one has that
  \begin{displaymath}
    \langle w, \frac{\td }{\td m}\frac{\td u}{\dz}_{\!
      +}\!(z)\rangle_{X',X}\ge 0\qquad\text{for every $w\in
      \partial\varphi(u(z))$.}
  \end{displaymath}
  Hence by~\eqref{eq:104} of Lemma~\ref{lem:differentiability}, the
  function $z\mapsto \norm{u(z)}_{X}^{2}$ is
  convex on $[0,\infty)$. Since $u$ is bounded, also the function $z\mapsto \norm{u(z)}_{X}^{2}$
  is bounded and hence monotonically decreasing. This implies that
  \begin{displaymath}
    \norm{u(z)}_{X}^{2}\le \norm{u(0)}_{X}^{2}=0\qquad\text{ for all $z\in
  [0,\infty)$,}
\end{displaymath}
completing the proof of this theorem.
\end{proof}

\subsection{Existence of weak solutions of the Dirichlet problem}
\label{subsec:existence}

In this section, we establish existence of bounded solutions $u$ of
Dirichlet problem~\eqref{eq:49} under the hypothesis that $A$ is an
$m$-accretive operator on a Banach space $X$, and for given initial
value $f\in D(A)$. Our approach to this combines arguments from
stochastic analysis with linear semigroup theory. In particular, it is
based on the spectral representation of the \emph{first hitting time}
density $\omega_{\tau}$ of a given generalized diffusion
$\{Z_{t}\}_{t\ge 0}$ on $[0,\infty)$ (see Theorem \ref{DTN:Levy Measure}). Thus, we
employ the same notation here as introduced in Subsection
\ref{subsec:hitting-times}.\medskip

Throughout this subsection, let $m\in \mathfrak{m}_{\infty}$ with associated measure $\mu_{m}$. Further,
let $A$ an $m$-accretive operator on $X$ and $\{e^{-tA}\}_{t\ge 0}$
the semigroup generated by $-A$. Further, let $\{Z_{t}\}_{t\ge 0}$ be
a generalized diffusion on $[0,\infty)$ associated with $m$,
$B_{m}=\tfrac{1}{2}\tfrac{\td}{\td m}\frac{\td}{\dz}$ the
infinitesimal generator of $\{Z_{t}\}_{t\ge 0}$, and $\omega_{\tau}$
the density~\eqref{eq:45} of the first hitting time
 $\tau$ of $z=0$ by $\{Z_{t}\}_{t\ge 0}$. Then, for given $f\in X$, we set
 \begin{equation}\label{DTN:Exact Solution}
   u(z):=\int_{0}^{\infty}(e^{-tA}f)\,
   \omega_{\tau}(t,z)\,\dt, \qquad \,\, z\in [0,\infty).
 \end{equation}
Then, our aim is to show that the function $u$
 given by~\eqref{DTN:Exact Solution} is a solution of the Dirichlet
 problem~\eqref{eq:49}.\medskip
 
Since the semigroup $\{e^{-t A}\}_{t\ge 0}$ is contractive, one easily
sees that the integral in~\eqref{DTN:Exact Solution} is finite. Moreover, one has that 
 \begin{equation}\label{EQN:BND_SOLN}
   \sup_{z\in [0,\infty)}\norm{u(z)}_{X}
   \le \sup_{z\in [0,\infty)}\int_{0}^{\infty}\norm{e^{-t A}f}_{X}\omega_{\tau}(t,z)\,dt \le \norm{f}_{X},
 \end{equation}
 showing that $u\in L^{\infty}(0,\infty;X)$.

 Next, we intend to show that for given $f\in D(A)$,
 the function $u$ given by~\eqref{DTN:Exact Solution} is a weak solution
 of the extension equation~\eqref{eq: REMARK:SOLN-3} satisfying
 $u(0)=f$. More specifically, we show that $u$ can be rewritten as follows. 

 \begin{theorem}
   \label{prop:existence}
   In addition to the hypotheses of this subsection, let $f\in D(A)$, and
   $u$ be given by~\eqref{DTN:Exact Solution}. Then $u$ belongs to $W^{1,1}_{loc}([0,\infty);X)$
    and satisfies
    \begin{enumerate}
    \item  \label{prop:existence-claim1} $u(0)=f$ in $X$,
    \item  \label{prop:existence-claim2} $u(z)\in D(A)$ for every $z\in [0,\infty)$,
    \item  \label{prop:existence-claim3} $Au\in L^{\infty}([0,\infty);X)$, and
    \item  \label{prop:existence-claim4} $u$ can be rewritten as
      \begin{equation}\label{EQN:SOLN}
       \begin{split}
          u(z)&=f-2z\int_{0}^{\infty}\left[f-e^{-tA}f\right]\,h_m(t)\,\dt\\
          &\hspace{4cm}+2\int_{0}^{z}\int_{(0,x]}Au(v)\,\td
          \mu_m(v)\,\dx
       \end{split}
     \end{equation}
     for every $z\in [0,\infty)$.
    \end{enumerate}
\end{theorem}

We note that by the Radon property of the measure $\mu_{m}$ it follows
that $L^{\infty}(0,\infty;X)$ is contained in
$L^{1}_{loc,\mu_{m}}([0,\infty);X)$. Thus, by Remark~\ref{REMARK:SOLN},
the representation~\eqref{EQN:SOLN} of $u$ implies that $u$ is a \emph{weak solution} of Dirichlet
problem~\eqref{eq:49}. The statement of Theorem~\ref{prop:existence}
establishes our main result Theorem~\ref{thm:1DP}.

The proof of
Theorem~\ref{prop:existence} proceeds in several steps. We begin by
showing that $u(z)\in D(A)$ for every $z\in (0,\infty)$ and that
statement~\ref{prop:existence-claim3} holds.
 
 \begin{lemma}\label{LEMMA:Au_z}
   In addition to the hypotheses of this subsection, let $f\in D(A)$, and
   $u$ be given by~\eqref{DTN:Exact Solution}. Then for
   every $z\in (0,\infty)$, $u(z)\in D(A)$ and
   \begin{equation}
     \label{eq:130}
     Au(z)=\int_{0}^{\infty}(Ae^{-tA}f)\, \omega_{\tau}(t,z)\,\dt.
   \end{equation}
Moreover, one has that 
\begin{equation}\label{EQN:BND_SOLN}
\sup_{z\in [0,\infty)}\norm{Au(z)}_{X}\le \norm{Af}_{X}.
\end{equation}
 \end{lemma}

\begin{proof}
  Let $f\in D(A)$, $u$ be given by~\eqref{DTN:Exact Solution}, and
  $z\in (0,\infty)$. Then by $f\in D(A)$, we have that
  $Ae^{-tA}f=e^{-tA}(Af)$, and since the family $\{e^{-tA}\}_{t\geq 0}$
  consists of contractive operators $e^{-tA} \in \mathcal{L}(X)$, one
  sees that
  \begin{equation}\label{eq:131}
   \begin{split}
     \int_{0}^{\infty}\norm{Ae^{-tA}f}_{X}\, \omega_{\tau}(t,z)\,\dt
     &=\int_{0}^{\infty}\|e^{-tA}(Af)\|_{X}\, \omega_{\tau}(t,z)\,\td t\\
     &\le \,\norm{Af}_{X}\,\int_{0}^{\infty}\omega_{\tau}(t,z)\,\dt\\
     & \le\, \norm{Af}_{X},
   \end{split}
  \end{equation}
  showing that the function $g(t):=e^{-tA}f\, \omega_{\tau}(t,z)$ is
  Bochner integrable. Since $A$ is a closed linear operator on a
  Banach space $X$, $g(t)\in D(A)$ for every $t\ge 0$, it follows from
  \cite[see Proposition 1.1.7]{MR2798103} that
  $u(z)=\int_{0}^{\infty}g(t)\,\dt$ belongs to $D(A)$
  and~\eqref{eq:130} holds. Now, thanks to~\eqref{eq:130}, the
  estimates in~\eqref{eq:131} show that~\eqref{EQN:BND_SOLN} holds.
\end{proof}

Our next step is to calculate for the function $u$ given by \eqref{DTN:Exact
  Solution} the integral on the left hand-side of \eqref{eq:132} below.
To do this, we employ the integral identity from Lemma~\ref{LEMMA:dou_int}.

\begin{lemma}\label{lem:ibp}
  In addition to the hypotheses of this subsection, let $f\in D(A)$,
  $u$ be given by~\eqref{DTN:Exact Solution}, and $\beta_m$ the
  antiderivative of the density $h_{m}$ of the L\'evy measure
  $\nu_{m}$ is defined as in \eqref{eq:118}. Then one has that
  \begin{equation}
    \label{eq:132}
    \begin{split}
    2\int_{0}^{z}\int_{(0,x]}Au(v)\,\td \mu_m(v)\,\dx&=2z\int_{0}^{\infty}(Ae^{-tA}f)\beta_m(t)\,\dt\\
    &\hspace{1cm}-\int_{0}^{\infty}(Ae^{-tA}f)\,\mathbb{P}_{z}(\tau>t)\,\dt,
  \end{split}
\end{equation}

  for every $z\in (0,\infty)$.
\end{lemma}

\begin{proof}
  Let $f\in D(A)$ and $z\in (0,\infty)$. Then, by~\eqref{eq:130},
  \begin{displaymath}
    \int_{0}^{z}\int_{(0,x]}Au(v)\,\td \mu_m(v)\,\dx
    =\int_{0}^{z}\int_{(0,x]}\, \int_{0}^{\infty}(Ae^{-tA}f)\, \omega_{\tau}(t,v)\,\dt\,\td\mu_m(v)\,\dx.
  \end{displaymath}
  Since 
  \begin{displaymath}
     \int_{0}^{z}\int_{(0,x]}\int_{0}^{\infty}\norm{Ae^{-tA}f}_{X}\,
     \omega_{\tau}(t,v)\,\dt\,\td\mu_m(v)\,\dx
     \le z\,\mu_m((0,z])\,\norm{Af}_{X},
   \end{displaymath}
   it follows from Fubini's theorem for Bochner Integrals (see, for
   instance, \cite[Theorem 1.1.9]{MR2798103}) and by Lemma~\ref{LEMMA:dou_int} that
   \begin{align*}
     &2\,\int_{0}^{z}\int_{(0,x]}Au(v)\,\td \mu_m(v)\,\dx\\
     &\hspace{2cm}=2\, \int_{0}^{z}\int_{(0,x]}\,
       \int_{0}^{\infty}(Ae^{-tA}f)\,\omega_{\tau}(t,v)\,\dt\,\td\mu_m(v)\,\dx\\
     &\hspace{2cm}=\int_{0}^{\infty}(Ae^{-tA}f)\, \left( 2\int_{0}^{z}\int_{(0,x]}\omega_{\tau}(t,v)\,\td
       \mu_m(v)\,\dx\right)\,\dt\\
     &\hspace{2cm}=\int_{0}^{\infty}(Ae^{-tA}f)\, \Big( 2\,z\, \beta_m(t) 
        -\mathbb{P}_{z}(\tau>t)\Big)\,\dt.
   \end{align*}
   From this, one sees that~\eqref{eq:132} holds.
\end{proof}

We are now ready to prove the remaining statements of Theorem~\eqref{prop:existence}.

\begin{proof}[Proof of Theorem~\ref{prop:existence}]
  Given $f\in D(A)$ and for $z\in (0,\infty)$, let  $u$ be given
  by~\eqref{DTN:Exact Solution}. We focus on proving the characterization~\eqref{EQN:SOLN}.
  By Lemma \ref{lem:ibp}, and since
  \begin{displaymath}
    \tfrac{\td}{\td t}e^{-tA}f+Ae^{-tA}f=0,
  \end{displaymath}
  it follows that
  \begin{equation}\label{eq:87}
    \begin{split}
      &2\int_{0}^{z}\int_{(0,x]}Au(v)\,\td \mu_m(v)\,\dx\\
      &\qquad =2z\int_{0}^{\infty}(Ae^{-tA}f)\beta_m(t)\,\td t-\int_{0}^{\infty}(Ae^{-tA}f)\mathbb{P}_{z}(\tau>t)\,\dt\\
      &\qquad
      =2z\int_{0}^{\infty}(Ae^{-tA}f)\beta_m(t)\,\dt+\int_{0}^{\infty}\mathbb{P}_{z}(\tau>t)\frac{\td
      }{\td t}e^{-tA}f\,\dt.
    \end{split}
  \end{equation}
  To complete the proof of the representation~\eqref{EQN:SOLN}, it
  remains to compute the integral terms on the left hand-side of~\eqref{eq:87}.
  Since $\omega_{\tau}(t,z)$ is the probability density of the first
  hitting time $\tau$, and since $t\mapsto \omega_{\tau}(t,z)$ is
  continuous on $(0,\infty)$, we have that
  \begin{displaymath}
	\frac{\td}{\td
          t}\mathbb{P}_{z}(\tau>t)=-\omega_{\tau}(t,z)\qquad\text{for
          every $t>0$ and $z\in (0,\infty)$.}
  \end{displaymath}
  By using this together with an integration by parts, one sees that
  \begin{equation}\label{eq:88}
      \int_{0}^{\infty}\mathbb{P}_{z}(\tau>t)\frac{\td}{\dt}e^{-tA}f\,\dt
      =-f+\int_{0}^{\infty}e^{-tA}f\,\omega_\tau(t,z)\,\dt
       =-f+u(z)
   \end{equation}
  completing the computation of one of the two integral terms on the left hand-side of~\eqref{eq:87}.  
  Before we can compute also the second integral term on the left
  hand-side of~\eqref{eq:87}, we require proving that the two limits
  \begin{equation}
    \label{eq:133}
    \lim_{t\rightarrow 0^+}\beta_m(t)\,\Big(f-e^{-tA}f\Big)
    =\lim_{t\rightarrow \infty}\beta_m(t)\,\Big(f-e^{-tA}f\Big)=0
  \end{equation}
  exist in $X$. Since $f\in D(A)$, and since $\{e^{-tA}\}_{t\ge 0}$ is
  contractive on $X$, one has that
  \begin{displaymath}
    \norm{f-e^{-tA}f}_{X}=\lnorm{\int_{0}^{t}e^{-sA}(-Af)\,\ds}_{X}\le t\,\norm{Af}_{X}.
  \end{displaymath}
  Thus, and by~\eqref{COR:LIMITS_INT} from Proposition~\ref{COR:LIMITS_INT}, one sees that
  \begin{displaymath}
    \lim_{t\rightarrow 0^+}\abs{\beta_m(t)}\,\norm{f-e^{-tA}f}_{X}
    \le \lim_{t\rightarrow 0^+}\norm{Af}_{X}\,\abs{t\beta_m(t)}=0.
  \end{displaymath}
  To see that the second limit in~\eqref{eq:133} holds as well, we use
  again that the semigroup $\{e^{-tA}\}_{t\ge 0}$ is
  contractive on $X$ and apply~\eqref{COR:LIMITS_INT}. Then, one
  easily sees that
  \begin{displaymath}
    \lim_{t\rightarrow \infty}\abs{\beta_m(t)}\,\norm{f-e^{-tA}f}_{X}
    \le \lim_{t\rightarrow \infty}\abs{\beta_m(t)}\, 2\, \norm{f}_{X}=0. 
  \end{displaymath}
  Now, by \eqref{Cor:int_by_prts} from
  Proposition~\ref{COR:LIMITS_INT}, an integration by parts yields that 
  \begin{align*}%
      &2\,z\int_{0}^{\infty}\Big(Ae^{-tA}f \Big)\,\beta_m(t)\,\dt\\
      &\hspace{2cm}=
      2\,z\int_{0}^{\infty}\frac{\td}{\td t}\Big(f-e^{-tA}f \Big)\,\beta_m(t)\,\dt\\
      &\hspace{2cm}=2\,z\int_{0}^{\infty}\Big(f-e^{-tA}f
        \Big)\,h_m(t)\,\dt
        -\lim_{t\rightarrow \infty}2\,z\, \Big(f-e^{-tA}f \Big)\beta_m(t)\\
      &\hspace{5cm}+\lim_{t\rightarrow 0^+}2\,z\, \Big (f-e^{-tA}f
        \Big)\beta_m(t)
  \end{align*}
  and so, \eqref{eq:133} implies that
  \begin{equation}
    \label{eq:89}
      2\,z\int_{0}^{\infty}\Big(Ae^{-tA}f \Big)\,\beta_m(t)\,\dt =
      2\,z\int_{0}^{\infty}\Big(f-e^{-tA}f \Big)\,h_m(t)\,\dt.
  \end{equation}
  Applying~\eqref{eq:88} and~\eqref{eq:89} to~\eqref{eq:87}, one finds
  the desired representation~\eqref{EQN:SOLN} of $u$. Further, by the
  integral representation~\eqref{EQN:SOLN}, it follows that $u$ belongs
  to $W^{1,1}_{loc}([0,\infty);X)$ and $u(0)=f$ in $X$. This
  completes the proof of Theorem~\ref{prop:existence}.
\end{proof}


\subsection{Proof of Theorem~\ref{thm:1}: a first hitting time
   approach}\label{General_CASE}
 In this subsection, we outline the proof of Theorem~\ref{thm:1} based
 on a first hitting time approach.

 \begin{proof}[Proof of Theorem~\ref{thm:1}]
   Let $f\in D(A)$, and $u$ be the unique weak solution of
   Dirichlet problem~\eqref{eq:49}. Then, by the continuity of $u$ and
   by~\eqref{EQN:BND_SOLN}, the function
   \begin{displaymath}
     g(x):=\int_{(0,x]}Au(v)\,\td \mu_m(v)\qquad \text{for every $x\in [0,\infty)$,}
   \end{displaymath}
   is right-continuous at $x=0$ and $g(0+)=0$. Therefore, we
   can conclude from \eqref{EQN:SOLN} that
   \begin{align*}
     \frac{\td u}{\dz}_{\! +}\!(0)&=\lim_{h\to
                                    0+}\frac{u(h)-u(0)}{h}\\
                                  &=\lim_{h\to 0+}\left[ (-2)
                                    \int_{0}^{\infty}(f-e^{-tA}f)\,h_m(t)\,\dt\right.\\
                                  &\hspace{4cm}+\left.
                                    \frac{2}{h}\int_{0}^{h}\int_{(0,x]}Au(v)\,\td
                                    \mu_m(v)\,\dx\right]\\
     &=(-2) \int_{0}^{\infty}(f-e^{-tA}f)\,h_m(t)\,\dt
   \end{align*}
  and so,
   \begin{align*}
     \Lambda_mf&=m(0+)Au(0)-\frac{1}{2}\frac{\td u}{\dz}_{\! +}\!(0)\\
               &=m(0+)Af+\int_{0}^{\infty}(f-e^{-tA}f)\,h_m(t)\,\dt\\
               &=\psi_{m}(A)f.
   \end{align*}
   According to Phillips' subordination theorem (Theorem~\ref{thm:Phillips}),	
   the domain $\mathcal{D}(A)$ of $A$ is a core of
   $\mathcal{D}(\psi_m(A))$. Thus, we have thereby shown that the two operators
   $\Lambda_m$ and $\psi_m(A)$ coincide.

   Further, let $\{\tilde{L}^{-1}_{t}\}_{t\ge 0}$ be the local inverse
   time at zero of the generalized diffusion process
   $\{Z_{t}\}_{t\ge 0}$ associated with $m$. Then by Knight's theorem
   (see Theorem~\ref{thm:2}), the Laplace transform determines
   uniquely (see~\eqref{eq:134}) that the convolution semigroup
   $\{\gamma_{t}\}_{t\ge 0}$ of sub-probability measures $\gamma_{t}$ on
   $[0,\infty)$ associated with $\psi_m$ has to be given by the
   push-forward measure \eqref{eq:135}.
 %
 Thus and by~\eqref{eq:51} from
 Proposition \ref{prop:semigroup-property-sub-probability}, one sees that
   %
   \begin{displaymath}
     e^{-t\psi_{m}(A)}f
     =\int_{[0,\infty)}e^{-sA}f\,\td \gamma_{t}(s)=\mathbb{E}\left (e^{-\tilde{L}^{-1}_{t}A}f \right )
   \end{displaymath}
   for every $t\ge 0$, $f\in X$, showing that~\eqref{eq:105}
   holds. This completes the proof of this theorem.
 \end{proof}

\subsection{Stability: Proof of Theorem~\ref{prop:limits_str}}
\label{subsec:stability}

This section is dedicated to outline the proof of
Theorem~\ref{prop:limits_str}.\medskip 

We begin this subsection with the following lemma.

\begin{lemma}
  \label{lem:convergence-of-measures}
  Let $m\in \mathfrak{m}_{\infty}$ and for every $n\in \N$,
  $m_{n}\in \mathfrak{m}_{\infty}$ be strings on $\R$ of infinite
  length respectively with associated measures $\mu_{m}$ and
  $\mu_{m_{n}}$. If $m_n(z)\to m(z)$ as $n\to \infty$ for every
  continuity point $z\in \R$ of $m$, then 
  one has that
    \begin{equation}
      \label{eq:147}
      \lim_{n\to\infty}L^{-1}_{m_{n},t}=L^{-1}_{m,t}\qquad\text{$\mathds{P}$-a.e. and
      for all $t\ge 0$.}
  \end{equation}
  where denote by $\{L^{-1}_{m_{n},t}\}_{t\ge 0}$ and
  $\{L^{-1}_{m,t}\}_{t\ge 0}$ the inverse local
  time process at zero of the generalized diffusions
  $\{Z_{m_{n},t}\}_{t\ge 0}$ and $\{Z_{m,t}\}_{t\ge 0}$ associated
  with $m_{n}$ and $m$.
\end{lemma}

For the proof of this lemma, we use the Portmanteau theorem applied to
the special case that the sequence $\{\mu_{m_{n}}\}_{n\ge 1}$ of
measures $\mu_{m_{n}}$ is induced by a string
$m_{n}\in \mathfrak{m}_{\infty}$.

\begin{theorem}[{\cite[Theorem 21.15]{MR3644418}}]\label{PROP:equiv_str}
  Let $m$ and $m_{n}$ for every $n\in \N$, be monotone increasing functions on
  $\R$ and $\mu_{m}$ and $\mu_{m_{n}}$ the associated measures to $m$
  and $m_{n}$, respectively. Then, the following statements are equivalent.
  \begin{enumerate}
  \item One has that $m_n(z)\to m(z)$ as $n\to \infty$ for every continuity point $z\in \R$ of $m$;
  \item One has that $\mu_{m_n}\to \mu_{m}$ as $n\to \infty$ vaguely.
  \end{enumerate}
\end{theorem}

With this result in mind, we can now outline the proof of
Lemma~\ref{lem:convergence-of-measures}.

\begin{proof}[Proof of Lemma~\ref{lem:convergence-of-measures}]
 Since each $m_{n}$ and $m$
  are strings on $\R$ of infinite length, \eqref{EQN:INV_LOC_TIME}
  yields that the associated inverse local
  time process $\{L^{-1}_{m_{n},t}\}_{t\ge 0}$ and $\{L^{-1}_{m,t}\}_{t\ge 0}$ at zero of the generalized diffusions
  $\{Z_{m_{n},t}\}_{t\ge 0}$ associated with $m_{n}$ and $\{Z_{m,t}\}_{t\ge 0}$ associated with $m$ admit the
  integral representation
 \begin{displaymath}
   L^{-1}_{m_{n},t}=\int_{[0,\infty)}L_{L_{t}^{-1}}(z)\,\td
   \mu_{m_{n}}(z)\quad\text{and}\quad 
L^{-1}_{m,t}=\int_{[0,\infty)}L_{L_{t}^{-1}}(z)\,\td
   \mu_{m}(z)
  \end{displaymath}
  $\mathds{P}$-a.e. and for all $t\ge 0$. Moreover, according to
  Theorem~\ref{PROP:equiv_str}, the pointwise convergence of
  $m_n(z)\to m(z)$ as $n\to \infty$ at every continuity point
  $z\in \R$ of $m$ yields that $\mu_{m_n}\to \mu_{m}$ as
  $n\to \infty$ vaguely. Thus, in order to establish the desired
  limit~\eqref{eq:147}, it is sufficient to show that
  \begin{equation}\label{EQN:local_vague}
    \lim_{n\to\infty}\int_{[0,\infty)}L_{t}(z)\,\td \mu_{m_{n}}(z)=
    \int_{[0,\infty)}L_{t}(z)\,\td \mu_{m}(z) 
  \end{equation} 
  $\mathds{P}$-a.s. for every $t\ge 0$. Note, the limit
  \eqref{EQN:local_vague} follows from the vaguely convergence
  $\mu_{m_n}\to \mu_{m}$ as $n\to \infty$ provided that for every
  $t\ge 0$, the map $L_{t}$ is continuous and has a compact support in
  $[0,\infty)$. To see that the latter holds, recall that the local
  time process $\{L_{t}(z)\}_{t\ge 0}$
  of 
  $\{B_{t}^{+}\}_{t\ge 0}$ is a jointly continuous mapping
  $L : [0,\infty)\times [0,\infty)\to [0,\infty)$ assigning
  $(t,z)\mapsto L_{t}(z)$. Thus, for every $t\ge 0$,
  $L_{t} : [0,\infty)\to [0,\infty)$ is continuous. To see that for
  every $t\ge 0$, the map $L_{t}$ has a compact support in
  $[0,\infty)$, set $M:=\sup_{s\in [0,t]}B^{+}_{s}$. By the continuity
  of the reflecting Brownian motion $\{B^{+}_{t}\}_{t\ge 0}$, the
  upper bound $M$ is finite and so, by applying the occupation times
  formula~\eqref{EQN:BROWNIAN_OCC_TIMES} to the function
  $g:=\mathds{1}_{(M,\infty)}$, we can conclude that
  \begin{displaymath}
    \int_{M}^{\infty}L_{t}(z)\,\td z=\int_{0}^{t}\mathds{1}_{(M,\infty)}(B^{+}_{s})\,\td s=0.
  \end{displaymath}
  This shows that for every $t\ge 0$, $L_{t}(z)=0$ for every $z>M$,
  that is, the function $z\mapsto L_{t}(z)$ has compact support in
  $[0,\infty)$. This proves~\eqref{EQN:local_vague}. 
  This completes the proof of this lemma.
\end{proof}

Now, we are ready to give the proof of Theorem~\ref{prop:limits_str}.

\begin{proof}[Proof of Theorem~\ref{prop:limits_str}]
  Let $A$ be an $m$-accretive operator on $X$, and for given strings
  $\{m_n\}_{n\ge 1}\subseteq \mathfrak{m}_{\infty}$ and
  $m \in \mathfrak{m}_{\infty}$, let $\psi_{m_{n}}$ be the associated
  complete Bernstein function given by~\eqref{eq:48}.
  Further, let $\{e^{-t\psi_{m_{n}}(A)}\}_{t\ge 0}$ and
  $\{e^{-t\psi_{m}(A)}\}_{t\ge 0}$ be the semigroups generated by
  $-\psi_{m_{n}}(A)$ and $-\psi_{m}(A)$ on $X$, respectively. Then, our first aim is to
  show that for every $f\in X$,
  \begin{equation}
    \label{eq:141}
   e^{-t\psi_{m_{n}}(A)}f\to e^{-t\psi_{m}(A)}f\quad\text{ in
     $X$ pointwise for every $t\ge 0$.}
 \end{equation}
 According to Theorem~\ref{thm:1}, the operator $\psi_{m_{n}}(A)$ and $\psi_{m}(A)$ 
 respectively coincide with their Dirichlet-to-Wenzell operator
 $\Lambda_{m_{n}}$ and $\Lambda_{m}$ given by~\eqref{eq:91}. In particular, by the
 integral representation~\eqref{eq:105} of Theorem \ref{thm:1}, one
 has that
  can be rewritten as
  \begin{displaymath}
e^{-t\psi_{m_{n}}(A)}f 
    =\mathds{E}\left(e^{-\tilde{L}^{-1}_{m_{n},t}\! A}f\right) 
    \quad\text{ and }\quad
    e^{-t\psi_{m}(A)}f 
    =\mathds{E}\left(e^{-\tilde{L}^{-1}_{m,t}\! A}f\right),
  \end{displaymath}
  where $\{\tilde{L}^{-1}_{m_{n},t}\}_{t\ge 0}$ and
  $\{\tilde{L}^{-1}_{m,t}\}_{t\ge 0}$ are respectively the local inverse times at
  zero of the generalized diffusions $\{Z_{m_{n},t}\}_{t\ge 0}$
  associated with $m_{n}$ and $\{Z_{m,t}\}_{t\ge 0}$
  associated with $m$. From this, one sees that
 \begin{displaymath}
    \norm{e^{-\psi_{m}(A)t}f-e^{-\psi_{m_{n}}(A)t}f}_{X}
      \le\mathbb{E}\lnorm{e^{-\tilde{L}^{-1}_{m,t}A}f-e^{-\tilde{L}^{-1}_{m_{n},t}A}f}_{X}.
 \end{displaymath}
 Since the semigroup $\{e^{-tA}\}_{t\ge 0}$ is contractive and
 strongly continuous$e^{-tA}$, limit~\eqref{eq:141} follows from Lebesgue's
 dominate convergence theorem and by limit~\eqref{eq:147}.

Now, for every $\lambda>0$, let
$R(\lambda,\tilde{\psi}_{m}(A)):=(\lambda\textrm{id}_{X}+\tilde{\psi}_{m}(A))^{-1}$
 be the
\emph{resolvent} operator of $\tilde{\psi}_{m}(A)$ on $X$ and for every $n\ge
1$, let $R(\lambda,\psi_{m_{n}}(A))$ be the resolvent operator of
$\psi_{m_{n}}(A)$. Then by~\cite[Theorem~1.10 in Chapter
II.]{MR2229872}, for every $f\in X$, $R(\lambda, \tilde{\psi}_{m}(A))f$ and
$R(\lambda,\psi_{m_{n}}(A))f$ admit the integral representations
  \begin{align*}
R(\lambda,\tilde{\psi}_{m}(A))f&=\int_{0}^{\infty}e^{-\lambda t}e^{-t
                         \tilde{\psi}_{m}(A)}f\,\dt
                 \intertext{ and }
                         R(\lambda,\psi_{m_{n}}(A))f&=\int_{0}^{\infty}e^{-\lambda
t}e^{-t\psi_{m_{n}}(A)}f\,\dt
  \end{align*} from where one can conclude that
  \begin{align*}
    &\norm{R(\lambda,\tilde{\psi}_{m}(A))f-R(\lambda,\psi_{m_{n}}(A))f}_{X}\\
    &\hspace{2cm}\le
\int_{0}^{\infty}\! e^{-\lambda t}\lnorm{e^{-t
\tilde{\psi}_{m}(A)}f-e^{-t\psi_{m_{n}}(A)}f}_{X}\dt.
\end{align*}
Thus, by~\eqref{eq:141}, and since the semigroups
$\{e^{-t\tilde{\psi}_{m}(A)}\}_{t\ge 0}$ and $\{e^{-t\psi_{m_{n}}(A)}\}_{t\ge
0}$ are contractive, Lebesgue's dominate
convergence theorem yields that 
\begin{displaymath}
  \lim_{n\to
    \infty}R(\lambda,\psi_{m_{n}}(A))f=R(\lambda,\tilde{\psi}_{m}(A))f
  \qquad\text{in $X$,}
  \end{displaymath}
  for every $\lambda>0$ and $f\in X$; that is, $\psi_{m_{n}}(A)\to \tilde{\psi}_{m}(A)$
  strongly in the \emph{resolvent sense}. By Trotter-Kato's first
  approximation theorem (see Remark~\ref{rem:Trotter-Kato}), this type
  of convergence is equivalent to the type of convergence stated in
  Theorem~\ref{prop:limits_str}.
\end{proof}



\section{Application}
\label{section:application}


In order to demonstrate the usefulness of the main results of this
paper, we discuss in this section one classical example.


\subsection{Limits of fractional powers $A^{\sigma}$ as $\sigma\to
  1-$}
\label{subsec:limits}

Let $A$ be an $m$-accretive operator on Banach space $X$. Then, we
intend to give a proof of the limit
\begin{equation}\label{eq:149}
	\lim_{\sigma\to 1-}
        A^\sigma =A \qquad \text{in the graph sense}
 \end{equation}
by using Theorem~\ref{prop:limits_str}. It is worth noting that this limit
has been well studied (see, for instance, \cite[Lemma 2.3]{MR115096}). Here, we provide an
alternative, which from our perspective is much simpler.\medskip

For given $0<\sigma<1$ and $f\in D(A)$, we begin by realizing the
fractional power $A^{\sigma}$ by the Dirichlet-to-Neumann operator
$\Lambda_{\tilde{m}_{\sigma}}$ associated with the following
incomplete Dirichlet problem
\begin{equation}\label{eq:117}
\begin{cases}
Au(z)-\left(\frac{\sigma}{1-\sigma}\right)\,z^\frac{2\sigma-1}{\sigma}\frac{\td^2
  u}{\td z^2}(z)&=0\qquad \text{for 
	$z\in (0,\infty)$,}\\
\hspace{4.15cm}u(0)&=f,
\end{cases}
\end{equation}
for the extension operator
\begin{displaymath}
  \mathcal{A}_{\tilde{m}_{\sigma}}=A+B_{\tilde{m}_{\sigma}}\qquad\text{
    with }\qquad
B_{\tilde{m}_{\sigma}}=-\Big(\frac{\sigma}{1-\sigma}\Big) z^\frac{2\sigma-1}{\sigma}\frac{\td^2 }{\td z^2}.
\end{displaymath}
By comparing Dirichlet problem~\eqref{eq:117} with~\eqref{eq:92} from
the introduction, then one realized that the extension equation
\begin{displaymath}
  \mathcal{A}_{\tilde{m}_{\sigma}}u=0\qquad\text{on $\mathcal{X}_{+}$}
\end{displaymath}
in problem~\eqref{eq:117} is a \emph{scaled version} (by the factor
$\sigma/(1-\sigma)$) of the extension equation
$\mathcal{A}_{m_{\sigma}}=0$ employed in Example~\ref{ex:1}. In fact,
this is one possibility to circumvent the fact that for the string
$m_{\sigma}$ introduced in~\eqref{eq:113} the limit as $\sigma\to 1-$
does not exist. The string
$\tilde{m}_{\sigma}$ corresponding to Dirichlet problem~\eqref{eq:117} is given by
\begin{equation}
  \label{eq:151}
\tilde{m}_{\sigma}(z)= \frac{\sigma}{1-\sigma} \, m_{\sigma}(z)=
\begin{cases}
\tfrac{1}{2}z^{\frac{1-\sigma}{\sigma}}& \qquad \text{if $z\ge 0$,}\\
0 & \qquad\text{if $z<0$,}
\end{cases}
\end{equation}
for every $z\in \R$, and the associate complete Bernstein function
\begin{displaymath}
  \psi_{\tilde{m}_{\sigma}}(\lambda)=\frac{\sigma^{\sigma-1}(1-\sigma)^\sigma}{2}
  \frac{\Gamma(1-\sigma)}{\Gamma(\sigma)}\,\lambda^{\sigma}\qquad\text{for every $\lambda\ge 0$.}
\end{displaymath}
Then, by our main Theorem~\ref{thm:1}, the
Dirichlet-to-Neumann operator $\Lambda_{\tilde{m}_\sigma}$ associated
with $\mathcal{A}_{\tilde{m}_{\sigma}}$ characterizes $A^{\sigma}$ up
to a multiplicative constant; namely, one has that
\begin{displaymath}
  \Lambda_{\tilde{m}_\sigma}=\frac{\sigma^{\sigma-1}(1-\sigma)^\sigma}{2}\frac{\Gamma(1-\sigma)}{\Gamma(\sigma)}A^{\sigma}
\end{displaymath}
(cf. the results mentioned in Section
\ref{subsec:represenation-formula}). Since for the family
$\{\tilde{m}_{\sigma}\}_{\sigma\in (0,1)}$ of strings, 
\begin{displaymath}
  \lim_{\sigma\to 1-}\tilde{m}_{\sigma}(z)=m(z):=
\begin{cases}
\frac{1}{2} & \text{if $z\ge 0$,}\\
0 & \text{if $z<0$}
\end{cases}
\end{displaymath}
the \emph{Heaviside step
  function}, by Theorem~\ref{prop:limits_str} and
Example~\ref{EX:dirac}, we get that
\begin{displaymath}
  \lim_{\sigma\to 1-}\frac{\sigma^{\sigma-1}(1-\sigma)^\sigma}{2}\frac{\Gamma(1-\sigma)}{\Gamma(\sigma)}A^\sigma
  =\frac{1}{2}A \qquad \text{in the graph sense.}
\end{displaymath}
From this and since
\begin{displaymath}
	\lim_{\sigma \rightarrow 1^{-}}\sigma^{\sigma-1}(1-\sigma)^\sigma\frac{\Gamma(1-\sigma)}{\Gamma(\sigma)}=1,
  \end{displaymath}
  we can conclude that limit~\eqref{eq:149} holds.

\begin{remark}[{Probabilistic justification}]
  Since $A^{\sigma}$ coincides up to a scalar multiple with the
  Dirichlet-to-Neumann operator associated with the extension operator
  $A+B_{m_{\sigma}}$ on $\mathcal{X}_{+}=X\times (0,\infty)$ and
  $B_{m_{\sigma}}$ is the generator of the
  $2\sigma^{\textrm{th}}$-powered process
  $\{Y_{t}^{2\sigma}\}_{t\ge 0}$ in $[0,\infty)$ of the Bessel process
  $\{Y_{t}\}_{t\ge 0}$ of dimension $2(1-\sigma)$, the
  limit~\eqref{eq:149} can be justified with probabilistic
  arguments. In fact, the limit~\eqref{eq:149} means that the process
  induced by $A^{\sigma}$ changes to the process generated by
  $A$. This makes sense if one considers the interaction of the Bessel
  process with the left boundary point $z=0$. For $\sigma \in (0,1)$,
  the set $\{0\}$ is reflecting for the Bessel process of dimension
  $2(1-\sigma)$. For the Bessel process of dimension $0$, the case
  when $\sigma=1$, the set $\{0\}$ is absorbing.
\end{remark}

To conclude this section we note that both strings $m_{\sigma}$
from~\eqref{eq:113} and $\tilde{m}_{\sigma}$ from~\eqref{eq:151} can't
be used in combination with Theorem~\ref{prop:limits_str} for proving
that
\begin{displaymath}
	\lim_{\sigma\to 0+}
        A^\sigma =\textrm{id}_{X} \qquad \text{in the graph sense}
 \end{displaymath}
holds for any given $m$-accretive operator $A$ on Banach space
$X$. The reason for this is that the pointwise limit
\begin{displaymath}
  \lim_{\sigma\to 0+}m_{\sigma}(z)=\lim_{\sigma\to 0+}\tilde{m}_{\sigma}(z)=m(z):=
\begin{cases}
\infty& \text{if $z>1$,}\\
0 & \text{if $z\le 1$,}
\end{cases}
\end{displaymath}
the \emph{Delta-function}, is not a string.
\appendix

%
%

\section{A Primer on Bessel processes}
\label{subsec:bessel-processes}
In this subsection, we give a brief reminder on the theory of Bessel
processes.  For a detailed review of this subject, we refer the
interested reader to~\cite[Chapter IX]{MR1725357} or \cite[Chapter
3]{MR3616274}. Throughout this paper, we denote by
$\{B_{t}\}_{t\ge 0}$ a standard \emph{Brownian motion} in $\R$ starting at
$x=0$.\medskip

We begin by recalling the definition of the \emph{squared Bessel process of
dimension $\delta\ge 0$}; this process $\{Y^{2}_{t}\}_{t\geq 0}$ is the unique, continuous,
non-negative process defined by the stochastic differential equation (SDE)
\begin{displaymath}
	Y^2_{t}=y^2+\delta t+\int_{0}^{t}2\sqrt{Y^2_{s}}\,\td B_{s}
        \qquad\text{for every $t\ge 0$, $y\in \R$.}
\end{displaymath}
Then, for $\delta\ge 0$, the \emph{$\delta$-Bessel process}
$\{Y_{t}\}_{t\ge 0}$ is obtained by taking the square root
\begin{equation}
  \label{eq:2}
  Y_{t}:=\sqrt{Y^2_{t}}\qquad \text{for every $t\ge0$.}
\end{equation}

In the case $\delta>1$, It\^{o}'s Lemma yields that the $\delta$-Bessel process
$\{Y_{t}\}_{t\ge 0}$ is the unique solution to the SDE
\begin{equation}\label{eq:3}
	Y_{t}=y+\frac{\delta-1}{2}\int_{0}^{t}\frac{1}{Y_{s}}\,\ds+B_{t} 
        \qquad\text{for every $t\ge 0$, $y\ge 0$.}
 \end{equation}
It is well known that equation~\eqref{eq:3} satisfies path-wise uniqueness
and has strong solutions due to the drift being monotone decreasing
(cf.,~\cite[Theorem 3.2]{MR3616274}).
For $\delta=1$, the squared Bessel process $\{Y_{t}\}_{t\geq 0}$ 
coincides with the squared Brownian motion and hence, the $1$-Bessel
process $\{Y_{t}\}_{t\ge 0}$ coincides with the \emph{reflecting
  Brownian motion} $\{B^{+}_{t}\}_{t\ge 0}$, which is
the unique solution of the SDE 
\begin{equation}
  \label{eq:4}
   B^{+}_{t}=y+B_{t}+L_{t} \qquad \text{for every $t\ge 0$, $y\ge 0$,}
\end{equation}
where $\{L_{t}\}_{t\ge 0}$ is a continuous, monotone, non-decreasing
process, with $L_{0}=0$, $\td L_{t}\ge 0$, and
$\int_{0}^{t}B^{+}_{s}\,\td L_{s}=0$ for all $t\ge 0$. Uniqueness of the
solution $\{B^{+}_{t}\}_{t\ge 0}$ to~\eqref{eq:4} follows from
the Skorokhod lemma (see~\cite[Lemma~2.1]{MR3616274}).

In the case $0<\delta<1$, the situation is far more delicate and the
$\delta$-Bessel process $\{Y_{t}\}_{t\ge 0}$ is no longer a
solution of equation~\eqref{eq:3}. Instead
(cf.,~\cite[Proposition 3.8 \& 3.12]{MR3616274}), for the process
$\{Y_{t}\}_{t\ge 0}$ given by~\eqref{eq:2}, there is a continuous
family $\{\ell_{t}^{a}\}_{a, t\ge 0}$, called \emph{diffusion local
 times}, satisfying the \emph{occupation times formula}
\begin{displaymath}
	\int_{0}^{t}\varphi(Y_{r})\,\dr=\int_{0}^{\infty}\varphi(a)\, \ell_{t}^{a}\,a^{\delta-1}\,\td a
\end{displaymath} 
for all $t\ge 0$ and bounded and Borel-measurable functions
$\varphi: \mathbb{R}_{+}\rightarrow \mathbb{R}_{+}$, and
$\{Y_{t}\}_{t\ge 0}$ satisfies
\begin{displaymath}
	Y_{t}=y+\frac{\delta-1}{2(2-\delta)}\int_{0}^{\infty}\frac{\ell_{t}^{a}-\ell_{t}^{0}}{a}a^{\delta-1}\,\td
       a+B_{t}\qquad\text{for all $t\ge 0$, $y\ge 0$.}
\end{displaymath}

The next result due to Donati-Martin, Roynette, Vallois and Yor~\cite
{MR2417969} is quite important for understanding the probabilistic
approach to the Dirichlet-to-Neumann operator in the fractional power case
as discussed in Section~\ref{sec:motiviation}.

\begin{theorem}\label{thm:Donati-yor} 
 For $\sigma \in (0,1)$, let $\delta=2(1-\sigma)$ and $\{Y_{t}\}_{t\ge 0}$ be
 the $\delta$-Bessel process starting at $y=0$. Then, the
 $2\sigma^{\textrm{th}}$-powered process $\{Y_{t}^{2\sigma}\}_{t\ge 0}$ is a
 submartingale with the Doob-Meyer decomposition
\begin{equation*}
 Y_{t}^{2\sigma}=2\sigma\int_{0}^{t}Y_{r}^{2\sigma-1}\td B_r+L_{t}
\end{equation*}
where $\{L_{t}\}_{t\ge 0}$ is a continuous non-decreasing
process, carried by the zeros of $\{Y_{t}\}_{t\ge 0}$; that is, $\{L_{t}\}_{t\ge 0}$ only increases when
$Y_{t}=0$. Further, the inverse local time process $\{L_{t}^{-1}\}_{t\ge 0}$ of
$\{L_{t}\}_{t\ge 0}$ given by
\begin{displaymath}
 L_{t}^{-1}:=\inf\Big\{r>0\,\Big\vert\,L_{r}>t\Big\},\qquad(t\ge 0),
\end{displaymath}
is an $\sigma$-stable subordinator satisfying
\begin{displaymath}
 \mathbb{E}\Big(e^{-u
 L_{t}^{-1}}\Big)=e^{-t\frac{\Gamma(1-\sigma)}{\Gamma(1+\sigma)}\frac{u^\sigma}{2^\sigma}}\qquad\text{for
 every $t$, $u>0$}
\end{displaymath}
and L\'{e}vy measure $\nu$ given by~\eqref{eq:18}.
\end{theorem}

%
%
%
%
\section{A distributional definition of $\frac{\td f}{\td m}$}
\label{sec:def-weak}

The first section of the appendix deals with some reminder on
\emph{distributional definition} of $\frac{\td f}{\td m}$ when
$m\in \mathfrak{m}_{\infty}$ is a string and $f : [0,\infty)\to X$ is a
Banach space $X$-valued measurable function. Throughout this section,
let $X$ be a Banach space with norm $\norm{\cdot}_{X}$ and $X'$ its
dual space with duality brackets
$\langle\cdot,\cdot\rangle_{X',X}$ and $m\in \mathfrak{m}_{\infty}$ a
string with associated measure $\mu_{m}$.\medskip

We begin by considering the \emph{smooth} case. 

\begin{definition}
  We say that a function $f : [0,\infty)\to X$ is \emph{$m$-differentiable} at
  $t\in (0,\infty)$ provide the limit
  \begin{displaymath}
    \frac{\td f}{\td m}(t):=\lim_{h\to
    0}\frac{f(t+h)-f(t)}{m(t+h)-m(t)}\qquad\text{exists in $X$.}
  \end{displaymath}
  Then, we call $\frac{\td f}{\td m}(t)$ the \emph{(classic)
    $m$-derivative} of $f$ at $t\in (0,\infty)$.
\end{definition}

Of course, one natural question is to ask, which functions are
$m$-differentiable. One example is given by the following result
generalizing the classic theorem of Lebesgue (see, e.g.,~\cite[Theorem~1.18]{MR3726909}).

\begin{theorem}
  \label{thm:monotone}
  Let $m\in \mathfrak{m}_{\infty}$ be a string with associated measure
  $\mu_{m}$. Then every monotone function $f :
  [0,\infty)\to \R$ is $\mu_{m}$-a.e. $m$-differentiable.
\end{theorem}

\begin{proof}
  By Lebesgue's theorem, the two sets 
  \begin{displaymath}
    N_{f}:=\Big\{t\in [0,\infty)\,\Big\vert\, f'(t)\text{ does not
      exist }\Big\}
  \end{displaymath}
  and
  \begin{displaymath}
    N_{m}:=\Big\{t\in [0,\infty)\,\Big\vert\, m'(t)\text{ does not
      exist }\Big\}
  \end{displaymath}
  are measurable subsets of $[0,\infty)$ with Lebesgue measure
  $\lambda(N_{f})=0$ and $\lambda(N_{m})=0$. 
  Further, we decompose $[0,r_{m})$ into three disjoint Borel-measurable sets
  \begin{displaymath}
    A_{1}=\Big\{t\in N_{m}^{c}\, \Big\vert\,m'(t)=0 \Big\},\quad A_{2}=\Big\{t\in N_{m}^{c}\,\vert\,m'(t)\neq
    0 \Big\},\quad
    A_{3}=N_{m}.
  \end{displaymath}
  By the partition $(A_{i})_{i=1}^{3}$ of $[0,\infty)$, one has
  that 
  \begin{displaymath}
    [0,\infty)=N_{f}\,\dot{\cup}\,N_{f}^{c}=\dot\bigcup_{i=1}^{3}(N_{f}^{c}\cap
    A_{i})\,\dot\cup\, \dot\bigcup_{i=1}^{3}(N_{f}\cap A_{i}).
  \end{displaymath}
  Since the support $E_{m}:=\textrm{supp}(\mu_{m})$ of the
  measure $\mu_{m}$ is the set where $m$ increases, one has that
  $\mu_{m}(A_{1})=0$ and $\mu_{m}(A_{3})=0$. Thus, by the monotonicity of
  $\mu_{m}$, we can conclude that $\mu_{m}(N_{f}\cap A_{1})=0$,
  $\mu_{m}(N_{f}^{c}\cap A_{1})=0$, $\mu_{m}(N_{f}\cap A_{3})=0$, and $\mu_{m}(N^{c}_{f}\cap A_{3})=0$.
  Therefore, it remains to focus on the two cases
  \begin{displaymath}
    N_{f}\cap A_{2}=\Big\{t\, \Big\vert\,m'(t)\neq 0\;\&\;f'(t)\text{ does not exist}\Big\}
  \end{displaymath}
  and
  \begin{displaymath}
    N_{f}^{c}\cap A_{2}=\Big\{t\, \Big\vert\,m'(t)\neq 0\;\&\;f'(t)\text{ does exist}\Big\}.
  \end{displaymath}
  By Lebesgue 's decomposition (cf.,~\cite[Theorem~B.67 \& Section 6.3]{MR3726909}),
  the Lebesgue-Stieltje measure $\mu_{m}$ has the unique decomposition
  \begin{displaymath}
    \mu_{m}=\mu_{m,ac}+\mu_{m,s},
  \end{displaymath}
  where $\mu_{m,ac}$ is absolutely continuous w.r.t. the Lebesgue
  measure $\lambda$ and $\mu_{m,s}$ is singular. Moreover,
  $\mu_{m,ac}$ is given by
  \begin{displaymath}
    \mu_{m,ac}(E)=\int_{E\cap N_{m}^{c}} m'(t)\,\td\lambda(t)
  \end{displaymath}
  for every Borel-measurable subset $E$ of $[0,\infty)$. But since
  $N_{f}$ has Lebesgue measure $\lambda(N_{f})=0$, we also have that
  \begin{displaymath}
    \mu_{m}(N_{f}\cap A_{2})=\mu_{m,ac}(N_{f})=0.
  \end{displaymath}
  Therefore, we have shown that the set 
  \begin{displaymath}
    N:=[0,\infty) \setminus (N_{f}\cap A_{2})
  \end{displaymath}
  has measure $\mu_{m}(N)=0$ and
  \begin{displaymath}
    \frac{\td f}{\td m}(t)\text{ exits at all
  $t\in [0,\infty)\setminus N=\Big\{t\in [0,\infty)\, \Big\vert\,m'(t)\neq 0\;\&\;f'(t)\text{ exists}\Big\}$.}
\end{displaymath}
\end{proof}

For the moment, suppose the
string $m\in \mathfrak{m}_{\infty}$ is a smooth real-valued function satisfying
$m'\ge c_{0}$ on $\R_{+}$ for some $c_{0}>0$. Then, the associated measure
$\mu_{m}$ to $m$ is absolutely continuous with respect to
the Lebesgue measure and if $f$ is differentiable at $t$, then one
easily sees that $f$ is, in particular, $m$-differentiable and 
\begin{displaymath}
  \frac{\td f}{\td m}(t)=\lim_{h\to
    0}\frac{f(t+h)-f(t)}{h}\frac{h}{m(t+h)-m(t)}=f'(t)\frac{1}{m'(t)}.
\end{displaymath}
Now, set $g(t):=\frac{\td f}{\td m}(t)$ and recall that Radon-Nikodym
derivative $\frac{\td\mu_{m}}{\dt}(t)=m'(t)$ on $(0,\infty)$. Therefore, an
integration by parts shows that
\begin{align*}
  \int_{0}^{\infty}g(t)\,\xi(t)\td\mu_{m}(t)
     &=\int_{0}^{\infty} f'(t)\frac{1}{m'(t)}\,\xi(t)\td\mu_{m}(t)\\
     &=\int_{0}^{\infty} f'(t)\,\xi(t)\dt\\
     &=-\int_{0}^{\infty} f(t)\,\xi'(t)\dt
\end{align*}
for every real-valued test function $\xi\in C^{\infty}_{c}((0,\infty))$. We
emphasize that since the measure $\mu_{m}$ associated to the given string
$m$ is a Radon measure on $[0,\infty)$, we have the following result.

\begin{lemma}
  \label{lem:variational-lemma}
  Let $m\in \mathfrak{m}_{\infty}$ be a string with associated measure
  $\mu_{m}$. If a function $f\in
  L^{1}_{loc,\mu_{m}}(0,\infty;X)$ satisfies
  \begin{equation}
    \label{eq:97}
    \int_{0}^{\infty}f(t)\,\xi(t)\,\td\mu_{m}(t)=0
  \end{equation}
  for all $\xi \in C^{\infty}_{c}((0,\infty))$, then $f(t)=0$ in $X$
  for a.e. $t\in (0,\infty)$.
\end{lemma}

\begin{proof}
  Let $\xi \in C^{\infty}_{c}((0,\infty))$ and suppose~\eqref{eq:97} holds.
  Then multiplying~\eqref{eq:97} by $x'\in X$ gives
  \begin{displaymath}
    \int_{0}^{\infty}\langle x',f(t)\rangle_{X',X}\,\xi(t)\,\td\mu_{m}(t)=0
  \end{displaymath}
  Since the measure $\mu_{m}$ associated to the given string $m$ is a
  Radon measure on $[0,\infty)$, and the scalar-valued function
  $t\mapsto \langle x',f(t)\rangle_{X',X}$ belongs to
  $L^{1}_{loc,\mu_{m}}(0,\infty)$, it follows
  from~\cite[Theorem~3.14]{MR924157} that
  $\langle x',f(t)\rangle_{X',X}=0$ for a.e. $t\in (0,\infty)$. Since
  $x'\in X'$ was arbitrary, and the dual space $X$ is separating
  elements of $X$ (Hahn-Banach), we can conclude that $f(t)=0$ in $X$
  for a.e. $t\in (0,\infty)$.
\end{proof}

Because of Lemma~\ref{lem:variational-lemma}, Definition~\ref{def:weak-dfdm-intro} of the
\emph{weak $m$-derivative} make sense.

\begin{definition}
  \label{def:m-continuous}
  For a given string $m\in \mathfrak{m}_{\infty}$ on $\R$, a function
  $f : [0,\infty)\to \R$ is called to be \emph{$m$-continuous} at
  $t_{0}\in (0,\infty)$ if for every $\varepsilon>0$ there is an
  $\delta>0$ such that for all $t\in [0,\infty)$ satisfying
  $\abs{m(t_{0})-m(t)}<\delta$, one has
  $\norm{f(t)-f(t_{0})}_{X}<\varepsilon$. Further, $f$ is called
  \emph{uniformly $m$-continuous} on $[0,\infty)$ if  for every $\varepsilon>0$ there is an
  $\delta>0$ such that for all $t_{1}$, $t_{2}\in [0,\infty)$ satisfying
  $\abs{m(t_{1})-m(t_{2})}<\delta$, one has
  $\norm{f(t_{1})-f(t_{2})}_{X}<\varepsilon$. Finally, $f$ is called
  \emph{absolutely $m$-continuous} on $[0,\infty)$ if  for every $\varepsilon>0$ there is an
  $\delta>0$ such that for every finite family
  $((a_{k},b_{k}))_{k=1}^{n}$ of disjoint sub-intervals of $[0,\infty)$ satisfying
  $\sum_{k=1}^{n}\abs{m(b_{k})-m(a_{k})}<\delta$, one has
  $\sum_{k=1}^{n}\norm{f(b_{k})-f(a_{k})}_{X}<\varepsilon$. 
\end{definition}

\begin{proposition}\label{prop:m-absolutely-continuous}
   Let $m\in \mathfrak{m}_{\infty}$ be a string with associated measure
  $\mu_{m}$. If $g\in
  L^{1}_{loc,\mu_{m}}([0,\infty);X)$ then  
  \begin{displaymath}
    f(t):=\int_{t_{0}}^{t}g(r)\,\td\mu_{m}(r),\quad t\in [0,\infty),
  \end{displaymath}
  is locally $m$-absolutely continuous on $[0,\infty)$ and at
  $\mu_{m}$-a.e. $t\in [0,\infty)$, $f(t)$ is
  $m$-differentiable with $\tfrac{\td f}{\td m}(t)=g(t)$. 
\end{proposition}

\begin{proof}
  Since $f$ satisfies
  \begin{displaymath}
    \norm{f(t)-f(s)}_{X}\le \int_{s}^{t}\norm{g(r)}_{X}\,\td\mu_{m}(r)
  \end{displaymath}
  for every $t$, $s\in [0,\infty)$ with $s<t$, it follows that $f$ is
  \emph{$m$-absolutely continuous}. 
  By \cite[Corollary~2.23]{MR1857292}, one has that for every $\mu_{m}$-Lebesgue point $t\in
 [0,\infty)$,
 \begin{displaymath}
\lim_{h\to 0+} \tfrac{1}{m(t+h)-m(t-h)}\int_{t-h}^{t+h}\norm{g(r)-g(t)}_{X}\,\td\mu_{m}(r)=0.
\end{displaymath}
Further, for such $t\in [0,\infty)$ and $h>0$, one has that
\begin{align*}
  \lnorm{g(t)-\frac{f(t+h)-f(t)}{m(t+h)-m(t)}}_{X}
  &=\lnorm{g(t)-\tfrac{1}{m(t+h)-m(t)}\int_{t}^{t+h}g(r)\,\td\mu_{m}(r)}_{X}\\
  &\le \tfrac{1}{m(t+h)-m(t)}\int_{t}^{t+h}\lnorm{g(t)-g(r)}_{X}\,\td\mu_{m}(r).
\end{align*}
Similarly, for such $t\in [0,\infty)$ and $h<0$, 
\begin{displaymath}
  \lnorm{g(t)-\frac{f(t+h)-f(t)}{m(t+h)-m(t)}}_{X}\le 
  \tfrac{1}{m(t)-m(t+h)}\int_{t+h}^{t}\lnorm{g(t)-g(r)}_{X}\,\td\mu_{m}(r).
\end{displaymath}
Now, by letting $h\to 0+$ in the above estimates yields that $f$ is
$m$-differentiable with $\tfrac{\td f}{\td m}(t)=g(t)$.
\end{proof}

Further, we have the following important lemma.

\begin{lemma}
  \label{lem:variational-lemma-xi-dash}
  If a function $f\in
  L^{1}_{loc}(0,\infty;X)$ satisfies
  \begin{equation}
    \label{eq:97}
    \int_{0}^{\infty}f(t)\,\xi'(t)\,\td\mu_{m}(t)=0
  \end{equation}
  for all $\xi \in C^{\infty}_{c}((0,\infty))$, then there is a $C\in
  X$ such that $f(t)= C$ for a.e. $t\in (0,\infty)$. 
\end{lemma}

\begin{proof}
  Let $\xi$ and $\eta\in C^{\infty}_{c}((0,\infty))$ such that $\eta\neq
  0$. Since $\eta-\xi'\in C^{\infty}_{c}((0,\infty))$, the function 
  \begin{displaymath}
    \chi(t):=\frac{\eta(t)-\xi'(t)}{\int_{0}^{\infty}\eta(r)\,\td\mu_{m}(r)},\qquad
    t\in (0,\infty),
  \end{displaymath}
  belongs to $C^{\infty}_{c}((0,\infty))$ and
  $\int_{0}^{\infty}\chi(r)\,\td\mu_{m}(r)=1$. Inserting $\xi$
  into~\eqref{eq:97} gives
  \begin{align*}
    0&=\int_{0}^{\infty}f(t)\,\xi'(t)\,\td\mu_{m}(t)\\
     &=\int_{0}^{\infty}f(t)\,\eta(t)\,\td\mu_{m}(t)
       -\int_{0}^{\infty}\left(f(t)\,\chi(t)\,\int_{0}^{\infty}\eta(r)\,
       \td\mu_{m}(r)\right)\,\td\mu_{m}(t)\\
    &=\int_{0}^{\infty}\left(f(t)-\int_{0}^{\infty}\, f(r)\,\chi(r)
       \td\mu_{m}(r)\right) \,\eta(t)\,\td\mu_{m}(t)
  \end{align*}
  Since $\eta\in C^{\infty}_{c}((0,\infty))$ was arbitrary, it follows
  from Lemma~\ref{lem:variational-lemma} that
  \begin{displaymath}
    f(t)=C:=\int_{0}^{\infty}\, f(r)\,\chi(r)
       \td\mu_{m}(r)\qquad\text{for a.e. $t\in (0,\infty)$.}
  \end{displaymath}
\end{proof}

Thanks to the above lemma, we can make the following statements.

\begin{proposition}
  \label{propo:chara-weak-m-derivative}
  Let $m\in \mathfrak{m}_{\infty}$ be a string with associated measure
  $\mu_{m}$. Further, let $f\in L^{1}_{loc}([0,\infty);X)$ and $g\in
  L^{1}_{loc,\mu_{m}}([0,\infty);X)$. Then, the following statements
  are equivalent.
  \begin{enumerate}
  \item One has that $g=\tfrac{\td f}{\td m}$ is the weak $m$-derivative
    of $f$;
   \item There is an $x\in X$ such that 
     \begin{equation}
       \label{eq:99}
       f(t)=x+\int_{0}^{t}g(r)\,\td\mu_{m}(r),\quad \text{a.e. on $[0,\infty)$.}
     \end{equation}  
  \end{enumerate}
\end{proposition}

\begin{proof}
  We show that $(1)$ implies $(2)$. For this, set 
  \begin{displaymath}
    w(t)=f(t)-\int_{0}^{t}g(r)\,\td\mu_{m}(r),\quad t\in [0,\infty).
  \end{displaymath}
  Then, by Fubini's theorem,
  \begin{align*}
    \int_{0}^{\infty}\left(\int_{0}^{t}g(r)\,\td\mu_{m}(r)\right)\,\xi'(t)\,\dt
    &=  \int_{0}^{\infty}\left(\int_{r}^{\infty}\xi'(t)\,\dt \right)
    g(r)\, \td\mu_{m}(r)\\
    &=- \int_{0}^{\infty}g(r)\, \xi(r) \td\mu_{m}(r)
  \end{align*}
  and so,
  \begin{displaymath}
    \int_{0}^{\infty} w(t) \,\xi'(t)\,\dt
    =\int_{0}^{\infty} f(t) \,\xi'(t)\,\dt-
    \int_{0}^{\infty}\left(\int_{0}^{t}g(r)\,\td\mu_{m}(r)\right)\,\xi'(t)\,\dt=0.
  \end{displaymath}
  Therefore, \eqref{eq:99} follows from
  Lemma~\ref{lem:variational-lemma-xi-dash}. The implication $(2)$ implies
  $(1)$ follows from Proposition~\ref{prop:m-absolutely-continuous}.
\end{proof}

\end{document}